\documentclass[11pt]{amsart}
\usepackage{amsmath}
\usepackage{amsxtra}
\usepackage{amscd}
\usepackage{amsthm}
\usepackage{amsfonts}
\usepackage{amssymb}
\usepackage{eucal}

\newtheorem{cor}[subsubsection]{Corollary}
\newtheorem{lem}[subsubsection]{Lemma}
\newtheorem{prop}[subsubsection]{Proposition}
\newtheorem{propconstr}[subsubsection]{Proposition-Construction}

\newtheorem{conj}[subsubsection]{Conjecture}
\newtheorem{thm}[subsubsection]{Theorem}

\theoremstyle{remark}

\newcommand{\propconstrref}[1]{Proposition-Construction~\ref{#1}}
\newcommand{\thmref}[1]{Theorem~\ref{#1}}
\newcommand{\secref}[1]{Sect.~\ref{#1}}
\newcommand{\lemref}[1]{Lemma~\ref{#1}}
\newcommand{\propref}[1]{Proposition~\ref{#1}}
\newcommand{\corref}[1]{Corollary~\ref{#1}}
\newcommand{\conjref}[1]{Conjecture~\ref{#1}}

\newcommand{\nc}{\newcommand}
\newcommand{\renc}{\renewcommand}

\nc\on{\operatorname}
\nc\ssec{\subsection}
\nc\sssec{\subsubsection}

\nc\Ind{\on{Ind}}
\nc\Res{\on{Res}}
\nc\Loc{\on{Loc}}
\nc{\Gr}{\on{Gr}_G}
\nc{\Fl}{\on{Fl}_G}
\nc{\IC}{\on{IC}}

\nc{\bU}{\mathbf U}
\nc{\bM}{\mathbf M}
\nc{\bW}{\mathbf W}
\nc{\bV}{\mathbf V}
\nc{\bB}{\mathbf B}
\nc{\bL}{\mathbf L}
\nc{\bC}{\mathbf C}
\nc{\bg}{\mathbf g}
\nc{\bz}{\mathbf z}

\nc{\tu}{\overset{\bullet}{\mathfrak u}{}}
\nc{\tb}{\overset{\bullet}{\mathfrak b}{}}
\nc{\tRg}{\overset{\bullet}{\mathcal R}_{\check{G}}}
\nc{\Rg}{{\mathcal R}_{\check{G}}}
\nc{\tbM}{\overset{\bullet}{\bM}{}}
\nc{\tM}{\overset{\bullet}{M}{}}
\nc{\tR}{\overset{\bullet}{R}{}}
\nc{\tCM}{\overset{\bullet}{{\mathcal M}}{}}
\nc{\tL}{\overset{\bullet}{L}{}}
\nc{\tCF}{\overset{\bullet}{{\mathcal F}}{}}
\nc{\tCS}{\overset{\bullet}{{\mathcal S}}{}}

\nc{\fu}{\mathfrak u}
\nc{\fb}{\mathfrak b}
\nc{\fh}{\mathfrak h}
\nc{\fn}{\mathfrak n}
\nc{\fg}{\mathfrak g}
\nc{\fl}{\mathfrak l}
\nc{\fp}{\mathfrak p}
\nc{\fs}{\mathfrak s}
\nc{\fF}{\mathfrak F}
\nc{\fP}{\mathfrak P}

\nc{\sF}{\mathsf F}
\nc{\sG}{\mathsf G}
\nc{\sJ}{\mathsf J}

\nc{\BC}{\mathbb C}
\nc{\BD}{\mathbb D}
\nc{\BG}{\mathbb G}
\nc{\BO}{\mathbb O}
\nc{\BN}{\mathbb N}
\nc{\BZ}{\mathbb Z}

\nc{\CO}{\mathcal O}
\nc{\CA}{\mathcal A}
\nc{\CC}{\mathcal C}
\nc{\CP}{\mathcal P}
\nc{\CM}{\mathcal M}
\nc{\CW}{\mathcal W}
\nc{\CV}{\mathcal V}
\nc{\CU}{\mathcal U}
\nc{\CL}{\mathcal L}
\nc{\CB}{\mathcal B}
\nc{\CI}{\mathcal I}
\nc{\CJ}{\mathcal J}
\nc{\CF}{\mathcal F}
\nc{\CN}{\mathcal N}
\nc{\CZ}{\mathcal Z}
\nc{\CS}{\mathcal S}
\nc{\CH}{\mathcal H}
\nc{\CE}{\mathcal E}
\nc{\CY}{\mathcal Y}
\nc{\CK}{\mathcal K}
\nc{\CG}{\mathcal G}

\nc{\QCoh}{\on{QCoh}}
\nc{\Coh}{\on{Coh}}

\nc{\Catg}{\mathsf{Hecke}(\bU_\ell,\check G)}
\nc{\Catgd}{\overset{\bullet}{\mathsf{H}}\mathsf{ecke}(\bU_\ell,\check{G})}
\nc{\Catgeom}{{\mathsf{Hecke}(\Gr,\check G)}}
\nc{\Catgeomk}{{\mathsf{Hecke}(\Gr,\check G)^{G^k}}}
\nc{\Catgeomd}{\overset{\bullet}{\mathsf{H}}\mathsf{ecke}(\Gr,\check G)}
\nc{\Catgeomdi}{\overset{\bullet}{\mathsf{H}}\mathsf{ecke}(\Gr,\check G)^I}
\nc{\Catgeomni}{\mathsf{Hecke}(\Gr,\check G)^{I^0}}
\nc{\Catgeomi}{\mathsf{Hecke}(\Gr,\check G)^{I}}
\nc{\Catgeomdni}{\overset{\bullet}{\mathsf{H}}\mathsf{ecke}(\Gr,\check G)^{I^0}}
\nc{\Catgeomdk}{\overset{\bullet}{\mathsf{H}}\mathsf{ecke}(\Gr,\check G)^{G^k}}

\nc{\Perv}{\mathsf{Perv}}
\nc{\Pervgr}{\mathsf{Perv}(\Gr)}
\nc{\Pervgrk}{\mathsf{Perv}(\Gr)^{G^k}}
\nc{\Pervgrg}{\mathsf{Perv}(\Gr)^{G[[t]]}}
\nc{\Pervgri}{\mathsf{Perv}(\Gr)^I}
\nc{\Pervgrw}{\mathsf{Perv}(\Gr)^{I^-,\psi}}
\nc{\Dgrw}{\mathsf{D}(\Gr)^{I^-,\psi}}
\nc{\fPervgrni}{{}^f{\mathsf{Perv}}(\Gr)^{I^0}}
\nc{\fPervflni}{{}^f{\mathsf{Perv}}(\Fl)^{I^0}}
\nc{\fPervfli}{{}^f{\mathsf{Perv}}(\Fl)^{I}}
\nc{\Pervgrnib}{\overline{\mathsf{Perv}}(\Gr)^{I^0}}
\nc{\Pervgrib}{\overline{\mathsf{Perv}}(\Gr)^I}
\nc{\Pervgrb}{\overline{\mathsf{Perv}}(\Gr)}
\nc{\fPervgri}{{}^f{\mathsf{Perv}}(\Gr)^I}
\nc{\fDflni}{{}^f{\mathsf{D}}(\Fl)^{I^0}}
\nc{\Dgri}{\mathsf{D}(\Gr)^I}
\nc{\Dgrni}{\mathsf{D}(\Gr)^{I^0}}
\nc{\Dgr}{\mathsf{D}(\Gr)}
\nc{\Dgrg}{\mathsf{D}(\Gr)^{G[[t]]}}
\nc{\Dgrone}{\mathsf{D}(\Gr)^{G^1}}
\nc{\Dfli}{\mathsf{D}(\Fl)^I}
\nc{\Dfl}{\mathsf{D}(\Fl)}
\nc{\Dflni}{\mathsf{D}(\Fl)^{I^0}}
\nc{\Pervgrni}{\mathsf{Perv}(\Gr)^{I^0}}
\nc{\fDgrni}{{}^f\mathsf{D}(\Gr)^{I^0}}
\nc{\fDgri}{{}^f\mathsf{D}(\Gr)^{I}}
\nc{\fDfli}{{}^f\mathsf{D}(\Fl)^{I}}
\nc{\Pervfl}{\mathsf{Perv}(\Fl)}
\nc{\Pervfli}{\mathsf{Perv}(\Fl)^I}
\nc{\Pervflni}{\mathsf{Perv}(\Fl)^{I^0}}
\nc{\Pervfdb}{\mathsf{Perv}(G/B)^B}
\nc{\fPervfd}{{}^f\mathsf{Perv}(G/B)^B}
\nc{\fPervfdn}{{}^f\mathsf{Perv}(G/B)^N}
\nc{\Dfd}{\mathsf{D}(G/B)}
\nc{\Dfdn}{\mathsf{D}(G/B)^N}
\nc{\Pervfdn}{\mathsf{Perv}(G/B)^N}
\nc{\Pervfd}{\mathsf{Perv}(G/B)}
\nc{\Pervfdw}{\mathsf{Perv}(G/B)^{N^-,\psi}}
\nc{\Dfdw}{\mathsf{D}(G/B)^{N^-,\psi}}
\nc{\Pervfdnw}{\mathsf{Perv}(G/N)^{N^-,\psi}}
\nc{\Dfdnw}{\mathsf{D}(G/N)^{N^-,\psi}}
\nc{\Sph}{{\mathsf{Sph}_G}}

\nc{\Pervsm}{\mathsf{Perv}\bigl({\mathcal Fl}^{\frac{\infty}{2}}\bigr){}}
\nc{\Pervsmk}{\mathsf{Perv}\bigl({\mathcal Fl}^{\frac{\infty}{2}}\bigr){}^{G^k}}
\nc{\Pervsmkone}{\mathsf{Perv}\bigl({\mathcal Fl}^{\frac{\infty}{2}}\bigr){}^{G^{k_1}}}
\nc{\Pervsmktwo}{\mathsf{Perv}\bigl({\mathcal Fl}^{\frac{\infty}{2}}\bigr){}^{G^{k_2}}}
\nc{\Pervsmg}{\mathsf{Perv}\bigl({\mathcal Fl}^{\frac{\infty}{2}}\bigr){}^{G[[t]]}}
\nc{\Pervsmi}{\mathsf{Perv}\bigl({\mathcal Fl}^{\frac{\infty}{2}}\bigr){}^I}
\nc{\Pervsmni}{\mathsf{Perv}\bigl({\mathcal Fl}^{\frac{\infty}{2}}\bigr){}^{I^0}}
\nc{\Pervsmw}{\mathsf{Perv}\bigl({\mathcal Fl}^{\frac{\infty}{2}}\bigr){}^{I^-,\psi}}
\nc{\Dsmni}{\mathsf{D}\bigl({\mathcal Fl}^{\frac{\infty}{2}}\bigr){}^{I^0}}

\nc{\Pervsmkb}{\overline{\mathsf{Perv}}\bigl({\mathcal Fl}^{\frac{\infty}{2}}\bigr){}^{G^k}}
\nc{\Pervsmgb}{\overline{\mathsf{Perv}}\bigl({\mathcal Fl}^{\frac{\infty}{2}}\bigr){}^{G[[t]]}}
\nc{\Pervsmib}{\overline{\mathsf{Perv}}\bigl({\mathcal Fl}^{\frac{\infty}{2}}\bigr){}^I}
\nc{\Pervsmnib}{\overline{\mathsf{Perv}}\bigl({\mathcal Fl}^{\frac{\infty}{2}}\bigr){}^{I^0}}
\nc{\Pervsmnif}{{}^f\mathsf{Perv}\bigl({\mathcal Fl}^{\frac{\infty}{2}}\bigr){}^{I^0}}

\nc{\Bun}{\on{Bun}}
\nc{\Rep}{\on{Rep}}
\nc{\ConvHecke}{\on{Conv}^{\on{Hecke}}}

\nc{\BunNb}{\overline{\Bun}_{N^-}}
\nc{\BunN}{\Bun_{N^-}}
\nc{\BunBm}{\Bun_{B^-}}

\renc{\mod}{\on{-mod}}
\nc{\modo}{\on{-}\ol{\on{mod}}}

\nc{\clambda}{{\check\lambda}}
\nc{\cLambda}{{\check\Lambda}}
\nc{\cmu}{{\check\mu}}
\nc{\cnu}{{\check\nu}}
\nc{\crho}{{\check\rho}}
\nc{\cG}{{\check G}}
\nc{\cB}{{\check B}}
\nc{\cT}{{\check T}}
\nc{\cN}{{\check N}}

\nc{\tw}{{\tilde{w}}}
\nc{\ol}{\overline}
\nc{\ul}{\underline}
\nc{\wt}{\widetilde}
\nc{\uV}{\underline{V}}
\nc{\uU}{\underline{U}}
\nc{\uCO}{\underline{\CO}}
\nc{\oCZ}{\overset{\circ}\CZ{}}
\nc{\oX}{\overset{\circ}{X}{}}
\nc{\oi}{\ol{i}}
\nc{\ofF}{\overset{\circ}\fF{}}
\nc{\ci}{\overset{\circ}{i}}
\nc{\tboxtimes}{\,\widetilde{\boxtimes}}
\nc{\hr}{\overset{\rightarrow}{h}{}}
\nc{\hl}{\overset{\leftarrow}{h}{}}
\nc{\semiinf}{{\frac{\infty}{2}}}
\nc{\fIC}{{\mathsf{IC}}}
\nc{\bPi}{{\mathbf \Pi}}
\nc{\uBC}{\underline{\BC}}
\nc{\starstar}{\overset{*}\star}
\nc{\starshriek}{\overset{!}\star}
\nc{\sD}{\mathsf D}

\begin{document}

\author{S.~Arkhipov, R.~Bezrukavnikov, A.~Braverman, D.~Gaitsgory, I.~Mirkovi\'c}

\title[Modules over the small quantum group and semi-infinite
flag manifold]{Modules over the small quantum group and semi-infinite
flag manifold}

\address{}

\email{serguei.arkhipov@yale.edu, bezrukav@math.northwestern.edu,\newline
braval@math.brown.edu, gaitsgde@math.harvard.edu,
mirkovic@math.umass.edu}

\begin{abstract}
We develop a theory of perverse sheaves on the semi-infinite flag
manifold $G((t))/N((t))\cdot T[[t]]$, and show that the subcategory of
Iwahori-monodromy perverse sheaves is equivalent to the regular
block of the category of representations of the small quantum group
at an even root of unity.
\end{abstract}

\dedicatory{To V.~Drinfeld on the occasion of his 50th birthday}

\date{April 2005}




\maketitle

\tableofcontents


\section*{Introduction}

\ssec{Motivation}

Let $G$ be a reductive group. The purpose of this paper is to show that a certain
remarkable abelian category $\CA$ can be realized in (at least) three seemingly different
contexts as a category of representations of some sort. This abelian category has a
significance, since it can be thought of as a "local geometric Langlands" category,
corresponding to an unramified local system. Let us try to explain this point, even
though the local geometric Langlands correspondence has not been yet properly
formulated. As a result, the discussion in this subsection will not be rigorous.

Let us recall that the global geometric Langlands correspondence aims to attach
to a local system $\sigma:\pi_1(X)\to \cG$ (here $X$ is a smooth and complete curve)
a perverse sheaf $\CF_\sigma$ on the stack $\Bun_G$, classifying principal
$G$-bundles on $X$; one requires $\CF_\sigma$ to satisfy the Hecke property with
respect to $\cG$.

\medskip

The perverse sheaf $\CF_\sigma$ should be thought of as a "higher" analogue of an unramified
automorphic function $f_\sigma$ with Langlands paramaters given by $\sigma$ (the
latter makes sense, of course, only when the ground field is finite). To simplify the discussion,
let us assume that the unramified automorphic representation $\pi_\sigma$, containing
$f_\sigma$, lies discretely in the corresponding $L_2$ space and, moreover, that all of its local
components are irreducible unramified principal series representations.

Let us now fix a point $x\in X$, and instead of just one automorphic function $f_\sigma$
let us consider the sub-space $(\pi_\sigma)_x\subset \pi_\sigma$, consisting of
vectors invariant with respect to $\underset{x'\neq x}\Pi\, G(\CO_{x'})$. This is a representation
of the locally compact group $G(\CK_x)$ (here for a place $x'\in X$, $\CO_{x'}$ and $\CK_{x'}$
denote the local ring and the local field at this point, respectively).

According to the Langlands philosophy, $(\pi_\sigma)_x$
should be completely determined by the local Galois representation $\sigma_x$. Since
$\sigma$ was assumed unramified, $\sigma_x$ boils down simply to the conjugacy class
of the image of the Frobenius element.

\medskip

Let us now try to guess what a geometric analogue of the vector space $(\pi_\sigma)_x$
might be. Let $^{\infty}\Bun_G$
be the moduli stack of principal $G$-bundles on $X$ with a full level structure at $x$
.

We propose that there should exist an (abelian) category $\CA$,
acted on by $G(\CK_x)$ by functors (here $G(\CK_x)$ is understood as the corresponding
group ind-scheme), and a functor from $\CA$ to the category of perverse sheaves on
$^{\infty}\Bun_G$, whose image consists of perverse sheaves that satisfy the Hecke property
with respect to $\sigma$ an $X-x$.  \footnote{This, rather crude, form of the guess for what the
local geometric Langlands correspondence might be, has been voiced independently by
many people, and we by no means claim primacy in this matter.}

\medskip

The above considerations on the function-theoretic level suggest the following
candidate for $\CA$. Namely, this should be the category of perverse sheaves on the
affine Grassmannian $\Gr=G(\CK_x)/G(\CO_x)$ that satisfy the Hecke property (cf.
\secref{Hecke category} for the precise definition).

Recall now that $(\pi_\sigma)_x$ could also be realized as an (irreducible, spherical)
principal series representation. Therefore, it is tempting to realize the category $\CA$
in terms of perverse sheaves on the semi-infinite flag manifold
$G(\CK_x)/N(\CK_x)\cdot T(\CO_x)$. This is the point of departure for
the present paper.

\medskip

Before we proceed to the description of the concrete problem that is posed and solved here,
let us mention one more incarnation of the category $\CA$. Namely, the Beilinson-Drinfeld
construction of Hecke eigensheaves via quantization of the Hitchin integrable system
suggests, that the category $\CA$ should be also equivalent to the category of modules
over the affine algebra at the critical level, with a fixed central character, corresponding to
some oper on the formal disc around $x$.

This category of representations can indeed be connected to $\CA$. In the
forthcoming work \cite{FG} a functor is defined from the D-module version of
category $\CA$ to a certain category of modules over the affine Kac-Moody algebra 
at the critical level with a fixed central character.  It is conjectured in \cite{FG} 
that this functor is an equivalence of categories. Moreover, it is proved that 
it is fully faithful, and in the next paper the authors of {\it loc.cit.} will 
show that it indeed is an equivalence of categories when resricted to the Iwahori 
equivariant subcategories.

What is unfortunately unavailable at the moment, is a direct link between
critical level representations and the cattegory of sheaves on $G(\CK_x)/N(\CK_x)\cdot T(\CO_x)$.
Such a link, which was forseen by Feigin and Frenkel in \cite{FF}
as a localization-type theorem for sheaves on $G(\CK_x)/N(\CK_x)\cdot T(\CO_x)$, was 
the source of many people's interest in the study of both categories.

\ssec{The present work}

The goal of this paper is to connect the category of Hecke eigen-sheaves on the
affine Grassmannian, denoted $\Catgeom$ (or rather its graded version, denoted
$\Catgeomd$), to the category of perverse sheaves on the semi-infinite flag manifold.
An immediate problem that one runs into is that the latter category
does not a priori makes sense:

The semi-infinite flag manifold, thought of as $G(\CK_x)/N(\CK_x)\cdot T(\CO_x)$,
does not carry an algebro-geometric structure that would allow for the theory
of perverse sheaves, or D-modules, in the way it is known today.

\medskip

We get around this difficulty as follows. We define an "artificial" category $\Pervsm$ that
possesses the natural properties that one expects from the yet non-existing
category of perverse sheaves on $G(\CK_x)/N(\CK_x)\cdot T(\CO_x)$. The approach
to the definition of $\Pervsm$, developed in this paper, was initiated in \cite{FM}, and it
uses a geometric object, denoted $\BunNb$, introduced by Drinfled.

The space $\BunNb$ is a finite-dimensional (or, rather, ind-finite dimensional)
approximation to $G(\CK_x)/N(\CK_x)\cdot T(\CO_x)$, and it has as an input a global curve
$X$. By definition, $\BunNb$ classifies principal $G$-bundles on $X$ endowed
with a possibly degenerate reduction to the maximal unipotent subgroup $N^-$,
and it contains the stack $\Bun_{N^-}$ classifying $N^-$-bundles on $X$ as an
open substack.

The realization of $\Pervsm$ via $\BunNb$ is natural from the geometric Langlands
perspective as well: the space $\BunNb$ is used to define geometric Eisenstein
series by taking the direct image under the natural projection to $\Bun_G$
(cf. \cite{BG}). Therefore, such incarnation of $\Pervsm$ implies the existence
of a functor from $\CA$ to the (derived) category of perverse sheaves on $^\infty\Bun_G$.

\medskip

Having defined the category $\Pervsm$, we have at our disposal a naturally defined
functor from $\Catgeomd$ to it. However, we do not have any real evidence as to
whether this functor should be an equivalence. Quite possibly, to make this functor an
equivalence, one has to modify both categories by imposing some Noetherianness
condions on the $\Catgeomd$ side, and restrictions on the behaviour "at the boundary"
on the $\Pervsm$ side.

The problem arising here is similar to the one in the definition of the Schwarz
space on $G(\CK_x)/N(\CK_x)$ in the function-theoretic context in \cite{BK}. Identifying
the image of $\Catgeomd$ inside $\Pervsm$ appears to be an interesting problem, and it is
closely related to giving a geometric definition of Fourier-transform functors of {\it loc. cit.}

However, if instead of the entire $\Catgeomd$ and $\Pervsm$ we work with the subcategories,
denoted $\Catgeomdni$ and $\Pervsmni$, respectively, consisting of Iwahori-monodromic
objects, the required Noetherian and boundary conditions are easy to spell out, simply by
requiring that our objects have finite length.

Thus, the main result of this paper, \thmref{main}, states that the category, denoted
$\Catgeomdni_{Art}$, consisting of Artinian and Iwahori-monodromic objects in
$\Catgeomd$, is equivalent to the subcategory of Artinian objects in $\Pervsmni$.

\medskip

The method of proof of \thmref{main} relies rather heavily on the specifics of
Iwahori-monodromic situation. Namely, we use the fact that both categories are
hereditary (i.e., in many ways similar to the usual category $\CO$). In particular,
they both have standard and costandard objects, numbered by elements of
the extended affine Weyl group $W_{aff}$, etc.

The hereditary structure on $\Pervsmni$ is evident basically from the stratification of
$G(\CK_x)/N(\CK_x)\cdot T(\CO_x)$ by Iwahori orbits. However, for $\Catgeomdni_{Art}$ this structure
is not so evident, and it comes from another crucial ingredient of this paper, namely,
the equivalence between $\Catgeomdni_{Art}$ and the regular block of the
category of representations of the small quantum group, corresponding to $G$,
at an even root of unity.

The latter equivalence results by combining the main result of \cite{ABG} that links
representations of the big quantum group and perverse sheaves on $\Gr$, and
\cite{AG}, where an explicit relation between the categories of representations of the
big and small quantum group is established.

\medskip

We should point out, however, that the present paper relies formally on neither \cite{ABG},
nor \cite{AG}. We supply purely geometric proofs for all the statements needed to establish the
hereditary property of $\Catgeomdni_{Art}$. But these statements would be rather hard
to guess, had we not had the equivalence with the quantum group as a guide.

As a result, we also obtain that the category of Artinian objects in $\Pervsmni$ is equivalent
to the category $\tu_\ell\mod_0$--the above mentioned regular block in the
category of $\tu_\ell$-modules. This is our \thmref{main quantum}, which concludes
the project of proving such an equivalence, initiated and advanced almost to
the end by M. ~Finkelberg. \footnote{An equivalence between $\tu_\ell\mod_0$ and the
would-be category of Iwahori-monodromic perverse sheaves on $G(\CK_x)/N(\CK_x)\cdot T(\CO_x)$ 
has also been guessed independently by several people, among them, Lusztig and
Feigin-Frenkel, but we could not find a precisely formulated conjecture in the literature.
Our formulation as well as the strategy of the proof are due to Finkelberg.}

\ssec{Contents}

Let us now discuss the organization and contents of the present paper.

\medskip

Section 1 reviews the theory of modules over the big and small quantum groups.

\noindent In Sect. 1.1 we recall the basic definitions related to corresponding
categories of representations, and the quantum Frobenius homomorphism.
In Sect. 1.2 we recall the realization of the category of representations of the small
quantum group as representations of the big quantum group, satisfying the Hecke
property. In Sect. 1.3 we recall the \cite{ABG} equivalence between the regular block
of the category of representations of the big quantum group and Iwahori-monodromic
perverse sheaves on the affine Grassmannian; we also introduce the category of
Hecke eigen-sheaves on the Grassmannian and discuss its relation to the category
of representations of the small quantum group.

\medskip

Section 2 reviews some basic properties of Iwahori-equivariant perverse sheaves
on the affine Grassmannian.
In Sect. 2.1 we give a geometric proof of an irreducibility result on
convolution of certain perverse sheaves, which translates by means of \cite{ABG}
to the Steinberg-type theorem for representations of the quantum group; some
ingredients of the proof will be used later on for a crucial irreducibility result in Sect. 5.3.
In Sect. 2.2 we discuss the baby Whittaker category on the affine Grassmannian and its
relation to a certain Serre quotient category of $\Pervgrni$; the discussion here largely
repeats the one in \cite{AB}. In Sect. 2.3 we apply the results of the previous subsection
to establish a crucial result about cosocles of some costandard objects in $\Pervgrni$;
this result will be essential for the proof of the main theorem.

\medskip

Section 3 is devoted to the study of baby (co)Verma modules over
the small quantum group, which are the building blocks of the
category of its representations. In Sect. 3.1 we translate the
properties of baby co-Verma modules into properties of the
corresponding modules over the big quantum group, satisfying the
Hecke property. In Sect. 3.2 we reprove the corresponding facts
(often by different methods) in the context of Iwahori-monodromic
perverse sheaves on $\Gr$.

\medskip

In Section 4 we discuss the main object of study of this paper, namely, the category
$\Pervsm$, which is a surrogate for the non-existing category of perverse sheaves on
$G(\CK_x)/N(\CK_x)\cdot T(\CO_x)$.

In Sect. 4.1 we discuss the underlying geometric
object---the stack $\BunNb$ along with its numerous variants.
In Sect. 4.2 we finally introduce the category $\Pervsm$,
the main technical ingredient being the factorizability property, observed in \cite{FFKM};
we show that that $\Pervsm$ by and large behaves in the way one expects from
the analogy with $G(\CK_x)/N(\CK_x)\cdot T(\CO_x)$. In Sect. 4.3 we study the most
basic objects in $\Pervsm$, namely, the spherical ones, and show that the resulting
category is semi-simple. Finally, in Sect. 4.4 we discuss the Iwahori-monodromic
subcategory of $\Pervsm$, and prove some results that are parallel to the corresponding
assertions about Iwahori-monodromic sheaves on $\Gr$.

\medskip

As was mentioned above, the category $\Pervsm$ must be acted on by the group ind-scheme
$G(\CK_x)$ by auto-functors. A rigorous incarnation of this phenomenon is the
action of perverse sheaves on $G(\CK_x)$ by Hecke functors (the latter are defined
on the level of the corresponding derived category). In Sect. 5 we study this convolution
action in our realization of $\Pervsm$ via $\BunNb$.

In Sect. 5.1 we define the convolution action and show that it indeed respects the
category $\Pervsm$. In Sect. 5.2 we establish a crucial semi-smallness result that
allows to pass from perverse sheaves on $\Gr$ to $\Pervsm$ (this is largely borrowed
from \cite{FM} and \cite{BG}). In Sect. 5.3 we refine the discussion of the previous
subsection and show that certain convolution diagrams give rise to small
(vs. semi-small) maps, thereby implying certain irreducibility properties. In Sect. 5.4
we establish another important technical result that describes the convolution of
standard objects.

\medskip

Finally, in Section 6 we state and prove the equivalence between
the subcategories of Artinian objects in $\Catgeomdni$ and
$\Pervsmni$. In Sect. 6.1 we define the required functor. In Sect.
6.2 we show that this functor is exact and reduce the equivalence
assertion to a computation of the image of baby co-Verma modules.
In Sect. 6.3 we perform the required calculation using some
information on cosocles of costandard objects in both categories.

\medskip

The conventions adopted in this paper regarding the quantum group follow those
of \cite{AG}. Conventions and notation concerning the affine Grasmannian and Drinfeld's
compactifications follow those of \cite{BG}. To fix the context we will work with varieties
and stacks over the ground field $\BC$, and holonomic D-modules (but we will
still call them perverse sheaves). If $\CY$ is a smooth variety, $\uBC_\CY$ will denote the
(cohomologically shifted) D-module, corresponding to the constant sheaf on it.

\medskip

\noindent{\bf Acknowledgements.}
As was mentioned earlier, the problem solved in this paper was both posed (and
the method of solution was suggested) by M.~Finkelberg back in 1998, when
the authors were at IAS, Princeton, for the special year on geometric representation theory.
We are grateful to him for the permission to publish many of his results and ideas.

We would also like to thank A.~Beilinson, V.~Drinfeld, B.~Feigin, E.~Frenkel and
D.~Kazhdan for sharing their ideas and stimulating discussions.

It is an honour for us to dedicate this paper to Vladimir Drinfeld. Along with numerous
other things in modern mathematics, the three main objects of study in this
paper--quantum groups, Hecke eigen-sheaves and $\BunNb$ were invented by him.

\section{Background: modules over the big and small quantum groups}

\ssec{Basics of quantum groups}

\sssec{Root data}

Let $G$ be a reductive group with connected center. Let $\cG$ be its Langlands
dual; by assumption the derived group of $\cG$ is simply connected.
\footnote{For what follows we could replace $G$ by an isogenous group such
that $[G,G]$ is simply connected.   In this case $\cG$ also has connected center.}

We will denote by $T$ (resp., $\cT$) the Cartan group of $G$
(resp., $\cG$), and by $W$ the Weyl group. We fix Borel subgroups
$B,B^-\subset G$ (resp., $\cB,\cB^-\subset \cG$) and think of $T$ (resp., $\cT$)
as a subgroup of $G$ (resp., $\cG$) equal to their intersection.

We will denote by $\cLambda$ (resp., $\Lambda$) the coweight (resp., weight
lattice) of $G$; by $\cLambda^+$ (resp., $\Lambda^+$) we will denote the
subset of dominant coweights (resp., weights). We will denote by
$\langle\cdot,\cdot,\rangle$ the pairing between the two. We will denote
by $W_{aff}$ the extended Weyl group $W\ltimes \cLambda$.

Let $\CI$ be the set of vertices of the Dynkin graph of $G$;
for $\imath\in \CI$ we will denote by $\check\alpha_\imath\in \cLambda$
(resp., $\alpha_\imath\in \Lambda$) the corresponding simple coroot
(resp., root). We will denote by $\cLambda^{pos}$ (resp., $\Lambda^{pos}$)
the sub-semigroup spanned by positive coroots (resp., roots).

\medskip

Let $(\cdot,\cdot):\on{Span}\{\alpha_\imath\}\otimes \on{Span}\{\alpha_\imath\}\to \BZ$ be
the {\it canonical} inner form. In other words, $||\alpha_\imath||^2=2d_\imath$, where
$d_\imath\in\{1,2,3\}$ is the minimal set of integers such that the matrix
$(\alpha_\imath,\alpha_\jmath):=
d_\imath\cdot \langle \alpha_\imath, \check\alpha_\jmath \rangle$ is symmetric.

We choose a symmetric $W$-invariant form
$(\cdot,\cdot)_\ell:\cLambda\otimes \cLambda\to \BZ$,
such that there exists a {\it sufficiently large positive even integer} $\ell$, divisible by all $d_\imath$,
such that
$$(\check\alpha_\imath,\clambda)_\ell=\ell_\imath\cdot \langle \alpha_\imath,\clambda\rangle,$$
$\forall \clambda\in \cLambda$, where $\ell_\imath=\frac{\ell}{d_\imath}$.

We will denote by $\phi_\ell$ the resulting map $\cLambda\to \Lambda$, and also
the map $T\to \cT$.

\sssec{The big quantum group}

As was mentioned earlier, our conventions regarding representations of the big quantum
group follow those of \cite{AG}. Let $\bU_\ell\mod$ be the category of representations of the
big quantum group, corresponding to $G$ and $\ell$. By definition, objects of this category are
finite-dimensional vector spaces, acted on by the algebraic group $T$, and the
operators $E_\imath$, $F_\imath$, $E^{(\ell_\imath)}_\imath$, $F^{(\ell_\imath)}_\imath$,
that satisfy the well-known relations. The category $\bU_\ell\mod$ has a natural
monoidal structure.

We will denote by $\bU_\ell\modo$ the ind-completion of $\bU_\ell\mod$. I.e.,
this is the category of infinite-dimensional vector spaces, acted on by the same
set of operators, which can be represented as unions of finite-dimensional
sub-representations.

\medskip

Let $\bB^-_\ell\mod$ be the category of representations of the "negative quantum
Borel". I.e., objects of this category are finite-dimensional vector spaces,
acted on by the algebraic group $T$, and the operators $F_\imath$,
$F^{(\ell_\imath)}_\imath$, which satisfy the same relations. This is also a
monoidal category and there exists a natural forgetful monoidal functor
$\Res^{\bU_\ell}_{\bB^-_\ell}:\bU_\ell\mod\to \bB^-_\ell\mod$. This functor
admits a right adjoint, denoted by $\Ind^{\bU_\ell}_{\bB^-_\ell}$.

In addition there exists a natural functor $\Rep(T)\to \bB^-_\ell\mod$, where we let
the operators $F_\imath$, $F^{(\ell_\imath)}_\imath$ act trivially on the
corresponding vector space.

For $\lambda\in \Lambda$ we let $\bW^\lambda\in \bU_\ell\mod$ be the dual Weyl
module defined as $\Ind^{\bU_\ell}_{\bB^-_\ell}(\bC^\lambda)$, where $\bC^\lambda$
is the $1$-dimensional representation of $T$, corresponding to $\lambda$.
It is known that $\bW^\lambda\neq 0$ if and only if $\lambda\in \Lambda^+$. It is
also known that $\bW^\lambda$ admits a unique irreducible submodule, denoted
$\bL^\lambda$, and this establishes a bijection between $\Lambda^+$ and the
set of irreducibles in $\bU_\ell\mod$.

\medskip

As every Artinian category, $\bU_\ell\mod$ splits into a direct sum of
indecomposable Artinian categories, called blocks. Slightly deviating
from the accepted conventions, we will denote by $\bU_\ell\mod_0$ the
direct summand of $\bU_\ell\mod$ that contains the irreducibles
$\bL^\lambda$ for $\lambda$ of the form
$$w(\rho)-\rho+\phi_\ell(\clambda),$$
$w\in W$, $\clambda\in \cLambda$.

We will denote by $\bU_\ell\modo_0$ the ind-completion of
$\bU_\ell\mod_0$, which is a direct summand in $\bU_\ell\modo$.

\sssec{Quantum Frobenius homomorphism}

Let $\Rep(\cG)$ denote the category of finite-dimensional
representations of $\cG$. Following \cite{Lu1} there exists
a monoidal functor
$$\on{Fr}:\Rep(\cG)\to \bU_\ell\mod,$$
defined as follows. For $V\in \Rep(\cG)$, the representation
$\on{Fr}(V)$ occurs on the same underlying vector space, denoted $\uV$,
and the action of $T$ is given via $\phi_\ell:T\to \cT$. The operators
$E_\imath,F_\imath$ act trivially, and $E^{(\ell_\imath)}_\imath,F^{(\ell_\imath)}_\imath$
act via the Chevalley generators $e_\imath,f_\imath\in \check\fg$.

It is known that the functor $\on{Fr}$ is fully faithful. Moreover, for
$\clambda\in \cLambda^+$
$$\on{Fr}(V^\clambda)\simeq \bL^{\phi_\ell(\clambda)},$$
where $V^\clambda$ denotes the corresponding irreducible representation of $\cG$.

\medskip

Let recall that a dominant weight $\lambda$ is called restricted if $\forall \imath\in \CI$
$$\langle \lambda,\check\alpha_\imath\rangle < \ell_\imath.$$

We have the following fundamental result:
\begin{thm} \label{steinberg}
If $\clambda\in \Lambda^+$ is restricted, then for every $\cmu\in \cLambda^+$
$$\on{Fr}(V^\cmu)\otimes \bL^\lambda \simeq \bL^{\lambda+\phi_\ell(\cmu)}.$$
\end{thm}

Since every $\lambda\in \Lambda^+$ can be written as $\clambda_1+\clambda_2$
with $\clambda_1$ restricted and $\clambda_2$ in the image of $\cLambda^+$,
the above theorem describes all irreducibles in $\bU_\ell\mod$. (Note that the
decomposition of a weight as above is unique modulo elements $\nu\in \Lambda$,
orthogonal to all roots, i.e., those for which $\bL^\nu$ is $1$-dimensional.)

\begin{cor}  \label{block preserved}
The functor $M\mapsto \on{Fr}(V)\otimes M:\bU_\ell\mod\to\bU_\ell\mod$
preserves $\bU_\ell\mod_0$.
\end{cor}

\sssec{The graded small quantum group}   \label{intro to baby Verma}

We define the category of representations of the
graded small quantum group $\tu_\ell\mod$ to consist of finite-dimensional
vector spaces, acted on by the algebraic group $T$ and the operators
$E_\imath$, $F_\imath$, satisfying the usual relations. This is also
a monoidal category, and we have a monoidal forgetful functor
$\Res^{\bU_\ell}_{\tu_\ell}:\bU_\ell\mod\to \tu_\ell\mod$.

In addition, we have a fully-faithful functor $\Rep(\cT)\to \tu_\ell\mod$.
By a slight abuse of notation we will denote by $\bC^\cmu$ the
$1$-dimensional module over $\tu_\ell$, corresponding to $\cmu\in \cLambda$.

Let $\tb^-_\ell\mod$ be the category of representations of the
corresponding "graded small negative Borel subgroup". I.e.,
this is the category of vector spaces, acted on by $T$ and the
$F_\imath$'s, satisfying the same relations. We will denote by
$\Res^{\tu_\ell}_{\tb^-_\ell}$ the forgetful functor
$\tu_\ell\mod\to \tb^-_\ell\mod$ and by $\Ind^{\tu_\ell}_{\tb^-_\ell}$
(resp., $\on{Coind}^{\tu_\ell}_{\tb^-_\ell}$) its right (resp., left) adjoint.
We also have a functor $\Rep(T)\to \tb^-_\ell\mod$.

\begin{lem}   \label{Frobenius algebra}
Both functors $\Ind^{\tu_\ell}_{\tb^-_\ell}$ and $\on{Coind}^{\tu_\ell}_{\tb^-_\ell}$
are exact and faithful and for a character $\lambda$ of $T$,
$$\Ind^{\tu_\ell}_{\tb^-_\ell}(\bC^\lambda)\simeq
\on{Coind}^{\tu_\ell}_{\tb^-_\ell}(\bC^{\lambda-\phi_\ell(2\crho)+2\rho}).$$
\end{lem}

We will denote the module
$\Ind^{\tu_\ell}_{\tb^-_\ell}(\bC^\lambda)$ by $\tM^\lambda$ and
call it the baby co-Verma module of highest weight $\lambda$.
 One easily shows that the socle of each
$\tM^\lambda$ is simple. We will denote the corresponding
irreducible by $\tL^\lambda$. Thus we obtain a bijection between
$\Lambda$ and the set of irreducibles in $\tu_\ell\mod$.

For $\cmu\in \cLambda$, we have:
$$\tM^{\lambda+\phi_\ell(\check\mu)}\simeq \bC^\cmu\otimes
\tM^\lambda \text{ and }
\tL^{\lambda+\phi_\ell(\check\mu)}\simeq \bC^\cmu\otimes
\tL^\lambda.$$

In addition, we have the following result:

\begin{prop}
If $\lambda$ is dominant and restricted, $$\tL^\lambda\simeq
\Res^{\bU_\ell}_{\tu_\ell}(\bL^\lambda).$$
\end{prop}

\medskip

Being Artinian, the category $\tu_\ell\mod$ also admits a decomposition into blocks.
We will denote by $\tu_\ell\mod_0$ the direct summand of $\tu_\ell\mod$ that contains
the irreducibles $\tL^\lambda$ for $\lambda$ of the form $w(\rho)-\rho+\phi_\ell(\clambda)$
$w\in W$, $\clambda\in \cLambda$.

\begin{lem}
The sub-category $\bU_\ell\mod_0\subset \bU_\ell\mod$ is the preimage
of $\tu_\ell\mod_0\subset \tu_\ell\mod_0$ under the forgetful functor
$\Res^{\bU_\ell}_{\tu_\ell}$,
\end{lem}

Finally, we will denote by $\tu_\ell\modo$ (resp., $\tu_\ell\modo_0$)
the ind-completion of $\tu_\ell\mod$ (resp., $\tu_\ell\mod_0$).

\sssec{The non-graded small quantum group}

We define the category $\fu_\ell\mod$ to consist of finite-dimensional
vector spaces, acted on by the group $T_\ell:=\on{ker}(\phi_\ell:T\to \cT)$,
and the operators $K_\imath\cdot E_\imath$, $F_\imath$, subject to the usual relations.
Note that $\fu_\ell\mod$ is {\it not} a monoidal category; however, we
have a well-defined functor of tensor product on the right by an object of
$\tu_\ell\mod$:
$$N\in \fu_\ell\mod,M\in \tu_\ell\mod\mapsto N\otimes \Res^{\tu_\ell}_{\fu_\ell}(M).$$

The following proposition describes the relation between the small quantum group
and the quantum Frobenius homomorphism:

\begin{prop}   \label{quantum Frob and small}  \hfill

\smallskip

\noindent{\em(1)}
For $M\in \tu_\ell\mod$ and $\clambda\in \cLambda$,
$$\Res^{\tu_\ell}_{\fu_\ell}\Bigl(\bC^\clambda\otimes M\Bigr)\simeq
\Res^{\tu_\ell}_{\fu_\ell}(M).$$

\smallskip

\noindent{\em(2)}
For $M$ as above the maximal trivial sub- (resp., quotient-) object $N'$
of $\Res^{\tu_\ell}_{\fu_\ell}(M)$ comes from a sub- (resp., quotient-) object $M'$ of $M$,
which is in the image of the functor $\Rep(\cT)\to \tu_\ell\mod$.

\smallskip

\noindent{\em(3)}
For $M\in \bU_\ell\mod$, $V\in \Rep(\cG)$,
$$\Res^{\bU_\ell}_{\tu_\ell}\Bigl(\on{Fr}(V)\otimes M\Bigr)\simeq
\underset{\cnu}\oplus\, \bC^\cnu\otimes \Res^{\bU_\ell}_{\tu_\ell}(M)\otimes \uV(\cnu),$$
where $\uV(\cnu)$ denotes the $\cnu$-weight space of $V$, and $\bC^\cnu$
the corresponding $1$-dimensional representation of $\tu_\ell$.

\smallskip

\noindent{\em(4)}
For an object $M\in \bU_\ell\mod$ the maximal trivial sub- (resp., quotient-) object $N'$
of $\Res^{\bU_\ell}_{\fu_\ell}(M)$, comes from a sub- (resp., quotient-) object $M'$ of $M$,
which is in the image of the functor $\on{Fr}$.

\end{prop}

Let $\fb^-_\ell\mod$ be the category consisting of finite-dimensional vector spaces,
acted on by the group $T_\ell$ and the operators $F_\imath$, satisfying
the usual relations. We have the evident functor $\Rep(T_\ell)\to \fb^-_\ell\mod$,
such that the analog of \lemref{Frobenius algebra} holds. For a character
$\ol{\lambda}:T_\ell\to \BC^*$ we will denote by $M^{\ol{\lambda}}$ the module
$\Ind^{\fu_\ell}_{\fb^-_\ell}(\BC^{\ol{\lambda}})$.

We have:

\begin{lem}
For a character $\lambda\in \Lambda$ we have:

\smallskip

\noindent{\em (1)} $\Res^{\tu_\ell}_{\fu_\ell}(\tM^\lambda)\simeq M^{\ol{\lambda}}$,
where $\ol\lambda$ is the restriction of $\lambda$ to $T_\ell$.

\smallskip

\noindent{\em (2)} The module $L^{\ol\lambda}:=\Res^{\tu_\ell}_{\fu_\ell}(\tL^\lambda)$
depends only on the class of $\lambda$ modulo $\phi_\ell(\cLambda)$, and is
irreducible. Moreover, these are all the irreducibles in $\fu_\ell\mod$.

\end{lem}

Let $\fu_\ell\mod_0$ be the direct summand of $\fu_\ell\mod$, that contains the
trivial representation.

\begin{lem}
The subcategory $\tu_\ell\mod_0\subset \tu_\ell\mod$ is the preimage of
$\fu_\ell\mod_0\subset \fu_\ell\mod$ under the forgetful functor.
\end{lem}

We will denote by $\fu_\ell\modo$ (resp., $\fu_\ell\modo_0$) the ind-completion of
$\fu_\ell\mod$ (resp., $\fu_\ell\mod_0$).

\medskip

In the sequel we will need the following assertion:

\begin{prop} \label{quant Frob for B}
There exists a fully-faithful functor $\on{Fr}_{B^-}:\Rep(\cB^-)\to \bB^-_\ell\mod$,
such that

\smallskip

\noindent{\em (1)} We have a commutative diagram of functors.
$$
\CD
\Rep(\cG) @>{\on{Fr}}>> \bU_\ell\mod \\
@V{\Res^{\cG}_{\cB^-}}VV  @V{\Res^{\bU_\ell}_{\bB^-_\ell}}VV \\
\Rep(\cB^-) @>{\on{Fr}_B}>> \bB^-_\ell\mod \\
@AAA   @AAA  \\
\Rep(\cT) @>{\phi_\ell}>> \Rep(T)
\endCD
$$

\smallskip

\noindent{\em (2)}  For $N\in \bB^-_\ell\mod$ the maximal sub-
(resp., quotient-) space of $N$, on which $\fb^-_\ell$ acts trivially, is
a sub- (resp., quotient-) module, which lies in the image of the
functor $\on{Fr}_B$.

\end{prop}

\sssec{Weyl group action}   \label{Weyl group action}

Following Lusztig, to every element $w$ of the Weyl group we can attach an
invertible operator acting functorially on the vector space underlying every
object of $\bU_\ell\mod$, or which is the same, an automorphism of the
forgetful functor $\bU_\ell\mod\to \on{Vect}$. This automorphism is well-defined
modulo elements of $T$.

This construction can be reformulated as follows. To every $w\in W$ we attach
a self-functor $\sF_w:\bU_\ell\mod\to \bU_\ell\mod$, that commutes with
the forgetful functor to vector spaces, and an isomorphism
$$w_\ell:\on{Id}_{\bU_\ell\mod}\Rightarrow \sF_w.$$

Restricting these data to the sub-category $\Rep(\cG)\subset \bU_\ell\mod$
we obtain that the pair $(\sF_w,w_\ell)$ gives rise to an element
$w_\cG\in \cG$ that normalizes $\cT$.

\begin{lem}  \hfill  \label{W twisting functors}

\smallskip

\noindent{\em (1)}
There exists a monoidal self-equivalence $\sF_w:\tu_\ell\mod\to \tu_\ell\mod$
that commutes with the restriction functor $\bU_\ell\mod\to \tu_\ell\mod$.

\smallskip

\noindent{\em (2)}
There exists a self-equivalence $\sF_w:\fu_\ell\mod\to \fu_\ell\mod$,
compatible with the functor tensor product functor
$\fu_\ell\mod\times \tu_\ell\mod\to \fu_\ell\mod$.

\end{lem}

We will return to the discussion of functors $\sF_w$ in \secref{further Weyl}.

\medskip

For an element $w\in W$ let $\tb^{w,-}_\ell$ be the corresponding subalgebra
of $\tu_\ell$. Let us denote by $^w\tM^\lambda$ the
$\tu_\ell$-module induced from the $\tb^{w,-}_\ell$-character $\bC^\lambda$.
For $w=1$ we recover $\tM^\lambda$.  We have:
$$^w\tM^{w(\lambda)}\simeq \sF_w(\tM^\lambda).$$

\medskip

As in \lemref{Frobenius algebra},
$$^w\tM^\lambda \simeq
\on{Coind}_{\tb^{w,-}_\ell}^{\tu_\ell}(\bC^{\lambda-w(\phi_\ell(2\check\rho)-2\rho)}).$$

In particular, the module
$^{w_0}\tM^{\lambda-\phi_\ell(2\check\rho)+2\rho}$ is isomorphic
to what is usually called the baby Verma module with highest
weight $\lambda$. Since all
$\on{Coind}_{\tb^{w,-}_\ell}^{\tu_\ell}(\bC^{\lambda})$ have
simple cosocles, we deduce that all twisted baby co-Verma modules
also have simple cosocles.

\ssec{Modules over $\fu_\ell$ as Hecke-proper modules over $\bU_\ell$}

\sssec{The Hecke categories}

Following \cite{AG}, we introduce the category $\Catg$ to consist of pairs
$$(M\in \bU_\ell\modo,\{\alpha_V,\,\,\forall\, V\in \Rep(\cG)\}),$$
where each $\alpha_V$ is a map of $\bU_\ell$-modules
$$\alpha_V: \on{Fr}(V)\otimes M\to M\otimes \underline{V}$$
(for $V\in \Rep(\cG)$, the notation
$\underline{V}$ stands for the underlying vector space),
such that
\begin{itemize}

\item For $V=\BC$, $\alpha_V: M\to M$ is the identity map.

\item For a map $V_1\to V_2$, the diagram
$$
\CD
\on{Fr}(V_1)\otimes M  @>{\alpha_{V_1}}>>  M\otimes \underline{V_1}  \\
@VVV        @VVV    \\
\on{Fr}(V_2)\otimes M  @>{\alpha_{V_2}}>>  M\otimes \underline{V_2}
\endCD
$$ commutes.

\item
A compatibility with tensor products holds in the sense that the map
$$\on{Fr}(V_1)\otimes \on{Fr}(V_2) \otimes M\to \on{Fr}(V_1\otimes V_2)
\otimes M\overset{\alpha_{V_1\otimes V_2}}
\longrightarrow  M\otimes  \underline{V_1\otimes V_2}
\to M\otimes \underline{V_1}\otimes \underline{V_2}$$
equals
$$\on{Fr}(V_1)\otimes \on{Fr}(V_2)\otimes M\overset{\on{id}\otimes\alpha_{V_2}}
\longrightarrow \on{Fr}(V_1)\otimes M\otimes \underline{V_2}\overset{\alpha_{V_1}}\simeq
M \otimes \underline{V_1}\otimes \underline{V_2}.$$
\end{itemize}

It was shown in \cite{AG} that the maps $\alpha_V$ are necessarily isomorphisms.

\medskip

Morphisms in this category between $(M,\alpha_V)$ and $(M',\alpha'_V)$
are $\bU_\ell$-module maps $M\to M'$ preserving the above structures.
Evidently, $\Catg$ is an abelian category.

The  main result of \cite{AG} is the following theorem:

\begin{thm} \label{main AG}
The category $\Catg$ is naturally equivalent to  $\fu_\ell\modo$.
\end{thm}

We recall that the functors $\Catg\rightleftarrows \fu_\ell\modo$ are defined as follows.
To $N\in \fu_\ell\modo$ we attach the object in $\bU_\ell\modo$ by taking
$\Ind^{\bU_\ell}_{\fu_\ell}(N)$. It satisfies the Hecke condition due to
\propref{quantum Frob and small}.

Vice versa, given an object $M$ of $\Catg$,
the restriction $\Res^{\bU_\ell}_{\fu_\ell}(M)$ is acted on by the algebra
$\CO_{\cG}$, and the corresponding object of $\fu_\ell\modo$ is
by definition the tensor product $\Res^{\bU_\ell}_{\fu_\ell}(M)\underset{\CO_{\cG}}\otimes \BC_1$,
where $\BC_1$ is the skyscraper at $1\in \cG$.

\medskip

A typical example of an object of $\Catg$ is obtained by taking
$\on{Fr}(R_\cG)\otimes M$, where $R_\cG$ is the algebra of functions
on $\cG$, regarded as a representation of $\cG$, and $M\in \bU_\ell\modo$.

We say that an object of $\Catg$ is finitely generated if it admits a surjection
from an object of the above form for $M\in \bU_\ell\mod$.
Evidently, the subcategory of finitely
generated objects of $\Catg$, denoted $\Catg_{f.g.}$, transforms under the
equivalence of \thmref{main AG} to $\fu_\ell\mod$. In particular, this subcategory
is Artinian, and $\Catg$ is the ind-completion of $\Catg_{f.g.}$.

\medskip

Consider the subcategory $\Catg_0$ of $\Catg$, equal to the preimage of
$\bU_\ell\mod_0$ under the forgetful functor. According to \cite{AG}, the
equivalence of \thmref{main AG} induces an equivalence between
$\Catg_0$ and $\fu_\ell\modo_0$. We will denote by $\Catg_{0,f.g.}$
the intersection of $\Catg_0$ with $\Catg_{f.g.}$; this category is
equivalent to $\fu_\ell\mod_0$.

\sssec{Hecke categories, graded version}

We define the category $\Catgd$ as follows. Its objects are
$\cLambda$-graded objects $\tM=\oplus\, M_\cnu$ of $\bU_\ell\mod$,
each endowed with a collection of grading-preserving maps
$\alpha_V,\, \forall V\in \Rep(\cG)$
$$\on{Fr}(V)\otimes \tM\simeq \tM\otimes \uV$$
(where the grading on the LHS is induced from that on $M$, and
on the RHS is diagonal with respect to the action of $\cT$ on $\uV$),
which satisfy the same conditions as in the definition of $\Catg$.

Maps in this category are grading preserving maps in $\bU_\ell\mod$
that intertwine the corresponding $\alpha_V$'s. We have the following
graded version of \thmref{main AG}:

\begin{thm} \label{AG graded}
The category $\Catgd$ is equivalent to $\tu_\ell\modo$. The forgetful functor
$\Catgd\to \Catg$ identifies under this equivalence with $\Res^{\tu_\ell}_{\fu_\ell}$.
\end{thm}

We will denote by $\tM\mapsto \tM\{\cmu\}$ the functor on $\Catgd$ given by
the shift of grading by $\cmu\in \cLambda$. Under the equivalence of
\thmref{AG graded} this functor transforms to the functor
$N\mapsto \bC^\cmu\otimes N$.

\medskip

Let $\Catgd_0$ be the preimage in $\Catgd$ of $\bU_\ell\mod_0$ under
the obvious forgetful functor. This subcategory goes over under the
equivalence of \thmref{AG graded} to $\tu_\ell\modo_0$.

\medskip

Let $\tR_\cG$ be the algebra of functions of $\cG$, regarded as a
$\cLambda$-graded representation of $\cG$ (the grading comes from the action of $\cG$
on itself on the right). A typical example of an object of
$\Catgd$ is $\on{Fr}(\tR_\cG)\otimes M$ for $M$ being a $\cLambda$-graded
object of $\bU_\ell\mod$.

We will denote by $\Catgd_{f.g.}$ (resp., $\Catgd_{0,f.g.}$) the corresponding
subcategories of finitely generated objects. These subcategories
transform under the equivalence of \thmref{AG graded} to $\tu_\ell\mod$
and $\tu_\ell\mod_0$, respectively.

\sssec{Action of the dual group}  \label{further Weyl}

Note that the equivalence of \thmref{main AG} makes it explicit
that the category $\fu_\ell\modo$ carries an action of the dual group
by auto-equivalences. The latter means that to every $N\in \fu_\ell\modo$
we can attach a family $^{\cG}N$ of objects of
$\fu_\ell\modo$, parametrized by $\cG$, such that the natural associativity
condition holds.

The corresponding family is defined in the language of $\Catg$ as follows.
For $(M,\{\alpha_V\})\in \Catg$ its fiber at
$\bg\in \cG$ is $(M,\{\bg\cdot \alpha_V\})$, where each $\bg\cdot \alpha_V$
is the composition of $\alpha$ with the automorphism induced by $\bg$
on $\uV$. We will use the notation $N\mapsto {}^{\bg}N$ for these functors.

\medskip

Consider now the case of $\tu_\ell\modo\simeq \Catgd$. In this case
we do not have an action of the entire $\cG$ on the category, but rather of the
normalizer of the Cartan subgroup $\cT$, due to the grading condition.

\begin{lem}
For a pair $(\sF_w,w_\ell)$ as above, the functors $\sF_w$
$$\fu_\ell\modo\to \fu_\ell\modo \text{ and }
\tu_\ell\modo\to \tu_\ell\modo$$
are naturally isomorphic to the functors
$N\mapsto {}^{w_\cG}N$,
where $w_\cG$ is the corresponding element in the normalizer
of $\cT$ in $\cG$.
\end{lem}

\sssec{Compatibility with duality}    \label{comp with duality}

Recall that both categories $\bU_\ell\mod$ and $\tu_\ell\mod$ carry
a canonical self anti-equivalence (contragredient duality),
$M\mapsto M^\vee$, compatible
with the forgetful functor $\Res^{\bU_\ell}_{\tu_\ell}$.
We would like to express the duality functor on $\tu_\ell\mod$
in terms of $\Catgd_{f.g.}$.

Thus, let $N$ be an object of $\tu_\ell\mod$ and $\tM\in \Catgd$ the object
corresponding to it under \thmref{AG graded}. Since $\tM\in \Catgd_{f.g.}$, it
can be represented as the cokernel of an arrow
$$\underset{i}\oplus\, \on{Fr}(\tR_G)\otimes M^1_i\{\cmu^1_i\}\to
\underset{j}\oplus\, \on{Fr}(\tR_G)\otimes M^2_j\{\cmu^2_j\},$$
where the indices $i$ and $j$ run over some finite sets, and $M^1_i,M^2_j$
are objects of $\bU_\ell\mod$.

An arrow as above comes from a system of maps in $\bU_\ell\mod$
$$M^1_i\to \on{Fr}(V^{i,j})\otimes M^2_j\otimes (\uV^{i,j})^*(\cmu_j^2-\cmu_i^1),$$
where $V^{i,j}$ are some finite-dimensional representations of
$\cG$. By adjunction, we obtain a system of maps
$$\on{Fr}\Bigl((V^{i,j})^*\Bigr)\otimes M^1_i \otimes \uV^{i,j}(\cmu_j^1-\cmu_i^2)\to M^2_j,$$
and applying the duality,
$$\Bigl(M^2_j\Bigr){}^\vee\to \on{Fr}\left(\Bigl((V^{i,j})^*\Bigr){}^\vee\right)\otimes
\Bigl(M^1_j\Bigr){}^\vee \otimes (\uV^{i,j})^*(\cmu_j^2-\cmu_i^1).$$

Note that for a representation $V$ of $\cG$,
$$\underline{V^\vee}(\cmu)\simeq \uV^*(-\cmu).$$ Hence, if we set
$U^{i,j}=\Bigl((V^{i,j})^*\Bigr)^\vee$, we obtain a system of maps map
$$\Bigl(M^2_j\Bigr){}^\vee\to \on{Fr}(U^{i,j})\otimes \Bigl(M^1_j\Bigr){}^\vee\otimes
(\uU^{i,j})^*(\cmu_1-\cmu_2),$$
which in turn gives rise to a map in $\Catgd$:
$$\underset{j}\oplus\, \on{Fr}(\tR_G)\otimes \Bigl(M^2_j\Bigr){}^\vee\{\cmu_j^2\}\to
\underset{i}\oplus\, \on{Fr}(\tR_G)\otimes \Bigl(M^1_i\Bigr){}^\vee\{\cmu_i^1\}.$$

Then $N^\vee$ corresponds to the object in $\Catgd$ equal to
the kernel of the latter map.

\ssec{Realization via the affine Grassmannian}

\sssec{}

Let $\Gr\simeq G((t))/G[[t]]$ be the affine Grassmannian corresponding to $G$,
and let $\Pervgr$ denote the category of perverse sheaves on it.

Let $\Sph$ denote the category of $G[[t]]$-equivariant perverse sheaves on $\Gr$.
We recall that $\Sph$ is naturally a monoidal category that acts on $\Pervgr$
by convolution functors:
$$\CF\in \Pervgr,\, \CS\in \Sph\mapsto \CF\star \CS.$$

Moreover, $\Sph$ possesses a natural commutativity constraint, and as
a tensor category it is equivalent to $\Rep(\cG)$. We will denote this equivalence
by $V\in \Rep(\cG)\mapsto \CV\in \Sph$. Under this equivalence, the irreducible
representation $V^\clambda$ goes over to $\CV^\clambda=\IC_{\clambda,\Gr}$,
where the latter is the IC sheaf on the closure of the orbit
$\Gr^\clambda=G[[t]]\cdot \clambda$.

\sssec{}

For $k\in \BN$ we will denote by $G^k$ the corresponding congruence
subgroup in $G[[t]]$, and by $\Pervgrk$ the category of $G^k$-equivariant
perverse sheaves on $\Gr$. For $k=0$ we recover $\Sph$; for $k>0$
this is a full subcategory of $\Pervgr$, stable under extensions, since
$G^k$ is pro-unipotent.

Let $I$ (resp., $I^0$) be the Iwahori subgroup of $G$ (resp., its unipotent
radical). We will denote by $\Pervgri$, $\Pervgrni$, $\Dgri$, $\Dgrni$
the corresponding categories of equivariant perverse sheaves and
triangulated categories.

\medskip

Recall that $I$-orbits on $\Gr$ are parametrized by $W_{aff}/W$, which
we will identify with the set of elements in $W_{aff}$, right-minimal with
respect to $W$. Any such element $\tw$ can be uniquely written as
$$\tw=w\cdot \clambda,$$ where
$w\in W$, $\clambda\in \cLambda^+$. The condition of being right-minimal with
respect to $W$ implies that whenever for some $\imath\in \CI$, we have
$\langle \alpha_\imath,\clambda\rangle=0$, then $w(\alpha_\imath)\in \Lambda^{pos}$.

\medskip

For $\tw$ as above we will denote by $\IC_{\tw,\Gr}$ the IC sheaf on the closure
of the corresponding $I$-orbit. By $\CW^{*,\tw}$ (resp., $\CW^{!,\tw}$) we will
denote the corresponding costandard (resp., standard) objects corresponding
to the extension by * (resp., !) of the constant perverse sheaf on this orbit.

Since $I$-orbits and $I^0$-orbits on $\Gr$ coincide, the irreducibles in
$\Pervgri$ are the same as in $\Pervgrni$.

\medskip

In the sequel we will also need some notation pertaining to the affine
flag variety $\Fl=G((t))/I$. We will denote by $\Pervfl$ (resp., $\Pervfli$)
the category of perverse (resp., $I$-equivariant) sheaves on $\Fl$, and
by $\Dfl$ (resp., $\Dfli$) the corresponding triangulated
category.

The category $\Dfli$ has a natural monoidal structure, and it acts
by convolution on $\Dfl$. In addition, we have a natural convolution
functor
$$\Dfl\times \Dgri\to \Dgr.$$

For $\tw\in W_{aff}$ we will denote by $j_{*,\tw}$ (resp., $j_{!,\tw}$) the
costandard (resp., standard) object in $\Pervfli$ attached to the
corresponding $I$-orbit on $\Fl$. We have:
$$j_{*,\tw_1}\star j_{*,\tw_2}=j_{*,\tw_1\cdot \tw_2} \text{ and }
j_{!,\tw_1}\star j_{!,\tw_2}=j_{!,\tw_1\cdot \tw_2},$$
whenever $l(\tw_1)+l(\tw_2)=l(\tw_1\cdot \tw_2)$, where $l(\cdot)$
is the length function on $W_{aff}$. Morover, if $\tw$ is right
$W$-minimal,
$$j_{*,\tw}\star \delta_{1,\Gr}\simeq \CW^{*,\tw} \text{ and }
j_{!,\tw}\star \delta_{1,\Gr}\simeq \CW^{!,\tw}.$$

\sssec{}

According to \cite{KT} combined with \cite{KL} (or, alternatively
by  \cite{ABG}, adapted to the even root of unity case), we have
the following:

\begin{thm} \label{ABG}
There exists an equivalence of categories
$$\Loc:\bU_\ell\mod_0\to \Pervgrni,$$
such that the functor
$$\bU_\ell\mod_0\times \Rep(\cG)\to \bU_\ell\mod_0:
M,V\mapsto \on{Fr}(V)\otimes M$$ identifies
with
$$\Pervgrni\times \Sph\to \Pervgrni: \CS,V\mapsto \CS\star \CV.$$
Moreover, the contragredient duality functor on $\bU_\ell\mod$ goes
over to Verdier duality on $\Pervgrni$.
\end{thm}

Let us describe the image of irreducibles under this equivalence.
If $\lambda\in \Lambda^+$ is such that $\bL^\lambda\in \bU_\ell\mod_0$,
we can uniquely write
$$\lambda=\phi_\ell(\clambda)+w^{-1}(\rho)-\rho,$$
where $\clambda\in \cLambda$. In this case $\tw:=w\cdot \clambda\in W_{aff}$
is right $W$-minimal. Then
$$\on{Loc}(\bL^\lambda)\simeq \IC_{\tw,\Gr} \text{ and }
\on{Loc}(\bW^\lambda)\simeq \CW^{*,\tw}.$$

\medskip

Note also that a weight $\lambda$ as above is restricted if and only if
the pair $(\clambda,w)$ satisfies the following:
$$
\begin{cases}
& \langle \alpha_i,\clambda\rangle=0 \text{ if } w(\alpha_i)\in \cLambda^{pos} \\
& \langle \alpha_i,\clambda\rangle=1 \text{ if } -w(\alpha_i)\in \cLambda^{pos}.
\end{cases}
$$
Hence, for each $w$, the corresponding element $\clambda$ is well-defined
modulo characters of $\cG/[\cG,\cG]$ (which are the same as cocharacters
of $Z(G)$). We will make such a choice and denote
the corresponding irreducible in $\Pervgri$ by $\CL^w$.  We will assume
that for $w=1$, $\CL^w=\delta_{1,\Gr}$. Note that $\CL^{w_0}\simeq
\IC_{w_0\cdot \crho',\Gr}$, where $\crho'$ is {\it some} element of $\cLambda$,
for which $\langle \alpha_\imath,\crho'\rangle=1$ for $\forall \imath\in \CI$.
Such $\crho'$ exists due to the assumption that the center of $G$ is connected.
Note that $2\crho'$ is not in general equal to $2\crho$, the latter being the
sum of positive coroots.

The following is a corollary of \thmref{steinberg} combined with the equivalence of \thmref{ABG}:

\begin{thm}  \label{geometric steinberg}
\hfill

\smallskip

\noindent{\em (1)}
For any $w$ and $\cmu\in \cLambda^+$,
the convolution $\CL^w\star \IC_{\cmu,\Gr}$ is
irreducible and isomorphic to $\IC_{w\cdot (\clambda+\cmu),\Gr}$,
if $\CL^w=\IC_{w\cdot\clambda,\Gr}$.

\smallskip

\noindent{\em (2)} Any irreducible object of $\Pervgri$ has the form
$\CL^w\star \IC_{\cmu,\Gr}$ for unique $w$ and $\cmu$.

\end{thm}

For completeness, in the next section we will give a purely geometric proof of this result.

\sssec{Hecke categories}   \label{Hecke category}

Let $\Pervgrb$ denote the ind-completion of $\Pervgr$.
Let $\Catgeom$ denote the category, whose objects are pairs
$$(\CS\in \Pervgrb,\,\{\alpha_V,\,\,\forall\, V\in \Rep(\cG)\}),$$
where each $\alpha_V$ is a map
$$\CS\star \CV\to \uV\otimes \CS,$$
such that the collection $\{\alpha_V\}$ satisfies the same compatibility
conditions as in the definition of $\Catg$. As in the case of the quantum
group, one shows that the maps $\alpha_V$ are then automatically
isomorphisms.

Morphisms between $(\CS^1,\{\alpha^1_V\})$ and $(\CS^2,\{\alpha^2_V\})$
are maps $\CS^1\to \CS^2$ that intertwine the data of $\alpha_V$. The
category $\Catgeom$ is evidently abelian.

\medskip

Let $\Rg$ be ind-object of $\Sph$, corresponding under the equivalence
$\Rep(\cG)\simeq \Sph$ to $R_\cG$. A typical example of an object of $\Catgeom$
is obtained by setting
$$\CS:=\CS^1\star \Rg$$
for $\CS^1\in \Pervgrb$, where the Hecke isomorphisms come from
the canonical isomorphisms
$$\Rg\star \CV\simeq \uV\otimes \Rg.$$

As in the case of $\Catg$, the category $\Catgeom$ is naturally acted
on by the group $\cG$.

\medskip

We say that an objet of $\Catgeom$ is finitely generated if it admits a surjection
from an object of the form $\CS^1\star \Rg$ with $\CS^1\in \Pervgr$. This
condition is equivalent to the fact that the functor of $Hom$ from this object commutes
with direct sums.

\begin{conj} \label{Noetherianness of Hecke}
A sub-object of a finitely generated object of $\Catgeom$ is finitely generated.
\end{conj}

We will denote by $\Catgeomk$ (resp., $\Catgeomni$, $\Catgeomi$) a
version of the above category, where $\CS$ is assumed to be an object
of the ind-completion of the corresponding category $\Pervgrk$
(resp., $\Pervgrni$, $\Pervgri$).

As we shall see shortly, a particular case of \conjref{Noetherianness of Hecke},
concerning $\Pervgrni$, follows easily from \thmref{geometric steinberg}.

\medskip

We introduce a graded version $\Catgeomd$ of $\Catgeom$ analogously to
the definition of $\Catgd$: its objects are pairs $(\tCS,\{\alpha_V\})$,
where $\tCS$ is a $\cLambda$-graded object of $\Pervgrb$, and the maps
$\alpha_V$ preserve the gradings on both sides. Similarly, we introduce
the categories $\Catgeomdk$, $\Catgeomdni$, $\Catgeomdi$.

All of these categories are acted on naturally by the normalizer of
$\cT$ in $\cG$.

\medskip

Let $\tRg$ denote the same thing as $\Rg$, where we regard it as graded
via the right-action of $\cG$ on $R_\cG$. A typical example of an object of
$\Catgeomd$ is obtained by taking $\CS^1\star \tRg\{\cmu\}$ for $\CS^1\in \Pervgrb$,
where $\{\cmu\}$ denotes the shift of the grading functor.

In what follows we will state the results explicitly for $\Catgeomd$
and its versions; the transcription to the case of $\Catgeom$ is
straightforward.

\sssec{} \label{char of simples on gr}

Consider now the category $\Catgeomdni$. Combining \thmref{ABG} with
\thmref{AG graded} we obtain:

\begin{thm}  \label{ABG for small}
The category $\Catgeomdni$ is equivalent to $\tu_\ell\modo_0$.
\end{thm}

In particular, we obtain:

\begin{cor}  \hfill  \label{irr in Hecke}

\smallskip

\noindent{\em (1)}
The irreducibles in $\Catgeomdni$ are of the form $\CL^w\star \tRg\{\cmu\}$
for some $w\in W$ and $\cmu\in \cLambda$.

\smallskip

\noindent{\em (2)}
Every finitely generated object in $\Catgeomdni$ is Artinian.

\end{cor}

We will now give a geometric proof of this fact, using \thmref{geometric steinberg}.

\begin{proof}

Let us first see that any map $\tCS\to \CL^w\star \tRg\{\cmu\}$ is necessarily
a surjection. (This would imply that $\CL^w\star \tRg\{\cmu\}$ is irreducible.)

With no restriction of generality, we can assume that $\tCS$ has the form
$\CS'\star \tRg\{\cmu'\}$ for some $\CS'\in \Pervgrni$, $\cmu'\in \cLambda$.
Moreover, we can assume that $\CS'$ is itself irreducible. Then, by
\thmref{geometric steinberg}, $\CS'\simeq \CL^{w'}\star \IC_{\clambda,\Gr}$
for some $w'\in W$, $\clambda\in \cLambda^{pos}$. Hence,
\begin{equation}  \label{geom induced}
\CS'\star \tRg\{\cmu'\}\simeq \underset{\cnu}\oplus\,
\CL^{w'}\star \tRg\{\cmu'+\cnu\} \otimes \uV^\clambda(\cnu).
\end{equation}

Again, by \thmref{geometric steinberg}, the existence of a map
$$\CL^{w'}\star \tRg\{\cmu'\} \otimes \uV^\clambda\to \CL^w\star \tRg\{\cmu\}$$
forces that $w'=w$ and the map factors through the direct summand corresponding
to $\cmu'+\cnu=\cmu$ and a linear functional $\uV^\clambda(\cnu)\to \BC$.
Such a map is manifestly surjective.

\medskip

The same argument shows that any irreducible object of $\Catgeomdni$ admits
a map from some $\CL^w\star \tRg\{\cmu\}$. This establishes the first point
of the corollary.

\medskip

To prove the second point, it suffices to show that the objects of the form
$\CS'\star \tRg\{\cmu\}$, $\CS'\in \Pervgrni$ have finite lengths. For
that we can assume that $\CS'$ is irreducible, and our assertion follows
from \eqref{geom induced}.

\end{proof}

Let $\Catgeomdni_{Art}$ denote the subcategory of Artinian (or, equivalently,
finitely generated) objects of $\Catgeomdni$. By \thmref{ABG for small}, it
is equivalent to the category $\tu_\ell\mod_0$. Hence, it also carries
a duality functor, denoted $\BD$.

Explicitly, this functor is determined by the fact that it is exact;
$$\BD(\CS'\star \tRg\{\cmu\})\simeq \BD(\CS')\star \tRg\{\cmu\};$$
it is extended to the entire $\Catgeomdni_{Art}$ by the procedure
described in \secref{comp with duality}.

This functor goes over to the functor $N\mapsto N^\vee$ on $\tu_\ell\mod_0$,
since the equivalence of \thmref{ABG} transforms contragredient duality to Verdier
duality.

\section{Some results on $\Pervgr$}

\ssec{Proof of \thmref{geometric steinberg}}

\sssec{}

It clear that point (1) of the theorem implies point (2).
Indeed, for any $\IC_{w\cdot \cnu,\Gr}\in \Pervgri$ define $\CJ\subset \CI$ to be
the subset of simple roots, for which $w(\alpha_\imath)\in \cLambda^{pos}$. Define
$$\clambda':=\cnu-\underset{\jmath\in \CJ}\Sigma\,  \langle \alpha_\jmath,  \cnu\rangle
\cdot\check\omega_\jmath -\underset{\imath\in \CI-\CJ}\Sigma\, (\langle \alpha_\imath,  \cnu\rangle-1)
\cdot\check\omega_\imath,$$
where $\check\omega_\imath$ are (some choice of) fundamental coweights.

Then $w\cdot \clambda'$ is left-minimal with respect to $W$, and
$$\IC_{w\cdot \clambda',\Gr}\simeq \CL^w\star \delta_{\check\eta,\Gr},$$
where $\check\eta$ is a co-character of $Z(G)$ and $\delta_{\check\eta,\Gr}$
is the $\delta$-function at the corresponding point of $\Gr$.

By point (1),
$$\CL^w\star (\IC_{\cnu-\clambda',\Gr}\star \delta_{-\check\eta,\Gr})
\simeq \IC_{w\cdot \cnu,\Gr}.$$

\sssec{}

The assertion of point (1) is equivalent to the fact that
$\on{End}(\CL^w\star \IC_{\cmu,\Gr})\simeq \BC$. By adjunction, this is
equivalent to the fact that if $\CL^w=\IC_{w\cdot\clambda,\Gr}$, then
$$\on{Hom}\left(\IC_{w\cdot\clambda,\Gr},\IC_{w\cdot\clambda,\Gr}\star \IC_{\cmu,\Gr}
\star \IC_{-w_0(\cmu),\Gr}\right)\simeq \BC.$$
By decomposing $\IC_{\cmu,\Gr}\star \IC_{-w_0(\cmu),\Gr}$ as a sum of irreducibles,
we arrive to the conclusion that it is enough to show that
\begin{equation} \label{when Hom}
\on{Hom}\left(\IC_{w\cdot\cnu,\Gr},\IC_{w\cdot\clambda,\Gr}\star \IC_{\cmu,\Gr}\right)\neq 0
\Rightarrow \cnu=\clambda+\cmu.
\end{equation}
Note that this would automatically imply that
\begin{equation} \label{w=w'}
\on{Hom}\left(\IC_{w'\cdot\cnu,\Gr},\IC_{w\cdot\clambda,\Gr}\star \IC_{\cmu,\Gr}\right)=0
\text{ for } w'\neq w,
\end{equation}
a fact that will be used later on.

\sssec{}

We will establish \eqref{when Hom} by analyzing the convolution diagram.
First, we need to recall why the convolution functor
$\Pervgr\times \Sph\to \Pervgr$ is exact.

\medskip

Let
$$\Gr\star \Gr\simeq G((t))\underset{G[[t]]}\times \Gr$$
be the convolution diagram, which we think of as fibered over
$\Gr$ by means of projection to the first factor, with typical fiber
$\Gr$, which we think of as the second factor. We will denote
by $\pi$ the map $\Gr\star \Gr\to \Gr$ given by multiplication.
This ind-scheme is acted on by $G((t))$, and the map $\pi$ is
evidently $G((t))$-equivariant.

For $G[[t]]$-orbits $\Gr^{\cmu_1}, \Gr^{\cmu_2}\subset \Gr$ we will denote by
$\Gr^{\cmu_1}\star \Gr^{\cmu_2}$ the corresponding locally closed subset
in $\Gr\star \Gr$, which is fibered over $\Gr^{\cmu_1}$ with typical fiber
$\Gr^{\cmu_2}$. We will denote by $\Bigl(\Gr\star \Gr\Bigr)^\cmu$
(resp., $\Bigl(\Gr^{\cmu_1}\star \Gr^{\cmu_2}\Bigr)^\cmu$) the preimage
of $\Gr^\cmu$ in $\Gr\star \Gr$ (resp., $\Gr^{\cmu_1}\star \Gr^{\cmu_2}$)
under the map $\pi$.

We recall the following dimension estimates:
\begin{equation} \label{dim est gr}
\dim(\Gr^\cmu)=\langle 2\rho,\cmu\rangle, \,\,
\dim\left(\Bigl(\Gr^{\cmu_1}\star \Gr^{\cmu_2}\Bigr)^\cmu\right)=
\langle \rho,\cmu_1+\cmu_2+\cmu\rangle.
\end{equation}
Hence, the dimsnion of the fibers of the map
$$\pi^\cmu_{\cmu_1,\cmu_2}:
\Bigl(\Gr^{\cmu_1}\star \Gr^{\cmu_2}\Bigr)^\cmu\to \Gr^\cmu$$
is $\leq \langle \rho,\cmu_1+\cmu_2-\cmu\rangle$.

\medskip

For $\CS\in \Pervgr$, $\CF\in \Sph$, we will denote by
$\CS\tboxtimes \CF$ the corresponding perverse sheaf on
$\Gr\star \Gr$, and by definition,
$$\CS\star \CF=\pi_!(\CS\tboxtimes \CF).$$

To prove the exactness of convolution, by Verdier duality, it suffices
to show that the *-restriction of $\CS\tboxtimes \CF$ to every
$\Bigl(\Gr^{\cmu_1}\star \Gr^{\cmu_2}\Bigr)^\cmu$ lives in the
cohomological degrees $-\leq \langle \rho,\cmu_1+\cmu_2-\cmu\rangle$.

\medskip

It is evident that the *-restriction of $\CS\tboxtimes \CF$ to
$\Gr^{\cmu_1}\star \Gr^{\cmu_2}$ lives in the cohomological degrees
$\leq 0$. Moreover, $\CS\tboxtimes \CF|_{\Gr^{\cmu_1}\star \Gr^{\cmu_2}}$
is a pull-back of a complex on the base $\Gr^{\cmu_1}$.

Observe now that the constant sheaf
$\uBC_{\Bigl(\Gr^{\cmu_1}\star \Gr^{\cmu_2}\Bigr)^\cmu}$, thought of
as a complex on $\Gr^{\cmu_1}\star \Gr^{\cmu_2}$, is universally
locally acyclic (ULA) with respect to $\pi^\cmu_{\cmu_1,\cmu_2}$.
Indeed, it is $G[[t]]$-equivariant, and this group acts transitively
on the base. Hence,
$$\CS\tboxtimes \CF|_{\Bigl(\Gr^{\cmu_1}\star \Gr^{\cmu_2}\Bigr)^\cmu}\simeq
\CS\tboxtimes \CF|_{\Gr^{\cmu_1}\star \Gr^{\cmu_2}}\otimes
\uBC_{\Bigl(\Gr^{\cmu_1}\star \Gr^{\cmu_2}\Bigr)^\cmu}$$ lives in the
cohomological degrees
$$\leq -\on{codim}\left(\Bigl(\Gr^{\cmu_1}\star \Gr^{\cmu_2}\Bigr)^\cmu,
\Gr^{\cmu_1}\star \Gr^{\cmu_2}\right)\leq -\langle \rho,\cmu_1+\cmu_2-\cmu\rangle,$$
which is what we needed.

\medskip

The same argument proves also the following. Let $\CY\subset \Gr^\cmu$ be
a locally closed subscheme. In order for $\CS\star \CF|_\CY$ to have a non-zero
$0$-th perverse cohomology, it is necessary that there exist $\cmu^1$ and $\cmu^2$,
such that the fibers of the map
$$\on{supp}\Bigl(h^0(\CS|_{\Gr^{\cmu_1}})\tboxtimes h^0(\CF|_{\Gr^{\cmu_2}})\Bigr)
\cap (\pi^\cmu_{\cmu_1,\cmu_2})^{-1}(\CY)\to \CY$$
are of dimension equal to $\langle \rho,\cmu_1+\cmu_2-\cmu\rangle$,
i.e., saturating the upper bound given above.

\sssec{}

Thus, to prove \eqref{when Hom}, we must show that the fibers of the map
$$\pi^{-1}\Bigl((I\cdot (w\cdot\cnu))\Bigr)
\cap \Bigl((I\cdot (w\cdot\clambda))\star \Gr^\cmu\Bigr)\to \Bigl(I\cdot (w\cdot\cnu)\Bigr).$$
have dimensions
$<\langle \rho, \clambda+\cmu-\cnu\rangle$ unless $\cnu=\clambda+\cmu$.
(In the latter case the map in question is clearly one-to-one.)

Consider the orbit of the group $\on{Ad}_{w\cdot w_0}\left(N((t))\right)$ in $\Gr$ passing
through the point $w\cdot \cnu$. Its preimage in
$$\Bigl(I\cdot (w\cdot \clambda)\Bigl)\star \Gr^\cmu\subset \Gr\star \Gr$$
is the union over the parameters $\cnu'$ of schemes
$$\Bigl(\bigl(\on{Ad}_{w\cdot w_0}N((t))\cdot w(\cnu')\bigr)
\cap \bigl(I\cdot (w\cdot \clambda)\bigr)\Bigr)
\star \Bigl(\bigl(\on{Ad}_{w\cdot w_0}N((t))
\cdot (w\cdot (\cnu-\cnu'))\bigr)\cap \Gr^\cmu\Bigr),$$
each of which is fibered over
\begin{equation}  \label{intersection of orbits}
\Bigl(N((t))\cdot w_0(\cnu')\Bigr)\cap \Bigl(\on{Ad}_{(w\cdot w_0)^{-1}}(I)\cdot
w_0(\clambda)\Bigr)\subset \Gr^\clambda,
\end{equation}
with the typical fiber $\on{Ad}_{w\cdot w_0}N((t))
\cdot (w\cdot (\cnu-\cnu'))\cap \Gr^\cmu$.

\medskip

Since the intersection $\on{Ad}_{w\cdot w_0}N((t))\cdot (w\cdot \cnu)
\cap I\cdot (w\cdot \cnu)$ consists of a single point, namely, $w\cdot \clambda$,
the preimage of this point in $\Bigl(I\cdot (w\cdot \clambda)\Bigl)\star \Gr^\cmu$
injects into the variety \eqref{intersection of orbits}. The dimension of this variety
is a priori $\leq$ than
$$\langle \rho,-\cnu'+\clambda\rangle=\langle \rho,\clambda+\cmu-\cnu\rangle-
\langle \rho, \cmu-\cnu+\cnu'\rangle.$$

The non-emptiness condition on $\Bigl(\on{Ad}_{w\cdot w_0}N((t))
\cdot (w\cdot (\cnu-\cnu'))\Bigr)\cap \Gr^\cmu$ implies that $\langle \rho,\cmu-\cnu+\cnu'\rangle\geq 0$,
and the equality is achieved only for $\cnu-\cnu'=\cmu$. Hence, it is sufficient to prove
that the variety in \eqref{intersection of orbits} has dimension equal to
$\langle \rho,-\cnu'+\clambda\rangle$ only for $\cnu'=\clambda$.

\medskip

Note that the condition on $\clambda$ implies that
$$\on{Ad}_{(w\cdot w_0)^{-1}}(I)\cdot w_0(\clambda)\subset N^-[[t]]'\cdot w_0(\clambda),$$
where $N^-[[t]]'$ is the preimage under $N^-[[t]]\to N^-$ of $[N^-,N^-]\subset N^-$.
Let $\Psi_0$ be a non-degenerate character on $N^-((t))$ with conductor $0$.
Again, by the condition on $\clambda$,
$$\on{Ad}_{-w_0(\clambda)}(N^-[[t]]')\subset \on{ker}(\Psi_0).$$

Hence, the required assertion follows from the next result:

\begin{prop} \label{non-degenerate intersection}
The intersection
$$\Bigl(N((t))\cdot \cmu\Bigr)\cap \Bigl(\on{ker}(\Psi_0)\cdot 1_{\Gr}\Bigr)\subset
\Bigr(N((t))\cdot \cmu\Bigl)\cap \Bigl(N^-((t))\cdot 1_{\Gr}\Bigr)\subset \Gr$$
has dimension $<\langle \rho,\cmu\rangle$ unless $\cmu=0$.
\end{prop}

\sssec{Proof of \propref{non-degenerate intersection}}

The assertion of the proposition is equivalent to the fact that the character $\Psi_0$
is non-trivial on every connected component of the intersection
$$\Bigl(N((t))\cdot \cmu\Bigr)\cap \Bigl(N^-((t))\cdot 1_{\Gr}\Bigr).$$

Let $\clambda\in \cLambda^+$ be a large. Then, then it is well-known that
$$\Bigl(N((t))\cdot (\clambda+\cmu)\Bigr)\cap \Bigl(N^-((t))\cdot \clambda\Bigr)
=\Gr^{\clambda+\cmu}\cap \Bigl(N^-((t))\cdot \clambda\Bigr).$$

Hence, it is sufficient to show that a character $\Psi_\clambda$ on
$N^-((t))$ with conductor $\clambda$ is non-constant on every connected
component of the intersection $\Gr^{\clambda+\cmu}\cap (N^-((t))\cdot \clambda)$.
But the latter readily follows from the (top cohomology part) of the Casselman-Shalika
formula, \cite{FGV}, Sect. 7.1.7.

\ssec{The baby Whittaker category}

\sssec{}

Let us denote by $I^-$ the group $\on{Ad}_{w_0}(I^0)\subset G[[t]]$,
and let $\psi:I^-\to \BG_a$ a non-degenerate character.
We introduce the (baby Whittaker) category $\Pervgrw$ as the that of $(I^-,\psi)$-equivariant
perverse sheaves on $\Gr$. \footnote{The term "baby Whittaker" refers to the fact that we are
imposing equivariance with respect to $I^-$, rather than with respect to the group
ind-scheme $N^-((t))$.}

If $\CS\in \Pervgrw$ and $\clambda\in \cLambda^+$, both *- and !- restrictions of
$\CS|_{\Gr^\clambda}$ can be non-zero only if $\clambda$ is regular.
Moreover, in this case, these restrictions are supported on the $I^-$-orbit
of the point $w_0(\clambda)\in \Gr$.

For $\clambda\in \cLambda^+$ we will denote by $\IC^\psi_{\clambda,\Gr}$ the
Goresky-MacPherson extention of the $(I^-,\psi)$-character sheaf on
the $I^-$-orbit of $w_0(\clambda+\crho')\in \Gr$. It is easy to see that the
$\IC^\psi_{\clambda,\Gr}$'s are the irreducibles of $\Pervgrw$.

\medskip

We will denote $\IC^\psi_{0,\Gr}$ simply by $\IC_{\Gr}^\psi$. It is easy to
see that $\IC^\psi_{0,\Gr}$ is in fact a {\it clean} extension of the
corresponding character sheaf on $I^-\cdot w_0(\crho')$. Indeed, all
$G[[t]]$-orbits in the closure of $\Gr^{\crho'}$ correspond to non-regular
coweights.

Using the same argument as in the proof of \thmref{geometric steinberg},
one shows:

\begin{thm}  \label{Cass-Shal for I}
$\IC_{\Gr}^\psi\star \IC_{\clambda,\Gr}\simeq \IC^\psi_{\clambda,\Gr}$.
\end{thm}

The same argument as in \cite{FGV}, Sect. 6, implies then the following:

\begin{cor}  \label{cleanness of Whit}  \hfill

\smallskip

\noindent{\em (1)}
The category $\Pervgrw$ is semi-simple and equivalent to $\Sph$ by means
of $\CF\mapsto \IC^\psi_{\Gr}\star \CF$.

\smallskip

\noindent{\em (2)}
$\IC^\psi_{\clambda,\Gr}$ equals both the !- and *-extension of the
corresponding character sheaf on $I^-\cdot w_0(\clambda+\crho')$.

\end{cor}

\sssec{}

Let $\Dgrw$ denote the $(I^-,\psi)$-equivariant derived category on $\Gr$.
The forgetful functor $\Pervgrw\to \Pervgr$ admits natural left and right adjoints,
denoted $\on{Av}_{!,I^-,\psi}$ and $\on{Av}_{*,I^-,\psi}$, respectively.

\begin{prop} \label{psi averaging}
The functors $\on{Av}_{!,I^-,\psi}[-\dim(\fn)]$ and $\on{Av}_{*,I^-,\psi}[\dim(\fn)]$,
when restricted to $\Dgrni$, are isomorphic. Both these functors are exact.
\end{prop}

\begin{proof}

Note that the character $\psi$ factors through the map $I^-\to N^-$;
we will denote by the same symbol $\psi$ the resulting character
of $N^-$. Let $\psi_{N^-}$ denote the corresponding character sheaf
on $N^-$.

It is clear that the restrictions of $\on{Av}_{!,I^-,\psi}$ and
$\on{Av}_{*,I^-,\psi}$ to $\Dgrone$ are the functors
$$\CS\mapsto \psi_{N^-}\starshriek \CS[\dim(\fn)] \text{ and }
\CS\mapsto \psi_{N^-}\starstar \CS[-\dim(\fn)],$$
respectively, where $\starshriek$ and $\starstar$ are the two
convolution functors
$${\mathsf D}(G)\times \Dgrone\to \Dgrone.$$

In particular, we have a map of functors
$$\on{Av}_{!,I^-,\psi}[-\dim(\fn)]|_{\Dgrone}\to \on{Av}_{*,I^-,\psi}[\dim(\fn)]|_{\Dgrone}.$$

\medskip

To show that the above map of functors is an isomorphism, when restricted further
to $\Dgrni$, it is sufficient to prove the corresponding fact for $\Dgri$.
 \footnote{A more efficient proof of this fact is given in \cite{BBM}}

Let $\Pervfdw$ (resp., $\Dfdw$) be the corresponding
$(N^-,\psi)$-equivariant category on $G/B$. We will denote by
$\psi_{G/B}$ its only irreducible, i.e., the
clean extension of the $(N^-,\psi)$-character perverse sheaf on $N^-\cdot 1_{G/B}$.

For $\CS\in \Dgri$ we have:
$$\psi_{N^-}\starshriek \CS\simeq \psi_{G/B}\starshriek \CS \text{ and }
\psi_{N^-}\starstar \CS\simeq \psi_{G/B}\starstar \CS,$$
but the map
$$\psi_{G/B}\starshriek \CS\to \psi_{G/B}\starstar \CS$$ is an isomorphism,
since the convolution map $\pi_I:G[[t]]\underset{I}\times \Gr\to \Gr$ is proper.

\medskip

The exactness assertion follows as well, since the functor
$\CS\mapsto \psi_{N^-}\starshriek \CS$ is left-exact, and
$\CS\mapsto \psi_{N^-}\starstar \CS$ is right-exact.

\end{proof}

Henceforth, we will denote the functor
$$\on{Av}_{!,I^-,\psi}[-\dim(\fn)]|_{\Dgrni}\simeq \on{Av}_{*,I^-,\psi}[\dim(\fn)]|_{\Dgrni}$$
simply by $\on{Av}_{I^-,\psi}$.

\sssec{Partial integrability}

We say that an object of $\Pervgrni$ (resp., $\Pervgri$, $\Pervflni$, $\Pervfli$) is
{\it partially integrable} if it admits a filtration, such that each subquotient is
equivariant with respect to some parahoric, contained in $G[[t]]$, and strictly containing $I$.
(The latter condition is equivalent to demanding that this subquotient
is equivariant with respect to some subminimal parabolic $P_\imath\subset G
\subset G[[t]]$.)  Let us denote the resulting Serre subcategories
by $^{PI}\Pervgrni$ (resp., $^{PI}\Pervgri$, $^{PI}\Pervflni$, $^{PI}\Pervfli$).

Note that an irreducible $\IC_{w\cdot\clambda,\Gr}\in \Pervgrni$ is non-partially
integrable if and only if $w=w_0$. Similarly, $\IC_{w,G/B}\in \Pervfd$ is
partially integrable unless $w=1$.

Let $\fPervgrni$ (resp., $\fPervgri$, $\fPervflni$, $\fPervfli$) be the resulting
quotient abelian category of $\Pervgrni$ (resp., $\Pervgri$, $\Pervflni$, $\Pervfli$),
and let $\fDgrni$ (resp., $\fDgri$, $\fDflni$, $\fDfli$) be the corresponding quotient
triangulated category.

The convolution functor descends to functors
$$\fDflni\times \Dgri\to {}\fDgrni \text{ and } \fDfli\times \Dgri\to {}\fDgri.$$

\medskip

\begin{prop}  \label{whit and part int}
The functor $$\on{Av}_{I^-,\psi}:\Pervgrni\to \Pervgrw$$ factors through $\fPervgrni$,
and the resulting functor $\fPervgrni\to \Pervgrw$ is exact and faithful.
\end{prop}

\begin{proof}

The fact that $\on{Av}_{I^-,\psi}$ annihilates all partially integrable objects
follows from the observation that the direct image of $\psi_{G/B}$ to
any partial flag variety $G/P_\imath$ is zero.

The fact that $\fPervgrni\to \Pervgrw$ is exact follows from the exactness
statement of \propref{psi averaging}. To show that it is faithful, it is
enough to prove the corresponding fact for $\fPervgri$. We argue as follows:

\medskip

Let $\on{Av}_{!,I^0}:\Dgr\to \Dgrni$ (resp., $\Dfd\to \Dfdn$) be the functor, left adjoint to
the tautological embedding. Let us denote by
$\Xi$ the object
\begin{equation}  \label{defn Xi}
\Xi:=\on{Av}_{!,I^0}[-\dim(\fn)](\psi_{G/B})\in \Dfdn.
\end{equation}

We have, tautologically:

\begin{lem}   \label{convolution with Xi}
The composition
$$\Pervgri\overset{\on{Av}_{I^-,\psi}}\longrightarrow\Pervgrw\to
\Pervgr \overset{\on{Av}_{!,I^0}[-\dim(\fn)]}\longrightarrow \Pervgrni$$
is isomorphic to the convolution functor
$$\CF\mapsto \Xi\star \CF.$$
\end{lem}

\medskip

It is known that $\Xi$ is the longest indecomposable projective in $\Pervfdn$, and
it admits two filtrations: one whose subquotients are the standard objects
$j_{!,w}$, $w\in W$, and another, whose subquotients are the costandard objects
$j_{*,w}$.

Note, however, that the arrows $j_{*,w}\to j_{*,1}$ and $j_{!,1}\to j_{!,w}$ become isomorphisms
on $\fPervfd\subset \fPervfli$. Hence, the image of $\Xi$ in $\fPervfdn$ is isomorphic to the
extension of $|W|$-many copies of $\delta_{1,G/B}$. Hence, the convolution with
$\Xi$, viewed as a functor $\fPervgri\to \fPervgrni$, is faithful.

\end{proof}

\medskip

\noindent{\it Remark.}
One can strengthen \propref{whit and part int} and prove the following more precise
assertion:

\medskip

Let $h^0$ be the algebra of functions on the scheme-theoretic preimage of
$0$ under $\fh^*\to \fh^*/W$. It is known that $h_0$ is isomorphic to the algebra
of endomorphisms of $\Xi$.

For an abelian category $\CC$ we will denote by
$\CC\otimes h^0$ the category of objects of $\CC$, endowed with an action of $h^0$.

Then the category $\fPervgrni$ is equivalent to $\Pervgrw\otimes h^0\simeq \Sph\otimes h^0$.

\ssec{Cosocles of costandard objects}

\sssec{}

In this subsection we will prove the following assertion:

\begin{prop} \label{cosocle of 0 costandard} \hfill

\smallskip

\noindent{\em(1)}
For a regular dominant element $\clambda\in \cLambda$, the cosocle of
$\CW^{*,\clambda}\in \Pervgri$ is simple and is isomorphic to
$\IC_{w_0\cdot \clambda,\Gr}$.

\smallskip

\noindent{\em(2)}
The kernel of $\CW^{*,\clambda}\to \IC_{w_0\cdot \clambda,\Gr}$
is partially integrable.
\end{prop}

\begin{proof}

First, we claim that if we have a surjection from $\CW^{*,\clambda}$
to an irreducible $\CS$, then this $\CS$ must be non-partially
integrable. Suppose the contrary, and let $\imath\in \CI$ be such that
$\CS$ is equivariant with respect to the corresponding sub-minimal
parahoric. Then the convolution $j_{!,s_\imath}\star \CS$ lives in
the cohomological degree $+1$. However,
$j_{s_\imath,!}\star \CW^{*,\clambda}\simeq \CW^{*,s_\imath\cdot \lambda}$
is still perverse. Hence,
$\on{Hom}_{\Dgri}(j_{s_\imath,!}\star \CW^{*,\clambda},j_{!,s_\imath}\star \CS)=0$,
which is a contradiction, since the convolution with $j_{s_\imath,!}$
is an auto-equivalence of $\Dgri$.

\medskip

To finish the proof of the proposition, it suffices to show that
$\IC_{w_0\cdot \clambda,\Gr}$ is the only
non-partially integrable irreducible that appears in the Jordan-H\"older
series of $\CW^{*,\clambda}$. Since the natural map
$$\CW^{*,\clambda}\simeq j_{*,\lambda}\star \delta_{1,\Gr}\to
j_{!,w_0}\star j_{*,\lambda}\star \delta_{1,\Gr}\simeq \CW^{*,w_0\cdot \clambda}$$
becomes an isomorphism in $\fPervgri$, by
\propref{whit and part int}, it suffices to show that the map
$$\on{Av}_{I^-,\psi}(\CW^{*,w_0\cdot\clambda})\to
\on{Av}_{I^-,\psi}(\IC^{w_0\cdot \clambda,\Gr})$$
is an isomorphism.

By \propref{psi averaging}, it would be sufficient to show that
$\on{Av}_{I^-,\psi}(\CW^{*,w_0\cdot\clambda})$ is an irreducible object of
$\Pervgrw$. However, evidently,
$$\on{Av}_{*,I^-,\psi}(\CW^{*,w_0\cdot\clambda})[\dim(\fn)]$$
is the *-extension of the corresponding character sheaf on $I^-\cdot \clambda$.
Hence, we are done by \thmref{Cass-Shal for I}(2).

\end{proof}

\sssec{}

We will now prove the following:

\begin{prop}  \label{cosocle of another costandard}
If $\clambda\in \cLambda^+$ is large, the object $\CW^{*,w_0\cdot \clambda}$ admits
$\IC_{\clambda-2\crho,\Gr}$ as a quotient.
\end{prop}

The rest of this subsection is devoted to the proof of this result. Let
$\on{Av}_{!,G[[t]]/I}$ be the functor $\Pervgri\to \Pervgrg\simeq \Sph$ left
adjoint to the forgetful functor $\Pervgrg\to \Pervgri$. Note that since
$G[[t]]/I=G/B$ is compact, the corresponding right adjoint $\on{Av}_{*,G[[t]]/I}$
is isomorphic to $\on{Av}_{!,G[[t]]/I}[2\dim(\fn)]$.

For a regular $\clambda\in \cLambda^+$, let us denote by
$\on{emb}^\clambda_{!,\Gr},\on{emb}_{*,\Gr}^\clambda$ the natural functors
$\Dfd\to \Dgr$, along with its $I$- and $I^0$-equivariant versions.
Evidently, these functors commute with $\on{Av}_{!,G[[t]]/I}$ in the natural sense.
Therefore,
\begin{equation}  \label{ext by 0 as averaging}
\on{Av}_{!,G[[t]]/I}(\CW^{*,w_0\cdot \clambda})\simeq
\on{emb}_{*,\Gr}^\clambda(\ul{\BC}{}_{G/B}[2\dim(\fn)]).
\end{equation}
In particular, we obtain that the object
$\on{emb}_{*,\Gr}^\clambda(\IC_{G/B})$ lives in the cohomological
degrees $\leq \dim(\fn)$. Therefore,
$\on{Hom}_{\Pervgri}(\CW^{*,w_0\cdot \clambda},\IC_{\clambda-2\crho,\Gr})$
identifies with
$$\on{Hom}_{\Pervgrg}\Bigl(h^{\dim(\fn)}\left(\on{emb}_{*,\Gr}^\clambda(\IC_{G/B})\right),
\IC_{\clambda-2\crho,\Gr}\Bigr).$$

\medskip

Thus, we have to show that the top=$\dim(\fn)$--degree cohomology of
$\on{emb}_{*,\Gr}^\clambda(\IC_{G/B})$ has
a quotient (or, which in this case is the same, a direct summand) isomorphic
to $\IC_{\clambda-2\crho,\Gr}$. Set $\cmu=\clambda-2\crho$.

Consider the cohomology
$$H_c\Bigl(N((t))\cdot w_0(\cmu'),
\on{emb}_{!,\Gr}^\clambda(\IC_{G/B})|_{N((t))\cdot w_0(\cmu')}\Bigr),$$
where, as usual, we regard $N((t))\cdot w_0(\cmu)$ as a sub-indscheme
in $\Gr$. By \cite{MV} (and duality) it would suffice to show that the above
cohomology in degree $-\dim(\fn)-\langle 2\rho, \cmu'\rangle$ is
$1$-dimensional if $\cmu'=\cmu$, and is $0$ for $\cmu<\cmu'\leq \clambda$.

By base change, the above cohomology can be rewritten as
$$H_c^{\langle 2\rho, \clambda-\cmu'\rangle-\dim(\fn)}\Bigl(
\left(N((t))\cdot w_0(\cmu')\right)\cap \Gr^\clambda,\ul{\BC}\Bigr).$$

Since $\clambda$ was assumed large, the intersection
$\left(N((t))\cdot w_0(\cmu')\right)\cap \Gr^\clambda$ equals
$$\left(N((t))\cdot w_0(\cmu')\right)\cap \left(N^-((t))\cdot w_0(\clambda)\right))\simeq
\left(N((t))\cdot 1_{\Gr}\right)\cap \left(N^-((t))\cdot w_0(\clambda-\cmu')\right)).$$
Hence, our assertion follows from \corref{cohomology of open fibers}.

\begin{cor}
For $\clambda$ large the map $\CW^{*,\clambda}\to \IC_{w_0\cdot \clambda,\Gr}$
lifts to a map $\CW^{*,\clambda}\to \CW^{!,w_0\cdot \crho'}\star \IC_{\clambda-\crho',\Gr}$.
\end{cor}

\begin{proof}
The existence of the map in question is equivalent, by adjunction, to the existence of a map
$j_{*,w_0\cdot (-w_0(\crho'))}\star \CW^{*,\clambda}\to \IC_{\clambda-\crho',\Gr}$. Note that
$-w_0(\crho')=2\crho-\crho'$. Hence, the assertion follows from the above proposition, since:

$$j_{*,w_0\cdot (2\crho-\crho')}\star \CW^{*,\clambda}\simeq
j_{*,w_0\cdot (2\crho-\crho')}\star j_{*,\clambda}\star \delta_{1,\Gr}\simeq
j_{*,w_0\cdot (\clambda+2\crho-\crho')}\star \delta_{1,\Gr}\simeq \CW^{*,\clambda+2\crho-\crho'}.$$
\end{proof}

\section{A study of baby Verma and co-Verma modules}

\ssec{Baby co-Verma modules via $\bU_\ell$}

\sssec{}

Let $\tbM^\lambda$ be the object of $\Catgd$, corresponding to $\tM^\lambda$. Our
present goal is to describe it explicitly. First, we will describe $\tbM^\lambda$ as an
object of $\bU_\ell\modo$. By definition,
$\bM^\lambda=\underset{\check\mu\in \cLambda}\oplus\, \bM^\lambda_{\check\mu}$, where
each  $\bM^\lambda_{\check\mu}$ is given by
$$\Ind^{\bU_\ell}_{\tu_\ell}(\bC^{-\cmu}\otimes \tM^\lambda)\simeq
\Ind^{\bU_\ell}_{\tu_\ell}(\tM^{\lambda-\phi_\ell(\cmu)}).$$
Hence, it sufficient to describe the modules of the form $\Ind^{\bU_\ell}_{\tu_\ell}(\tM^\lambda)$.

\medskip

By construction,
$$\Ind^{\bU_\ell}_{\tu_\ell}(\tM^\lambda)\simeq \Ind^{\bU_\ell}_{\tb^-_\ell}(\bC^\lambda),$$
which, in turn, is isomorphic to
\begin{equation} \label{two step induction}
\Ind^{\bU_\ell}_{\bB_\ell^-}\left(\Ind^{\bB_\ell^-}_{\tb^-_\ell}(\bC^\lambda)\right)\simeq
\Ind^{\bU_\ell}_{\bB_\ell^-}\left(\bC^\lambda\otimes \Ind^{\bB_\ell^-}_{\tb^-_\ell}(\BC)\right).
\end{equation}

By \propref{quant Frob for B},
$\Ind^{\bB_\ell^-}_{\tb^-_\ell}(\BC)\simeq \on{Fr}_{B^-}
(\CO_{\check B^-/\check T})$.

\begin{propconstr} \label{functions on B/T}
As a $\check B^-$-module, $\CO_{\check B^-/\check T}$ is isomorphic to the
direct limit
$$\underset{\check\lambda\in \cLambda^+}{\underset{\longrightarrow}{lim}}\,
\Res^{\cG}_{\cB^-}\Bigl((V^{\check\lambda})^*\Bigr)\otimes \fl^{\check\lambda},$$
$\fl^\clambda$ denotes highest weight line of $V^\clambda$, regarded as a
$1$-dimensional representation of $\cT$ (and, hence, also of $\cB^-$).
\end{propconstr}

\begin{proof}

By adjunction, to specify a map of $\cB^-$-modules
\begin{equation} \label{mat coeff}
\Res^{\cG}_{\cB^-}\Bigl((V^{\check\lambda})^*\Bigr)
\otimes \fl^{\check\lambda}\to \CO_{\check B^-/\check T},
\end{equation}
is equivalent to giving a map $(\uV^\clambda)^*\to (\fl^\clambda)^*$,
compatible with the $T$-action. The latter corresponds to the
natural embedding of $\fl^\clambda$ into $\uV^\clambda$.

\medskip

To define the inductive system, we choose a compatible system
of isomorphisms $\fl^\clambda\otimes \fl^\cmu\simeq \fl^{\clambda+\cmu}$.
Such a system fixes as the maps $V^{\check\lambda}\otimes V^{\check\mu}\to
V^{\check\lambda+\check\mu}$ (which are otherwise defined up to a scalar).

Suppose that $\check\mu\in \cLambda^+$ is another dominant weight of $\check G$.
We define the map
$$\Res^\cG_{\cB^-}\Bigl((V^{\check\lambda})^*\Bigr)\otimes \fl^{\check\lambda}\to
\Res^\cG_{\cB^-}\Bigl((V^{\check\lambda+\check\mu})^*\Bigr)
\otimes \fl^{\check\lambda+\check\mu}$$ as
the composition
\begin{equation} \label{Cartan mult}
\Res^\cG_{\cB^-}\Bigl((V^{\check\lambda})^*\Bigr)\otimes \fl^{\check\lambda}\to
\Res^\cG_{\cB^-}\Bigl((V^{\check\lambda}\otimes V^{\check\mu})^*\Bigr)\otimes
\fl^\cmu \otimes \fl^\clambda \to
\Res^\cG_{\cB^-}\Bigl((V^{\check\lambda+\check\mu})^*\Bigr)\otimes
\fl^{\clambda+\cmu},
\end{equation}
where the first arrow comes from the map of $\check B^-$-modules
$(\fl^\cmu)^*\to \Res^\cG_{\cB^-}\Bigl((V^{\check\mu})^*\Bigr)$, and
the second arrow comes from the map $(V^{\check\lambda}\otimes V^{\check\mu})^*\to
(V^{\check\lambda+\check\mu})^*$. These maps define the inductive system
stated in the Proposition-Construction.

By construction, the map of \eqref{Cartan mult} is compatible with the maps of
\eqref{mat coeff} for $\check\lambda$ and $\check\mu$. Hence, the resulting
inductive limit maps to $\CO_{\check B^-/\check T}$. The fact that this map is an
isomorphism is an easy verification.

\end{proof}

\sssec{}

>From the above Proposition we obtain the following description of $\bM^\lambda_\cnu$:

\begin{cor}  \label{descr of baby Verma}
Choose a trivialization of the $\cT$-torsor given by $\{\fl^\clambda\}$. Then
$$\bM^\lambda_\cnu\simeq \underset{\check\lambda\in \cLambda^+}{\underset{\longrightarrow}{lim}}\,
\on{Fr}\Bigl((V^\clambda)^*\Bigr)\otimes \bW^{\lambda+\phi_\ell(\clambda-\cnu)},$$
where the maps in
the inductive system are given by
\begin{align*}
&\on{Fr}(V^\clambda)^*\otimes \bW^{\lambda+\phi_\ell(\clambda-\cnu)}\to
\on{Fr}(V^\clambda)^*\otimes \on{Fr}(V^\cmu)^*\otimes \bW^{\lambda+
\phi_\ell(\clambda+\cmu-\cnu)}\to \\
&\on{Fr}(V^{\clambda+\cmu})^*\otimes
\bW^{\lambda+\phi_\ell(\clambda+\cmu-\cnu)},
\end{align*}
where the first arrow comes from the canonical map
$$\on{Fr}(V^\cmu)\otimes \bW^{\lambda+\phi_\ell(\clambda-\cnu)}\to
\bW^{\lambda+\phi_\ell(\clambda+\cmu-\cnu)}.$$
\end{cor}

\sssec{Hecke property}

Let us now describe how the Hecke isomorphisms
$$\on{Fr}(V)\otimes \bM^\lambda_\cmu\to \underset{\cnu\in \cLambda^+}\oplus\,
\bM^\lambda_{\cmu-\cnu}\otimes \uV(\cnu)$$
look like in terms of the identification of \corref{descr of baby Verma}.

For a coweight $\clambda\in \cLambda^+$, large with respect to the weights of $V$,
we have a canonical isomorphism of $\cG$-modules
$$V\otimes V^\clambda\simeq \underset{\cnu}\oplus\,
V^{\clambda+\cnu}\otimes \uV(\cnu).$$
Hence, we obtain a map of $\cB^-$-modules
$$V\to \underset{\cnu}\oplus\, V^{\clambda+\cnu} \otimes (V^\clambda)^* \otimes
\uV(\cnu)\to (V^\clambda)^*\otimes \fl^{\clambda+\cnu}\otimes \uV(\cnu).$$

Applying the functor $\Ind^{\bU_\ell}_{\bB^-_\ell}$, for
$\lambda\in \Lambda^+$,  we obtain a map in $\bU_\ell\mod$:
$$\on{Fr}(V)\otimes \bW^\lambda \to \underset{\cnu}\oplus\,
\on{Fr}\Bigl((V^\clambda)^*\Bigr)\otimes
\bW^{\lambda+\phi_\ell(\clambda+\cnu)}\otimes \uV(\cnu).$$

\begin{prop}   \label{Hecke for baby Verma}
The Hecke morphisms for $\tbM^\lambda$ are equal in terms of the inductive
system to
\begin{align*}
&\on{Fr}(V)\otimes \on{Fr}\Bigl((V^{\clambda'})^*\Bigr)\otimes
\bW^{\lambda+\phi_\ell(\clambda'-\cmu)} \to \\
&\underset{\cnu}\oplus\,
\on{Fr}\Bigl((V^{\clambda'})^*\otimes
(V^\clambda)^*\Bigr)\otimes \bW^{\lambda+\phi_\ell(\clambda+\clambda'-\cmu+\cnu)}
\otimes \uV(\cnu) \to \\
&\underset{\cnu}\oplus\,
\on{Fr}\Bigl((V^{\clambda+\clambda'})^*\Bigr)
\otimes \bW^{\lambda+\phi_\ell(\clambda+\clambda'-\cmu+\cnu)} \otimes \uV(\cnu).
\end{align*}
\end{prop}

\begin{proof}

By the construction of the isomorphism in \corref{descr of baby Verma}, it suffices to
show that the isomorphism
$$\Res^{\cG}_{\cB^-}(V)\otimes
\CO_{\cB^-/\cT}\to \underset{\nu}\oplus\, \CO_{\cB^-/\cT}\otimes
\uV(\cnu)$$
looks in terms of the identification given by \propconstrref{functions on B/T}
as a system of morphisms
\begin{align*}
&\Res^{\cG}_{\cB^-}(V)\otimes \Bigl(V^{\clambda'})^*\Bigr)\otimes \fl^{\clambda'}\to
\underset{\cnu}\oplus\,
\Res^\cG_{\cB^-}\Bigl((V^{\clambda'})^*\otimes
(V^\clambda)^*\Bigr)\otimes \fl^{\clambda+\clambda'+\cnu}
\otimes \uV(\cnu)\to \\
&\to \underset{\cnu}\oplus\,
\Res^\cG_{\cB^-}\Bigl((V^{\clambda+\clambda'})^*\Bigr)
\otimes \fl^{\clambda+\clambda'+\cnu} \otimes \uV(\cnu).
\end{align*}

The latter is a straightforward verification.

\end{proof}

\sssec{Baby co-Verma as a quotient}   \label{Baby as a quotient}

Let us briefly discuss another realization of $\tM^\lambda$ (or, equivalently, $\tbM^\lambda$)
in terms of the big quantum group.

For an element $\cmu\in \cLambda^*$, let $\overset{\circ}{\uV}{}^\cmu$ be the hyperplane
in $\uV^\cmu$ orthogonal to $\fl^{-\cmu}\subset (\uV^\cmu)^*$. This subspace
is preserved by $\cB^-$, and in particular, it admits a well-defined weight decomposition
with respect to $\cT$.

For $\lambda\in \Lambda^+$ consider the canonical map of $\bU_\ell$-modules:
$\on{Fr}(V^\cmu)\otimes \bW^{\lambda}\to \bW^{\lambda+\phi_\ell(\cmu)}$
After the restriction to $\tu_\ell$, it gives rise to a map
$$\underset{\cnu}\oplus\, \bC^\cnu\otimes \on{Res}^{\bU_\ell}_{\tu_\ell}(\bW^\lambda)\otimes
\uV^\cmu(\cnu)\to
\on{Res}^{\bU_\ell}_{\tu_\ell}(\bW^{\lambda+\phi_\ell(\cmu)}).$$

For $\clambda\in \cLambda^+$ consider the canonical map
$\Res^{\bU_\ell}_{\tu_\ell}(\bW^{\lambda+\phi_\ell(\mu)})\to \tM^{\lambda+\phi_\ell(\mu)}$.

\begin{prop}  \label{Baby Verma as a quotient}
The composition
$$\underset{\cnu}\oplus\, \bC^\cnu\otimes \on{Res}^{\bU_\ell}_{\tu_\ell}(\bW^{\lambda})\otimes
\overset{\circ}{\uV}{}^\cmu(\cnu)
\to \on{Res}^{\bU_\ell}_{\tu_\ell}(\bW^{\lambda+\phi_\ell(\cmu)})
\to \tM^{\lambda+\phi_\ell(\mu)}$$
vanishes. For a fixed $\lambda$ and all sufficiently large $\cmu$ the complex
\begin{equation} \label{surj system 1}
\underset{\cnu}\oplus\, \bC^\cnu\otimes
\on{Res}^{\bU_\ell}_{\tu_\ell}(\bW^{\lambda+\phi_\ell(\clambda)})\otimes
\overset{\circ}{V}{}^\cmu(\cnu)\to
\on{Res}^{\bU_\ell}_{\tu_\ell}(\bW^{\lambda+\phi_\ell(\cmu+\clambda)})
\to \tM^{\lambda+\phi_\ell(\mu+\clambda)}\to 0
\end{equation}
is exact for all sufficiently dominant $\clambda$.
\end{prop}

\begin{proof}

The first assertion of the proposition is evident. To prove the second one we proceed
as follows.
Let $\BC_{\cB^-}$ be the sky-scraper coherent sheaf at the point $\cB^-$ in the
flag variety $\cG/\cB^-$. It admits a left resolution of the form
$$0\to \CP_{\dim(\cG/\cB^-)+1}\to \CP_{\dim(\cG/\cB^-)}\to
\CP_{\dim(\cG/\cB^-)-1}\to...\to \CP_{1}\to \CP_0\to \BC_{\cB^-}\to 0,$$
where $\CP_0\simeq \CO_{\cG/\cB^-}$, and the sheaves $\CP_{i}$ for
$i=1,...,\dim(\cG/\cB^-)$ are isomorphic to  $\CO(-\cmu_i)\otimes V^i$
for $\cmu_i\in \cLambda^+$; $V^i$ are some vector spaces. Moreover,
the weight $\cmu_1$ may be chosen arbitrarily large; and the vector space
$V^1$ surjects by construction onto
$\overset{\circ}{V}{}^{\cmu_1}$.

By pulling back this complex from $\cG/\cB^-$ to $\cG$, it gives rise to a complex
\begin{align*}
&P\to \bC^{-\cmu_{\dim(\cG/\cB^-)}}\otimes\Res^{\cG}_{\cB^-}(R_\cG)\otimes
V^{\dim(\cG/\cB^-)}\to
...\to \bC^{-\cmu_i}\otimes\Res^{\cG}_{\cB^-}(R_\cG)\otimes V^i \to...  \\
&...\to\bC^{-\cmu_1}\otimes\Res^{\cG}_{\cB^-}(R_\cG)\otimes V^1\to
\Res^{\cG}_{\cB^-}(R_\cG)\to R_{\cB^-}\to 0
\end{align*}
of $\cB^-$-modules, where $R_{\cB^-}$ denotes the regular representation of $\cB^-$.
By construction, the arising element in
$\on{Ext}_{\cB^-}^{\dim(\cG/\cB^-)+1}(R_{\cB^-},P)$ vanishes.

Let us tensor this complex with the $\bB^-_\ell$-module $\bC^{\lambda+\phi_\ell(\cmu_1+\clambda)}$,
where $\clambda$ is such that all the weights of the form
$\lambda+\phi_\ell(\cmu_1+\clambda-\cmu_i)$ become dominant. Then,
\begin{align*}
&\on{R}^i\Ind^{\bU_\ell}_{\bB^-_\ell}
\Bigl(\bC^{\lambda+\phi_\ell(\clambda+\cmu_1-\cmu_i)}\otimes
\on{Fr}_{\cB^-}\left(\Res^{\cG}_{\cB^-}(R_\cG)\right)\Bigr)
\simeq \\
&\simeq\on{Fr}(R_\cG)\otimes \on{R}^i\Ind^{\bU_\ell}_{\bB^-_\ell}
\left(\bC^{\lambda+\phi_\ell(\clambda+\cmu_1-\cmu_i)}\right)=0
\end{align*}
for $i>0$.

Hence, we obtain that the sequence of $\bU_\ell$-modules
$$\on{Fr}(R_\cG)\otimes \bW^{\lambda+\phi_\ell(\clambda)}\otimes
\overset{\circ}V{}^{\mu_1}\to
\on{Fr}(R_\cG)\otimes \bW^{\lambda+\phi_\ell(\cmu_1+\clambda)}\to
\Ind^{\bU_\ell}_{\bB^-_\ell}(\bC^{\lambda+\phi_\ell(\cmu_1+\clambda)}\otimes R_{\cB^-})\to 0$$
is exact. However, the above sequence of maps is obtained from \eqref{surj system 1}
for $\cmu=\cmu_1$ by applying the functor $\Ind^{\bU_\ell}_{\tu_\ell}\circ \Res^{\tu_\ell}_{\fu_\ell}$,
which is exact and faithful.

\end{proof}

\sssec{The case of twisted baby co-Verma modules}

For $w\in W$ let $\sF_w,w_\ell,w_{\cG}$ be as in \secref{Weyl group action}.
>From \secref{further Weyl} we obtain the following description of
the object $^w\tbM^{w(\lambda)}$ of the category $\Catgd$, corresponding to
$^w\tM{}^{w(\lambda)}$:

\begin{cor}
As an object of  $\bU_\ell\modo$, $^w\tbM^{w(\lambda)}_\cmu$ is isomorphic to
$\tbM{}^{\lambda}_{w(\cmu)}$. The Hecke property morphisms
$$\on{Fr}(V)\otimes {}^w\tbM^{w(\lambda)}_\cmu\to
\underset{\cnu}\oplus\, {}^w\tbM{}^{w(\lambda)}_\cmu\otimes \uV(\nu)$$
are obtained from those of $\tbM{}^{\lambda}_{w(\cmu)}$ by applying the
element $w_{\cG}:V(\nu)\to V(w(\nu))$.
\end{cor}

In addition, we have an analogue of \propref{Baby Verma as a quotient}.
Let $^w\overset{\circ}{\uV}{}^{\cmu}$ be the subspace of
$\uV^\cmu$ obtained by translating $\overset{\circ}{\uV}{}^\cmu$
by means of $w_\cG$.

\begin{cor}
We have a complex
$$\underset{\cnu}\oplus\,
\bC^\cnu\otimes\on{Res}^{\bU_\ell}_{\tu_\ell}(\bW^{\lambda+\phi_\ell(\clambda)})\otimes
{}^w\overset{\circ}{\uV}{}^{\cmu}(\cnu) \to
\on{Res}^{\bU_\ell}_{\tu_\ell}(\bW^{\lambda+\phi_\ell(\cmu+\clambda)})
\to {}^w\tM^{w(\lambda+\phi_\ell(\mu+\clambda))}\to 0,$$
which is exact when for a fixed $\lambda$, the coweights $\cmu$
and $\clambda$ are large enough.
\end{cor}

\sssec{The non-graded version}

For $\lambda\in \Lambda$ recall that $M^\lambda$ denotes the restriction of $\tM{}^\lambda$
to $\fu_\ell$ (the small, non-graded quantum group).
Let $\bM^\lambda$ be
the corresponding object of $\Catg$. From \corref{descr of baby Verma} and
we obtain a description of $\bM^\lambda$ as an object of
$\bU_\ell\mod$. Namely,
\begin{equation} \label{real baby Verma}
\bM^\lambda\simeq \underset{\cmu}\oplus\,
\underset{\check\lambda\in \cLambda^+}{\underset{\longrightarrow}{lim}}\,
\on{Fr}\Bigl((V^\clambda)^*\Bigr)\otimes \bW^{\lambda+\phi_\ell(\clambda)+\phi_\ell(\cmu)}.
\end{equation}

The Hecke isomorphisms for $\bM^\lambda$ are given by disregarding
the grading in the isomorphisms for $\tbM{}^\lambda$, given by
\propref{Hecke for baby Verma}.

\medskip

In addition, we can realize $M^\lambda$ as a quotient of modules, restricted
from $\bU_\ell$, using \propref{Baby Verma as a quotient}:
\begin{equation} \label{non-graded baby as a quotient}
M^\lambda\simeq \on{coker}\left(\on{Res}^{\bU_\ell}_{\tu_\ell}(\bW^{\lambda+\phi_\ell(\clambda)})
\otimes \overset{\circ}{\uV}{}^\cmu \to
\on{Res}^{\bU_\ell}_{\tu_\ell}(\bW^{\lambda+\phi_\ell(\cmu+\clambda)})\right).
\end{equation}

\sssec{$\check G$-action on baby co-Verma modules}

By \secref{further Weyl}, to any $\bg\in \cG$ we can attach a module
$^\bg M^\lambda\in \fu_\ell\mod$. Explicitly, $^\bg M^\lambda$ corresponds
to the object $^\bg\bM^\lambda\in \Catg$, where the latter is obtained from
$\bM^\lambda$ by modifying the Hecke isomorphism using
$\bg$ acting on $\underline{V}$ for $V\in \cG\mod$.  Equivalently,
$^\bg M^\lambda$ can be realized as
\begin{equation} \label{twisted non-graded baby as a quotient}
\on{coker}\left(\on{Res}^{\bU_\ell}_{\tu_\ell}(\bW^{\lambda+\phi_\ell(\clambda)})
\otimes {}^\bg\overset{\circ}{\uV}{}^\cmu \to
\on{Res}^{\bU_\ell}_{\tu_\ell}(\bW^{\lambda+\phi_\ell(\cmu+\clambda)})\right),
\end{equation}
where $^\bg\overset{\circ}{\uV}{}^\cmu$ is the $\bg$-translate of
$\overset{\circ}{\uV}{}^\cmu$ inside $\uV^\cmu$.

By \secref{further Weyl} if $\bg$ belongs to the normalizer of the torus
$\cT\subset \cG$, $^\bg M^\lambda$ is isomorphic to
$\Res^{\tu_\ell}_{\fu_\ell}({}^w\tM^{w(\lambda)})$ for the corresponding
$w\in W$.

\begin{prop}
If $\bg\in \cB^-$, then $^\bg M^\lambda$ is isomorphic to $M^\lambda$.
For $\lambda=0$ the above condition is "if and only if".
\end{prop}

\begin{proof}

The description of $^\bg M^\lambda$ given by \eqref{twisted non-graded baby as a quotient}
makes it clear that if $\bg\in \cB^-$, then $^\bg M^\lambda\simeq M^\lambda$.
To show the inclusion in the opposite direction we argue as follows:

It is easy to see that the subset of elements of $\cG$, which stabilize the isomorphism
class of $M^\lambda$ is a Zariski-closed subgroup of $\cG$. Hence, we must
show that this subgroup does not contain any parabolic strictly containing $\cB^-$.
Therefore, it is enough
to show that the none of the modules $\Res^{\tu_\ell}_{\fu_\ell}({}^w\tM^0)$
for $w\neq 1$ is isomorphic to $M^0$. This is equivalent to $\tM{}^0$ being
non-isomorphic to $^w\tM^{\phi_\ell(\cmu)}$ for $\cmu\in \cLambda$, $1\neq w\in W$.

Note that the socle of $^w\tM^{\phi_\ell(\cmu)}$ is isomorphic to $\bC^{\phi_\ell(\cmu)}$. Hence,
if $^w\tM^{\phi_\ell(\cmu)}\simeq \tM^0$, then $\cmu=0$. However, it is clear that
$^w\tM^0$ is non-isomorphic to $\tM^0$, because, for example,
$-\phi_\ell(2\check\rho)+2\rho$, which appears as a weight of $\tM^0$, is not among
the weights of $^w\tM^0$.

\end{proof}

One can show that the condition of the proposition is in fact "if
and only if" for any $\lambda$ belonging to the regular block.
This is because, as we shall see later, baby co-Verma modules with
the same $w$, but different parameters $\lambda$, can be obtained
from one another by (invertible) convolution functors.

\ssec{Baby co-Verma modules via perverse sheaves on the affine
Grassmannian}

\sssec{}

For an element $\tw\in W_{aff}$, let $\lambda\in \Lambda$ be the
corresponding weight in the regular block. That is, if $\tw=w\cdot \clambda$, then
$\lambda=\phi_\ell(\clambda)+w^{-1}(\rho)-\rho$.

Let $\tCM^{\tw}=\underset{\cmu}\oplus\, \tCM^\tw_\cmu$ be the object of
$\Catgeomdni$, corresponding to the object $\tbM^\lambda\in \Catgd$. By
\corref{descr of baby Verma},
as an object of $\Pervgrnib$,
\begin{equation} \label{Baby Verma on Grassmannian}
\tCM^{w\cdot \clambda}_\cmu:=\underset{\clambda' \in \cLambda^+}{\underset{\longrightarrow}{lim}}\,
\CW^{*,w\cdot (\clambda+\clambda'-\cmu)}\star \IC_{-w_0(\clambda'),\Gr}.
\end{equation}

The maps in this inductive system come from the canonical maps
\begin{align*}
&\CW^{*,w\cdot \cmu'}\star \IC_{\clambda',\Gr}\simeq
j_{*,w\cdot \cmu'}\star \IC_{\clambda',\Gr}\to
j_{*,w\cdot \cmu'}\star \CW^{*\clambda'}\simeq \\
&j_{*,w\cdot \cmu'}\star j_{*,\clambda'}\star \delta_{1,\Gr}\simeq
j_{*,w\cdot (\cmu'+\clambda')}\star \delta_{1,\Gr}\simeq
\CW^{*,w\cdot (\cmu'+\clambda')}.
\end{align*}

The Hecke morphisms $\tCM^\tw_\cmu\star \CV\to \underset{\cnu}\oplus\, \uV(\cnu)
\otimes \tCM^\tw_{\cmu-\cnu}$ for $V\in \cG\mod$, are given by
translating the morphisms of \propref{Hecke for baby Verma} into the
geometric context. Namely, let $\clambda$ be a weight large, compared to $V$.
Then the sought-for morphism is
\begin{align*}
&\CW^{*,w\cdot (\clambda'-\cmu)}\star  \IC_{-w_0(\clambda'),\Gr}\star \CV
\to \\
&\CW^{*,w\cdot (\clambda'-\cmu)}\star
\left(\IC_{\clambda,\Gr}\star \CV\right)\star \left(\IC_{-w_0(\clambda),\Gr}\star
\IC_{-w_0(\clambda'),\Gr}\right) \simeq \\
&\CW^{*,w\cdot (\clambda'-\cmu)}\star
\left(\underset{\cnu}\oplus\, \uV(\cnu)\otimes \IC_{\clambda+\cnu,\Gr}\right)
\star \left(\IC_{-w_0(\clambda),\Gr}\star \IC_{-w_0(\clambda'),\Gr}\right) \to \\
&\underset{\cnu}\oplus\, \uV(\cnu)\otimes
\CW^{*,w(\clambda+\clambda'-\cmu+\cnu)}\star
\left(\IC_{-w_0(\lambda),\Gr}\star \IC_{-w_0(\clambda'),\Gr} \right)\to \\
&\underset{\cnu}\oplus\, \uV(\cnu)\otimes  \CW^{*,w(\clambda+\clambda'-\cmu+\cnu)}\star
\IC_{-w_0(\clambda+\clambda'),\Gr}.
\end{align*}

\medskip

Evidently, if $\tw=w\cdot \lambda$ is such that for some $w'\in W$, $l(w')+l(w)=l(w'\cdot w)$,
we have
\begin{equation} \label{Baby Vermas from one another}
j_{!,w'}\star \tCM^{\tw}\simeq \tCM^{w'\cdot \tw},
\end{equation}
that it, the objects  $\tCM^{\tw}$ for different $\tw$ are obtained from one-another
by convolution.

Note also that for $\clambda,\cmu\in \cLambda$ with $\cmu$ dominant and $\clambda$
dominant and regular,
$$l(w\cdot \cmu\cdot w^{-1})+l(w\cdot \lambda)=
l(w\cdot (\cmu+\clambda)).$$ Hence, we obtain:

\begin{cor}   \label{lambda invariance}
For $\cmu\in \cLambda^+$ there are canonical isomorphisms
$$j_{*,w\cdot \cmu\cdot w^{-1}}\star \tCM^{w\cdot \clambda}\simeq
\tCM^{w\cdot (\clambda+\cmu)}\simeq
\tCM^{w\cdot \clambda}\{-\cmu\},$$
respecting the Hecke isomorphisms.
\end{cor}

Assume now that $w\in W$, $\clambda,\cmu\in \cLambda$ are such that
$w\cdot (\clambda+\cmu)$ is right $W$-minimal, i.e.
$\CW^{*,w\cdot (\clambda+\cmu)}$ is well-defined. By
\eqref{Baby Verma on Grassmannian}, we have a map
\begin{equation} \label{map of costandard to baby}
\CW^{*,w\cdot (\clambda+\cmu)} \to \tCM^{w\cdot \clambda}_{-\cmu},
\end{equation}
such that for $\cmu'\in \cLambda^+$, the diagram
$$
\CD
j_{*,w^{-1}\cdot\cmu'\cdot w}\star \CW^{*,w\cdot (\clambda+\cmu)} @>>>
j_{*,w^{-1}\cdot\cmu'\cdot w}\star \tCM^{w\cdot \clambda}_{-\cmu} \\
@V{\sim}VV @V{\sim}VV \\
\CW^{*,w\cdot (\clambda+\cmu+\cmu')}  @>>> \tCM^{w\cdot \clambda}_{-\cmu-\cmu'}
\endCD
$$
commutes.

Convolving \eqref{map of costandard to baby} on the right with $\IC_{\mu',\Gr}$
we obtain the map $$\CW^{*,w\cdot (\clambda+\cmu)}\star \IC_{\cmu',\Gr}\to
\tCM^{w\cdot \clambda}_{-\cmu}\star \IC_{\cmu',\Gr}\simeq
\underset{\cnu}\oplus\, \uV^{\cmu'}(\cnu)\otimes
\tCM^{w\cdot \clambda}_{-\cnu-\cmu}.$$

The above description of the Hecke morphisms implies also the following:

\begin{cor} \label{Miura} \hfill

\smallskip

\noindent {\em (1)} The diagram
$$
\CD
\CW^{*,w\cdot (\clambda+\cmu)}\star \IC_{\cmu',\Gr}  @>>>
\CW^{*,w\cdot (\cmu+\cmu'+\clambda)} \\
@VVV   @VVV \\
\underset{\cnu}\oplus\, V^{\cmu'}(\cnu)\otimes \tCM^{w\cdot \clambda}_{-\cnu-\cmu}
@>>> \tCM^{w\cdot \clambda}_{-\cmu'-\cmu},
\endCD
$$
commutes, where the bottom horizontal arrow is the projection on the direct summand,
corresponding to $\cnu=\cmu'$.

\smallskip

\noindent {\em (2)}
The object $\tCM^{w\cdot \clambda}\in \Catgeomdni$ is universal with respect to the properties
that (a) it satisfies \corref{lambda invariance}, (b) it receives a map
as in \eqref{map of costandard to baby} for some $\cmu$, such that (a) and (b)
render the above diagram is commutative.

\end{cor}

If we put $w=1$ and $\cmu=-\clambda$, the map in \eqref{map of costandard to baby}
identifies with
\begin{equation} \label{0 socle map}
\delta_{1,\Gr}\to \tCM^\clambda_{-\clambda}.
\end{equation}
Thus, we obtain a characterization of $\tCM^\clambda$ in terms of $\delta_{1,\Gr}$.

\medskip

Finally, we note that the normalizer of the torus $\cT\subset \cG$
acts on $\Catgeomdni$ by self-equivalences, modifying the Hecke morphisms.
The functors, corresponding to elements of $\cT$ are (non-canonically) isomorphic
to identity. For $w\in W$ we will denote by $^w\tCM^{\tw}$ the
object of $\Catgeomdni$ obtain in this way from $\tCM^\tw$;
it corresponds to the object $^w\tM^{w(\lambda)}\in \tu_\ell\mod$.

\sssec{}

We will now list several facts about the objects $\tCM^{\tw}$, most of which are formal
consequences of the corresponding properties of $\tM{}^\lambda$, but  we
will give geometric proofs for completeness.

Let $\tw=w\cdot \cmu$ be an element of $W_{aff}$, and let $\clambda$ be such that
$\CL^w\simeq\IC_{w\cdot \clambda,\Gr}$, in particular, $w\cdot\clambda$ is restricted.
Then we have:

\begin{prop} \label{socle of geometric baby Verma}
The socle of $\tCM^\tw$ is isomorphic to $\CL^w\star \tRg\{\clambda-\cmu\}$.
\end{prop}

\begin{proof}

By \secref{char of simples on gr}, every irreducible in $\Catgeomdni$
is of the form $\CL^{w'}\star \tRg\{\cmu'\}$ for some $w'\in W$ and $\cmu'\in \cLambda$.
Suppose that such an irreducible maps to $\tCM^{w\cdot \cmu}$. By adjunction, this
means that we have a map
$$\CL^{w'}\to \CW^{*,w\cdot (\clambda'+\cmu+\cmu')}\star \IC_{-w_0(\clambda'),\Gr}$$
in $\Pervgrni$ for some $\clambda'\in \cLambda^+$.

The latter can be rewritten as an element in
$$\on{Hom}(\CL^{w'}\star \IC_{\clambda',\Gr},
\CW^{*,w\cdot (\clambda'+\cmu+\cmu')}).$$
By \thmref{geometric steinberg}, and taking into account that the socle of
$\CW^{*,w\cdot (\clambda'+\cmu+\cmu')}$ is isomorphic to
$\IC_{w\cdot (\clambda'+\cmu+\cmu'),\Gr}$, this implies $w'=w$ and $\cmu'=\clambda-\cmu$.

We also obtain that the above $\on{Hom}$ is $1$-dimensional. I.e.,
$\CL^w\star \tRg\{\clambda-\cmu\}$ is the only irreducible that can map to $\tCM^\tw$,
and it appears in the socle with multiplicity $1$.

\end{proof}


\begin{prop}   \label{geometric cosocles}  \hfill

\smallskip

\noindent{\em (1)}
The object $\tCM^1\in \Catgeomdi$ is finitely generated.
\footnote{Here and in the sequel, the superscript "1" in $\tCM^1$ stands for
the unit element in $W_{aff}$.}
Its cosocle is isomorphic to $\CL^{w_0}\star \tRg\{\crho'\}$.
Moreover,
all the constituents in $\on{ker}(\tCM^1\to \CL^{w_0}\star\tRg\{\crho'\})$ are partially integrable.

\smallskip

\noindent{\em (2)}
There exists a surjection $\tCM^{w_0}\twoheadrightarrow \tRg\{2\crho\}$.
\end{prop}

Using \eqref{Baby Vermas from one another}, from point (1) of the proposition we obtain:

\begin{cor}
Every $\tCM^\tw$ is finitely generated (as an object of $\Catgeomdi$).
\end{cor}

Before giving the proof of the proposition, we need to introduce the following construction.
Consider a $\cLambda$-graded object of $\Pervgrnib$ given by $$\tRg{}':=\underset{\cmu}\oplus\, \underset{\clambda\in \cLambda^+}{\underset{\longrightarrow}{lim}}\,
\IC_{\clambda-\cmu,\Gr}\star \IC_{-w_0(\clambda),\Gr},$$
where the maps in the inductive system are given as follows. If $\clambda'=\clambda+\cnu$,
$\cnu\in \cLambda^+$,
\begin{align*}
&\IC_{\clambda-\cmu,\Gr}\star \IC_{-w_0(\clambda),\Gr}\to
\Bigl(\IC_{\clambda-\cmu,\Gr}\star \IC_{\cnu,\Gr}\Bigr)
\star \Bigl(\IC_{-w_0(\cnu),\Gr}\star \IC_{-w_0(\clambda),\Gr}\Bigr)\to \\
&\IC_{\clambda-\cmu+\cnu,\Gr}\star \IC_{-w_0(\cnu+\clambda),\Gr}.
\end{align*}

\begin{propconstr}  \label{new reg}
The object $\tRg{}'$ is a Hecke eigen-sheaf,  and as such, it
is canonically isomorphic to $\tRg$.
\end{propconstr}

\begin{proof}

Since all the appearing perverse sheaves are spherical, we can work
in the tensor category of $\Rep(\cG)$ instead of $\Pervgrni$. The Hecke
morphisms are given as follows. Let $\clambda\in \cLambda^+$ be large compared
to $V$. Then the sought-for map is the composition:
\begin{align*}
&V\otimes V^{\clambda'-\cmu}\otimes (V^{\clambda'})^*\to
V\otimes V^\clambda \otimes
V^{\clambda'-\cmu}\otimes (V^{\clambda'})^*\otimes (V^{\clambda})^*\to \\
&\underset{\cnu}\oplus\, \Bigl(V^{\clambda'-\cmu}\otimes V^{\clambda+\cnu}\Bigr)\otimes
\Bigl((V^{\clambda'})^*\otimes (V^{\clambda})^*\Bigr)\otimes \uV(\nu)\to
\underset{\cnu}\oplus\, V^{\clambda+\lambda'+\cnu-\cmu}\otimes
(V^{\clambda+\cnu})^*\otimes \uV(\nu).
\end{align*}

To see that $\tRg{}'$ is isomorphic to $\tRg$, it is enough to notice that
$$\on{Hom}\Bigl(V,\underset{\clambda\in \cLambda^+}{\underset{\longrightarrow}{lim}}\,
V^{\clambda-\cmu}\otimes (V^{\clambda})^*\Bigr)\simeq
\underset{\longrightarrow}{lim}\, \on{Hom}(V\otimes V^\clambda,V^{\clambda-\cmu}).$$
When $\clambda$ is large with respect to $V$, the latter inductive system stabilizes
to $(\uV)^*(\cmu)$.

\end{proof}

Now we are ready to prove \propref{geometric cosocles}.

\begin{proof}

First, we claim that $\tCM^1$ cannot map to any partially integrable object
of the category $\Catgeomdni$. Indeed, if $\tCS$ were partially integrable and we
had a non-zero map $\tCM^1\to \tCS$, we would have a non-zero map in $\Pervgri$:
$$\CW^{*,\lambda}\star \IC_{\cmu,\Gr}\to \CS'$$
for some $\clambda,\cmu\in \cLambda^+$ and $\CS'\in {}^{PI}\Pervgri$.
By adjunction we would then have a map
$$\CW^{*,\lambda}\to \CS'\star \IC_{-w_0(\cmu),\Gr}=:\CS'',$$
with $\CS''$ being also partially integrable. But the latter is impossible by
\propref{cosocle of 0 costandard}(1).

\medskip

Let us now construct a map $\tCM^1\to \CL^{w_0}\star \tRg\{\crho'\}$.
According to \propref{cosocle of 0 costandard}(1)  and \thmref{geometric steinberg},
for every dominant and regular $\clambda$ we have a canonical map
$$\CW^{*,\clambda}\to \IC_{w_0\cdot \clambda}\simeq \IC_{w_0\cdot \crho'}\star
\IC_{\clambda-\crho',\Gr}.$$
In addition, for $\cmu\in \cLambda^+$ the diagram
$$
\CD
\CW^{*,\clambda}\star \IC_{\cmu,\Gr} @>>> \CW^{*,\clambda+\cmu} \\
@VVV   @VVV \\
\IC_{w_0\cdot \crho'}\star (\IC_{\clambda-\crho',\Gr}\star \IC_{\cmu,\Gr})  @>>>
\IC_{w_0\cdot \crho'}\star \IC_{\clambda+\cmu-\crho',\Gr}
\endCD
$$
is easily seen to commute.

This defines the map of between the inductive systems:
$$\tCM^{1}_\cmu\simeq \underset{\clambda}{\underset{\longrightarrow}{lim}}\,
\CW^{*,\clambda-\cmu}\star \IC_{-w_0(\clambda),\Gr} \to
\underset{\clambda}{\underset{\longrightarrow}{lim}}\,
\CL^{w_0}\star (\IC_{\clambda-\cmu-\crho'}\star \IC_{-w_0(\clambda),\Gr}),$$
and the latter identifies with $\CL^{w_0}\star \tRg\{\crho'\}$, by
\propconstrref{new reg}.

It is straightforward to check that the above map respects the Hecke morphisms,
i.e., we obtained the desired map in $\Catgeomdni$.  Moreover, from
\propref{cosocle of 0 costandard}(2) it follows that the kernel of the map
$\tCM^{1}\to \CL^{w_0}\star \tRg\{\crho'\}$ is partially integrable.

\medskip

To finish the proof of the first part of the proposition, it remains to show
that the map
\begin{equation} \label{from W to M}
\CW^{*,\lambda}\star \tRg\{\clambda\}\to \tCM^1
\end{equation}
is surjective for some (and, in fact, every) regular $\clambda$. By construction, the composition
$$\CW^{*,\lambda}\star \tRg\{\clambda\}\to \tCM^1\to \CL^{w_0}\star \tRg\{\crho'\}$$
is surjective. Hence, by the above, the cokernel of \eqref{from W to M} is partially
integrable, and hence, is zero.

\medskip

To prove the second assertion of the proposition, recall from
\propref{cosocle of another costandard} that for $\clambda$ large we have a map
\begin{equation} \label{aux 1}
\CW^{*,w_0\cdot \clambda}\to \IC_{\clambda-2\crho},
\end{equation}
defined up to a scalar. Moreover, from the construction of this map one deduces that the
square
$$
\CD
\CW^{*,w_0\cdot \clambda}\star \IC_{\cmu,\Gr} @>>>
\CW^{*,w_0\cdot (\clambda+\cmu)} \\
@VVV   @VVV  \\
\IC_{\clambda-2\crho} \star \IC_{\cmu,\Gr}  @>>> \IC_{\clambda+\cmu-2\crho,\Gr}
\endCD
$$
commutes (up to a scalar). We can normalize the maps in \eqref{aux 1} to make
such diagrams commutative.

This gives us a map of inductive systems
$$\tCM^{w_0}_\cmu\simeq \underset{\clambda}{\underset{\longrightarrow}{lim}}\,
\CW^{*,w_0\cdot (\clambda-\cmu)}\star \IC_{-w_0(\clambda),\Gr} \to
\underset{\clambda}{\underset{\longrightarrow}{lim}}\,
\IC_{\clambda-\cmu-2\crho}\star \IC_{-w_0(\clambda),\Gr},$$
and the latter identifies with $\tRg\{2\crho\}$.

\end{proof}

\sssec{A dual description}

Recall that over the small quantum group, the baby Verma modules
can be expressed through the baby co-Verma modules and a twist by
elements of the Weyl group. We would like to establish this fact
geometrically as well. By \eqref{Baby Vermas from one another}, it
suffices to consider the case of just $\tCM^1$.

\begin{prop}   \label{dual descr}
We have an isomorphism
$$\BD(\tCM^1)\simeq ({}^{w_0}\tCM^{w_0})\{2\crho\}.$$
\end{prop}

Since the convolution functors commute with Verdier duality, from
\eqref{Baby Vermas from one another},we obtain:

\begin{cor}  \label{duality on all baby}
$$\BD(\tCM^{w\cdot \cmu})\simeq {}^{w_0}\tCM^{w\cdot w_0\cdot (w_0(\cmu)+2\crho)}.$$
\end{cor}

Combining this with \propref{socle of geometric baby Verma}, we also obtain:
\begin{cor}
The cosocle of every $\tCM^{w\cdot \cmu}$ is simple and isomorphic to
$$\IC_{w\cdot w_0\cdot \clambda,\Gr}\star \tRg\{w_0(\clambda)-\cmu+2\crho\},$$
where $\clambda\in \cLambda^+$ is such that $w\cdot w_0\cdot \clambda$ is restricted.
\end{cor}

\begin{proof}

By \propref{geometric cosocles}(1) and \propref{socle of geometric baby Verma},
it is enough to construct a map
\begin{equation} \label{map to dual realization}
\tCM^1\to \BD\left(({}^{w_0}\tCM^{w_0})\{2\crho\}\right),
\end{equation}
such that the composition
\begin{equation} \label{composition to cosocle 1}
\tRg\to \tCM^1\to \BD\left(({}^{w_0}\tCM^{w_0})\{2\crho\}\right)
\end{equation}
equals (up to a scalar) the map, obtained by duality from \propref{geometric cosocles}(2),
and such that the composition
\begin{equation} \label{composition to cosocle 2}
\tCM^1\to \BD\left(({}^{w_0}\tCM^{w_0})\{2\crho\}\right)\to \CL^{w_0}\{\crho'\}\star \tRg
\end{equation}
equals the map of \propref{geometric cosocles}(1).

\medskip

By \corref{lambda invariance} and duality
$$j_{!,w_0(\clambda)}\star \BD\left(({}^{w_0}\tCM^{w_0})\{2\crho\}\right)\simeq
\BD\left(({}^{w_0}\tCM^{w_0})\{2\crho\}\right)\{w_0(\clambda)\}.$$
Hence, by \corref{Miura}, to construct a map as in \eqref{map to dual realization},
we must construct a map $\tRg\to \BD\left(({}^{w_0}\tCM^{w_0})\{2\crho\}\right)$
and check the commutativity of the corresponding diagram. By duality, the above
amounts to a map $({}^{w_0}\tCM^{w_0})\{2\crho\}\to \tRg$. By the definition of
the twisting functors, the existence of the latter map follows from
\propref{geometric cosocles}(2). We need to check the commutativity of the following
diagram:
$$
\CD
({}^{w_0}\tCM^{w_0})\{2\crho+\clambda\} @>{\sim}>>
j_{!,\clambda}\star \left(({}^{w_0}\tCM^{w_0})\{2\crho\}\right) @>>>
j_{!,\clambda}\star \tRg \\
@VVV &  & @VVV \\
\underset{\cnu}\oplus\, \uV^\clambda(\cnu)\otimes
({}^{w_0}\tCM^{w_0})\{2\crho+\cnu\} @>{\sim}>>
({}^{w_0}\tCM^{w_0})\{2\crho\}\star \IC_{\clambda,\Gr} @>>> \tRg\star \IC_{\clambda,\Gr}.
\endCD
$$
Recalling the definition of the arrows, we arrive to the following diagram, defined for
$\cmu$ large:
$$
\CD
j_{!,\clambda}\star \CW^{*,w_0\cdot \cmu} @>>> j_{!,\clambda}\star \IC_{\cmu-2\crho,\Gr} \\
@V{\sim}VV   @VVV \\
\CW^{*,w_0\cdot (\cmu+w_0(\clambda))} @>>> \IC_{\cmu+w_0(\clambda)-2\crho,\Gr},
\endCD
$$
where the right vertical arrow is the composition
$$j_{!,\clambda}\star \IC_{\cmu-2\crho,\Gr}\to \IC_{\clambda,\Gr}\star  \IC_{\cmu-2\crho,\Gr}\to
\IC_{\cmu+w_0(\clambda)-2\crho,\Gr},$$
where the second arrow is obtained by adjunction from
$$\IC_{\cmu-2\crho,\Gr}\to \IC_{-w_0(\clambda),\Gr}\star \IC_{\cmu+w_0(\clambda)-2\crho,\Gr}.$$
The commutativity of the latter diagram follows from the construction of the arrow
in \propref{cosocle of another costandard}.

\medskip

By construction, the condition on the composed map from \eqref{composition to cosocle 1}
is satisfied. It remains to verify the condition in \eqref{composition to cosocle 2}.
The latter amounts to showing that the arrow
$$\IC_{w_0\cdot (\cmu+w_0(\clambda)),\Gr} \to
\CW^{*,w_0\cdot (\cmu+w_0(\clambda))}
\simeq
j_{!,\clambda}\star \CW^{*,w_0\cdot \cmu} \to j_{!,\clambda}\star \IC_{\cmu-2\crho,\Gr}$$
equals (up to a scalar) the map
\begin{align*}
&\IC_{w_0\cdot (\cmu+w_0(\clambda)),\Gr} \simeq
\IC_{w_0\cdot\crho',\Gr}\star \IC_{\cmu+w_0(\clambda)-\crho',\Gr}\to \\
&\IC_{w_0\cdot\crho',\Gr}\star \left(\IC_{\clambda-\crho',\Gr}\star \IC_{\cmu-2\crho,\Gr}\right)\simeq
\IC_{w_0\cdot \clambda,\Gr}\star \IC_{\cmu-2\crho,\Gr}\to j_{!,\clambda}\star \IC_{\cmu-2\crho,\Gr},
\end{align*}
where the last arrow comes by duality from \propref{cosocle of 0 costandard},
and the second arrow is obtained by adjunction from
$$\IC_{\cmu+w_0(\clambda)-\crho',\Gr}\star
\IC_{-w_0(\clambda)+w_0(\crho'),\Gr}\to \IC_{\cmu-2\crho,\Gr}.$$

By construction, both these maps are non-zero. Now our assertion follows from the fact,
that $\IC_{w_0\cdot \clambda,\Gr}$ is the only non-partially integrable constituent of
$\CW^{!,\clambda}$, which implies that $\IC_{w_0\cdot (\cmu+w_0(\clambda)),\Gr}$
appears with multiplicity one in the Jordan-H\"older series of
$j_{!,\clambda}\star \IC_{\cmu-2\crho,\Gr}$.

\end{proof}

\sssec{Non-graded version and presentation as a quotient}

Our present goal is to prove geometrically that $\tCM^{\tw}$ can
be presented as a quotient, as in \secref{Baby as a quotient}. For that it will
be convenient to consider the corresponding non-graded version,
$\CM^{\tw}\in \Catgeomi$. If $\tw=w\cdot \clambda$ and $\tw'=w\cdot \clambda'$,
then, evidently, $\CM^\tw\simeq \CM^{\tw'}$.

For $\bg\in \cG$ we will denote by $^\bg\CM^\tw$
the corresponding twist of $\CM^\tw$; for $\bg=w\in W$ we recover the objects
$^w\CM^\tw$. We will denote by $^{\cG}\CM^\tw$ the universal family of
$^\bg\CM^\tw$ over $\CO_{\cG}$.

\begin{lem}   \label{B^- action on baby} \hfill

\noindent{\em (1)}
As an object of $\Pervgrib$, $\CM^\tw$ admits a unique action of the
algebraic group $\cB^-$, such that

\begin{itemize}

\item If $\tw=w\cdot \clambda$, in terms of \eqref{Baby Verma on Grassmannian},
the image of $\CW^{*,\clambda+\clambda'}$ in $\CM^\tw$ transforms according
to the $\cB^-$-character $-\clambda'$.

\item
The Hecke isomorphisms
$$\CM^\tw\star \CV\simeq \uV\otimes \CM^\tw$$
intertwine the action of $\cB^-$ on the left-hand side, obtained by
transport of structure and the diagonal action of $\cB^-$ on the
right-hand side.

\end{itemize}

\noindent{\em (2)}
$$\on{Hom}_{\cB^-}\Bigl(\bC^{-\cmu}\otimes \Res^{\cG}_{\cB^-}(V^{\clambda'}),
\CM^\tw\Bigr)\simeq \CW^{w\cdot (\clambda+\cmu)}\star \IC_{-w_0(\clambda'),\Gr},$$
if $w\cdot (\clambda+\cmu)\in W_{aff}$ is right $W$-minimal, and $0$ otherwise.

\end{lem}

Note that the action of $\cT\subset \cB^-$ on $\CM^\tw$ as an object of $\Pervgrib$
comes from the grading on $\tCM^\tw$.

\medskip

The first assertion of the lemma means that, as an object of $\Catgeomi$,
$\CM^\tw$ is $\cB^-$-equivariant, i.e., that the $\CO_\cG$-family $^{\cG}\CM^\tw$
acquires a $\cB^-$-action, covering that on $\CO_\cG$. Alternatively,
a structure of a $\cB^-$-equivariant object on some $\CN\in \Catgeomi$
is a $\cB^-$-action on $\CN$ as an object of $\Pervgrib$, which is
compatible with the Hecke morphisms in the natural sense.

Let us denote this category by $\Catgeomi_{\cB^-}$. Let us consider also the
category $\Catgeomi_{\cG}$ of $\cG$-equivariant objects of $\Catgeomi$;
it is canonically equivalent to $\Pervgrib$.

\medskip

Let us now recall the following general construction. Let $\CN$ be a an object of
$\Catgeomi_{\cB^-}$. We claim that it gives rise to a functor
$$\QCoh(\cG/\cB^-)\to \Catgeomi.$$

Indeed, given $\CK\in \QCoh(\cG/\cB^-)$,
which we will view as a $\cB^-$-equivariant $\CO_{\cG}$-module, consider the tensor product
$\CK\underset{\CO_\cG}\otimes {}^\cG\CN$.
This is an object of $\Catgeomi$, endowed with an action of $\cB^-$, and we set
$$\CK\ast \CN:=\Bigl(\CK\underset{\CO_\cG}\otimes \CN\Bigr)^{\cB^-}.$$
The underlying object of $\Pervgrib$ is given by $\Bigl(\CK\otimes \CN)^{\cB^-}$.

Suppose now that $\CK$ is an object of $\QCoh(\cG/\cB^-)_{\cG}$. Then, by construction
$\CK\ast \CN$ belongs to $\Catgeomi_{\cG}\simeq \Pervgrib$.

\begin{lem}
The functor $\CN\mapsto \CO_{\cG/\cB^-}\ast \CN:\Catgeomi_{\cB^-}\to \Pervgrib$
is the right adjoint of the forgetful functor
$$\Pervgrib\simeq \Catgeomi_{\cG}\to \Catgeomi_{\cB^-}.$$
\end{lem}

Note that the functor $\CN\mapsto \CO_{\cG/\cB^-}\ast \CN:\Catgeomi_{\cB^-}\to \Pervgrib$
can be tautologically rewritten as $\CN\mapsto \CN^{\cB^-}$. The following is a
translation of the Borel-Bott-Weil theorem:

\begin{prop}  \label{BBW}
Assume that $\tw$ is right $W$-minimal. Then
$\CO_{\cG/\cB^-}\ast \CM^\tw\simeq \CW^{*,\tw}$, and for $i>0$
$$\on{R}^i\Bigl(\cB^-,\CM^\tw\Bigr)=0.$$
\end{prop}

\begin{proof}

The first assertion of the proposition is immediate from \lemref{B^- action on baby}(2).
To prove the second assertion, note that if $\CN$ is any Artinian $\cB^-$-equivariant
object of $\Catgeomi$, there exists $\cmu\in \cLambda^+$ large enough, so that
$\on{R}^i\Bigl(\cB^-,\bC^\cmu\otimes \CN\Bigr)=0$ for $i>0$.
This follows from the fact that the functor of derived $\cB^-$-invariants has a finite
cohomological dimension, and any Artinian
object of $\Catgeomi_{\cB^-}$, admits a left resolution, whose terms are of the form
$\CF\star \Rg\otimes U$, where $\CF\in \Pervgri$, and $U$
is a finite-dimensional representation of $\cB^-$.

Hence, for a given $w\in W$ and $\clambda'\in \cLambda^+$ large enough,
$$\on{RInv}(\cB^-,\CM^{w\cdot \clambda'})\simeq \CW^{*,w\cdot\clambda'}.$$

Note that the functor $\on{RInv}(\cB^-,\cdot):{\mathsf D}(\Catgeomi_{\cB^-})\to \Dgri$
commutes with the action of $\Dfli$ by convolutions. It suffices to remark
that if $\clambda'-\clambda\in \cLambda^+$, then
$\CM^{w\cdot \clambda}\simeq j_{!,w\cdot (\clambda-\clambda')\cdot w^{-1}}\star
\CM^{w\cdot \clambda'}$, and if $w\cdot \clambda$ is right $W$-minimal, then also
$\CW^{w\cdot \clambda}\simeq j_{!,w\cdot (\clambda-\clambda')\cdot w^{-1}}\star
\CW^{w\cdot \clambda'}$.

\end{proof}

\begin{cor}  \label{BBW negative}
Let $\tw\in W_{aff}$ be right {\em maximal} with
respect to $W$. Then for $i\neq \dim(\fn)$,
$H^i(\cB^-,\CM^\tw)=0$, and
$$H^{\dim(\fn)}(\cB^-,\CM^\tw)\simeq \CW^{!,\tw\cdot w_0}.$$
\end{cor}

\begin{proof}

Let $\tw=w\cdot \clambda$, and let $\cmu\in \cLambda^+$ be such that
$w\cdot (\clambda+\cmu)$ is {\it left} minimal with respect to $W$.
Then:
$$H^\bullet(\cB^-,\CM^\tw)\simeq j_{!,w\cdot (-\cmu)\cdot w^{-1}}\star
H^\bullet(\cB^-,\CM^{w\cdot (\clambda+\cmu)})\simeq
j_{!,w\cdot (-\cmu)\cdot w^{-1}}\star \CW^{*,w\cdot (\clambda+\cmu)}.$$

The latter is isomorphic to
$$j_{!, \tw}\star \delta_{1,\Gr}\simeq \CW^{!,\tw\cdot w_0}[-\dim(\fn)].$$

\end{proof}

Let now $...\to \CP_1\to \CP_0\to \BC_{\cB^-}$ be a left resolution of
the skyscraper on $\cG/\cB^-$, as in \secref{Baby as a quotient},
where each $\CP_i$ has the form $\CO(-\cmu_i)\otimes U^i$,
where $U^i$ are vector spaces.

Let $\tw$ be $w\cdot \clambda$. Tensoring by the line bundle $\CO(\clambda')$,
we can ensure that $\clambda+\clambda'-\cmu_i$ are such that
$w\cdot (\clambda+\clambda'-\cmu_i)$ is right $W$-minimal for
$i=0,...,\dim(\fn)$. We obtain that the complex
$$\Bigl(\CP_1\otimes \CO(\clambda')\Bigr)\ast \CM^\tw\to
\Bigl(\CP_0\otimes \CO(\clambda')\Bigr)\ast \CM^\tw\to \BC_{\cB^-}\ast\CM^\tw\to 0$$
is exact. However, $\BC_{\cB^-}\ast\CM^\tw\simeq \CM^\tw$, and
$$\CO(\cmu)\ast \CM^{w\cdot \clambda}\simeq \CW^{*,w\cdot (\cmu+\clambda)},$$
by \propref{BBW}, provided that $w\cdot (\cmu+\clambda)$ is right $W$-minimal.
Thus, we arrive to the same conclusion as in \propref{Baby Verma as a quotient}.
\sssec{Hereditary property}

In this subsection we will prove the following:

\begin{thm}  \label{geom hered}
$Ext^i_{\Catgeomni}(\BD(\tCM^{\tw}),\tCM^{\tw'})=0$ for $i>0$ and any $\tw,\tw'\in W_{aff}$, and
$Hom(\BD(\tCM^\tw),\tCM^{\tw'})$ is zero if $\tw\neq \tw'$, and $1$-dimensional otherwise.

\end{thm}

This theorem follows immediately from \thmref{ABG for small} due
to the corresponding property of baby co-Verma modules over the
small quantum group. Here we will discuss a geometric proof of
this fact, which the rest of this subsection is devoted to. In the
course of the proof we will introduce another important
object--the Wakimoto sheaf.

\medskip

Let $\Catgeomni_{\cB}$ be the category of $\cB$-equivariant objects in
$\Catgeomni$. By \propref{dual descr}, $\BD(\CM^{\tw})$ is naturally
an object of $\Catgeomni_{\cB}$.

For $\tw=w\cdot \cmu$ consider the following object $\on{Wak}^{\tw}$ of
$\Catgeomni_{\cB}$. It is defined as
$$\underset{\cnu\in \cLambda}{\underset{\longrightarrow}{lim}}\,
\CW^{*,w\cdot (\cnu+\cmu)}\star \Rg\{\cnu\},$$
where the maps in the inductive system, defined for $\cnu'-\cnu=\clambda\in \cLambda^+$,
are given by
\begin{align*}
&\CW^{*,w\cdot (\cnu+\cmu)}\star \Rg\{\cnu\}\to
\CW^{*,w\cdot (\cnu+\cmu)}\star \Rg\{\cnu+\clambda\}\otimes \uV^\clambda(\clambda)\to \\
&\to \CW^{*,w\cdot (\cnu+\cmu)}\star \CV^\clambda \star \Rg\{\cnu+\clambda\}\to
\CW^{*,w\cdot (\cnu+\cmu+\clambda)}\star \Rg\{\cnu+\clambda\}.
\end{align*}

\medskip

Note that the forgetful functor $\Catgeomni_{\cB}\to \Catgeomni$
admits a natural right adjoint given by
$\CN\mapsto {}\CO_{\cB}\underset{\CO_\cG}\otimes {}^\cG\CN$.
Similarly, the functor $\Catgeomni_{\cB}\to \Catgeomdni$ admits
a right adjoint $\CN\mapsto {}\CO_{\cN}\underset{\CO_\cG}\otimes {}^\cG\CN$.

\begin{lem}
$$\on{Wak}^{\tw}\simeq\CO_{\cN}\underset{\CO_\cG}\otimes {}^\cG(\tCM^{\tw}).$$
\end{lem}

Hence, we obtain
$$Ext^i_{\Catgeomdni}(\BD(\tCM^{\tw}),\tCM^{\tw'})\simeq
Ext^i_{\Catgeomni_\cB}(\BD(\CM^{\tw}),\on{Wak}^{\tw'}).$$
By the Artinian property and taking into account \corref{duality on all baby},
to prove \thmref{geom hered}, it is sufficient to show that
\begin{equation} \label{from Weyl to Wak}
Ext^i_{\Catgeomni_\cB}({}^{w_0}\CM^{w\cdot \cmu},\CW^{w'\cdot (\clambda+\cmu')}\star \Rg\{\clambda\})=0,
\end{equation}
unless $i=0$, $w'=w\cdot w_0$ and $\cmu'=w_0(\cmu)+2\crho$, whenever
$\clambda$ is deep in $\cLambda^+$.

\begin{lem}
For $\CN\in \Catgeomni_\cB$ and $\CF\in \Catgeomni$,
$$RHom_{\Catgeomni_\cB}(\CN,\CF\star \Rg)\simeq
RHom_{\Catgeomni}\Bigl(\on{RCoinv}(\cB,\CN),\CF\Bigr).$$
\end{lem}

\begin{proof}

It is sufficient to prove the assertion in the case when $\CN\simeq \CF^1\star \tR\otimes U$,
where $U$ is a representation of $\cB$. In this case, it amounts to the following
adjunction, which is a corollary of the Serre duality on $\cG/\cB$:
$$RHom_\cB(U,\Res^\cG_\cB(V))\simeq
RHom_\cG\Bigl(\on{RCoinv}(\cB,\CO_\cG\otimes U),V\Bigr).$$

\end{proof}

\begin{lem}
For $\clambda$ deep in the dominant chamber,
$$\on{RCoinv}\Bigl(\cB,({}^{w_0}\CM^{w\cdot \cmu})\{-\clambda\}\Bigr)\simeq
\CW^{!,w\cdot w_0(\clambda-w_0(\cmu)+2\crho)}.$$
\end{lem}

\begin{proof}

First,
$$\on{RCoinv}\Bigl(\cB,({}^{w_0}\CM^{w\cdot \cmu})\{-\clambda\}\Bigr)\simeq
\on{RCoinv}\Bigl(\cB^-,\CM^{w\cdot \cmu}\{-w_0(\clambda)\}\Bigr).$$

Note that for $\CN\in \Catgeomni_{\cB^-}$,
$$\on{RCoinv}(\cB^-,\CN)\simeq \on{RInv}(\cB^-,\CN)\{2\crho\}[\dim(\fn)].$$
Hence, the expression in the lemma is isomorphic to
$$\on{RInv}(\cB^-,\CM^{w\cdot \cmu}\{-w_0(\clambda)+2\crho\})[\dim(\fn)] \simeq
\CW^{!,w\cdot w_0\cdot (w_0(\cmu)+\clambda+2\crho)},$$
by \corref{BBW negative}.

\end{proof}

\medskip

Thus, we obtain that the expression in \eqref{from Weyl to Wak} is isomorphic to
$$Ext^i_{\Catgeomni}(\CW^{!,w\cdot w_0(\clambda+w_0(\cmu)+2\crho)},
\CW^{*,w'\cdot (\cmu'+\clambda)}),$$
for which the vanishing assertion is manifest.

\sssec{An application: 2-sided BGG resolution}

We will use the geometric interpretation of baby co-Verma modules
to prove the following result:

\begin{thm}
There exists an exact complex $\CB_{\semiinf}$ of objects of $\Catgeomdi$,
whose $n$-th term is $$\underset{\tw\in W_{aff},l^\semiinf(\tw)=n}\bigoplus\, \tCM^\tw.$$
\end{thm}

We remind that for $\tw=w\cdot \clambda\in W_{aff}$, its semi-infinite length
$l^\semiinf(\tw)$ is defined as $l(w\cdot (\clambda+\cmu))-l(\cmu)$ for some
(or all) large $\cmu\in \cLambda^+$.

Of course, using the equivalence between $\Catgeomdi$ and
$\tu_\ell\modo$, we obtain the corresponding exact complex
consisting of baby co-Verma modules over $\tu_\ell$. The rest of
this subsection os devoted to the proof of this theorem.

\medskip

Let $\CB_{\Gr}$ be the Cousin complex on $\Gr$. I.e., this is an exact complex of
perverse sheaves on $\Gr$, living in positive degrees, whose $n$-th term is given by
$$\underset{\tw\in W_{aff}/W,l(\tw)=n}\bigoplus\, \CW^{*,\tw}.$$

For $\cmu\in \cLambda^+$ consider the complex $\CB_{\Gr}\star \IC_{-w_0(\cmu),\Gr}[l(\cmu)]$.
This complex is acyclic, since convolution with $\IC_{-w_0(\cmu),\Gr}$ is an exact functor.

We claim that for $\cmu'=\cmu+\cnu$ with $\cnu\in \cLambda^+$ we have a map of complexes
$$\CB_{\Gr}\star \IC_{-w_0(\cmu),\Gr}[l(\cmu)]\to
\CB_{\Gr}\star \IC_{-w_0(\cmu'),\Gr}[l(\cmu')].$$

For $w\cdot \clambda\in W_{aff}/W$, the map
$$\CW^{*,w\cdot (\clambda+\cmu)}\star \IC_{-w_0(\cmu),\Gr}\to
\CW^{*,w\cdot (\clambda+\cmu')}\star \IC_{-w_0(\cmu'),\Gr}$$ has been constructed
in the definition of the inductive system that defines $\CM^{w\cdot \clambda}$.

\medskip

To check that this map respects the differential, we must show the following.
Let $\tw=w\cdot \clambda$ and $\tw'=w'\cdot \clambda'$ be such that
$l(\tw')=l(\tw)+1$, and the orbit $I\cdot \tw$ is in the closure
of $I\cdot \tw'$. Then we claim that for $\cmu\in \cLambda^+$, the orbit
$I\cdot (w\cdot (\clambda+\cmu))$ is in the closure of
$I\cdot (w\cdot (\clambda'+\cmu))$, and the square
$$
\CD
\CW^{*,w\cdot \clambda}\star \IC_{\cmu,\Gr} @>>>
\CW^{*,w'\cdot \clambda'}\star \IC_{\cmu,\Gr} \\
@VVV   @VVV   \\
\CW^{*,w\cdot (\clambda+\cmu)} @>>> \CW^{*,w'\cdot (\clambda'+\cmu)}
\endCD
$$
commutes, where the horizontal arrows are the canonical maps, corresponding
to adjoining orbits.

\begin{lem}
If $I\cdot \tw\subset \ol{I\cdot \tw'}$, then as elements of $W_{aff}$,
$\tw\leq \tw'$.
\end{lem}

\begin{proof}

We need to show that $I\cdot \tw_{\Fl}\subset \ol{I\cdot \tw'_{\Fl}}$,
where the subscript $_{\Fl}$ means that we are dealing with an
orbit in $\Fl$ (vs. $\Gr$).

Since the projection $\Fl\to \Gr$ is proper, there exists some
$\tw_1\in W_{aff}$, such that $\tw_1<\tw'$, and $\tw_1=\tw \on{mod}\, W$.
We have:
$$l(\tw)\leq l(\tw_1)\leq l(\tw')-1.$$
Since $l(\tw)=l(\tw')-1$, we obtain that $\tw_1=\tw$.

\end{proof}

Note that by the lemma, we obtain that $w\cdot (\clambda+\cmu)\leq w'\cdot (\clambda'+\cmu)$,
and hence we do have a containment $$I\cdot (w\cdot (\clambda+\cmu))\subset
\ol{I\cdot (w'\cdot (\clambda'+\cmu))}.$$

Also, by the lemma, the map $\CW^{*,w\cdot \clambda}\to
\CW^{*,w'\cdot \clambda'}$ is obtained from the map
\begin{equation} \label{aux 3}
j_{*,w\cdot \clambda}\to j_{*,w'\cdot \clambda'}
\end{equation}
by convolving with $\delta_{1,\Gr}$. Note that the map
$j_{*,w\cdot (\clambda+\cmu)}\to j_{*,w'\cdot (\clambda'+\cmu)}$ is obtained from
\eqref{aux 3} by convolving on the right with $j_{*,\cmu}$.

To prove the commutativity of the above diagram, it suffices to notice that the left vertical
arrow is equal to the composition
\begin{align*}
&\CW^{*,w\cdot \clambda}\star \IC_{\cmu,\Gr}\simeq
j_{*,w\cdot \clambda}\star \IC_{\cmu,\Gr}\to
j_{*,w\cdot \clambda}\star \CW^{*,\cmu}\simeq \\
&j_{*,w\cdot \clambda}\star j_{*,\cmu}\star \delta_{1,\Gr}\simeq
j_{*,w\cdot (\clambda+\cmu)} \star \delta_{1,\Gr}\simeq
\CW^{*,w\cdot (\clambda+\cmu)},
\end{align*}
and similarly for the right vertical arrow.

\section{Sheaves on semi-infinite flags}

\ssec{Drinfeld's spaces and factorization}

\sssec{}

Let $X$ be a global curve. Let $\Bun_G$ denote the moduli stack of principal
$G$-bundles on $X$. Let us recall the definition of the Drinfeld space $\BunNb$.
\footnote{The esposition in the section substantially relies on the results of
\cite{BG}, \cite{BFGM} and \cite{FGV}, and certain familiarity with these papers
will be assumed.}

First we define a bigger space $\BunNb'$ that classifies the data of a $G$-bundle $\CP_G$
on $X$, and its generalized reduction to $N$, i.e.,
a collection of non-zero maps defined for each $\lambda\in \Lambda^+$
$$\kappa^\lambda:V^\lambda_{\CP_G}\to \CO_X,$$
where $V^\lambda$ is the corresponding Weyl module over $G$, and
$V^\lambda_{\CP_G}$ is the associated vector bundle. The collection $\kappa^\lambda$
is required to satisfy the Pl\"ucker relations, cf. \cite{FM,BG}. We will denote by $\fp$ the
tautological projection $\BunNb'\to \Bun_G$.

We have a natural action of $T$ on $\BunNb$: an element $t\in T$
multiplies each $\kappa^\lambda$ by $\lambda(t)$. It is easy to see that the
map $\BunNb'/T\to \Bun_G$ is proper.

\medskip

If $[G,G]$ is simply-connected, then $\BunNb'$ is the sought-for Drinfeld space
$\BunNb$. Otherwise we proceed as follows.

\begin{lem}  \label{isogeny}
Let $G_1\to G_2$ be an isogeny, i.e., a homomorphism, whose kernel is
contained in $Z(G_1)$ and whose image contains $[G_2,G_2]$.
Then the natural map $\BunNb'(G_1)\to \BunNb'(G_2)$ is a closed
embedding.
\end{lem}

\begin{proof}

First, it is easy to see that if we have a short exact sequence
$$1\to G'\to G\to T'\to 1,$$
where $T'$ is a torus then $\BunNb'(G')\to \BunNb'(G)$ is an
isomorphism. This reduces the assertion of the lemma to the case when
$G_1\to G_2$ is surjective with finite kernel. Let $k$ be the index if
$\cLambda_1$ in $\cLambda_2$.

Since each of $\BunNb'(G_i)/T_i$, $i=1,2$ is proper over $\Bun_{G_i}$,
the map $$\BunNb'(G_1)/T_1\to \BunNb'(G_2)/T_2$$ is proper. Hence,
it remains to see that the map $\BunNb'(G_1)\to \BunNb'(G_2)$
is injective on the level of $S$-points for any base $S$.

\medskip

Let $(\CP_{G_2}, \{\kappa^\lambda_2\})$ be an
$S$-point of $\BunNb'(G_2)$, and let $(\CP_{G_1}, \{\kappa^\lambda_1\})$
be its lift to a point of $\BunNb'(G_1)$. Then the image of
$\kappa^\lambda_1$ in $V^\lambda_{\CP_{G_1}}$ is fixed by the
condition that
$$(\kappa^\lambda_1)^{\otimes k}=\kappa^{k\cdot \lambda}_2:
V^{k\cdot \lambda}_{\CP_{G_1}}\simeq V^{k\cdot \lambda}_{\CP_{G_2}}\to \CO_X.$$

Hence, when $\CP_{G_1}$ is fixed, any two choices of systems $\{\kappa^\lambda_1\}$
differ by an element of $T_{1,2}:=\on{ker}(T_1\to T_2)\simeq \on{ker}(G_1\to G_2)$.
However, two such lifts are isomorphic as points of $\BunNb'(G_1)$, via the
automorphism of $\CP_{G_1}$ given by the same element of $T_{1,2}$.

\medskip

Finally, if $\CP'_{G_1}$ is another principal $G_1$-bundle that reduces to
$G_2$, there exists a principal $T_{1,2}$-bundle $\CP_{T_{1,2}}$, such that
$\CP'_{G_1}\simeq \CP_{G_1}\overset{T_{1,2}}\times \CP_{T_{1,2}}$. Then
for every $\lambda$ as above,
$$V^\lambda_{\CP'_{G_1}}\simeq V^\lambda_{\CP_{G_1}}\otimes \CP^\clambda_{T_{1,2}},$$
where $\CP^\clambda_{T_{1,2}}$ is the line bundle associated with $\CP_{T_{1,2}}$
and the character $\lambda$.

However, the data of $\kappa^\lambda_1$ for $V^\lambda_{\CP'_{G_1}}$
identifies the line sub-sheaf
$$(\CP^\clambda_{T_{1,2}})^{-1}\subset (V^\lambda)^*_{\CP_{G_1}}\otimes
(\CP^\clambda_{T_{1,2}})^{-1}$$
with $\CO_X$, thereby giving a trivialization of $\CP_{T_{1,2}}$.

\end{proof}

\bigskip

For an arbitrary group $G$ we can find a group $G'$ with a surjective isogeny
$G'\to G$, such that (a) $\on{ker}(G'\to G)$ is connected, and (b) $[G',G']$ is
simply connected.

We define $\BunNb$ as the image of $\BunNb(G')=\BunNb'(G')$ in $\BunNb'$
under $$\BunNb'(G')\to \BunNb'(G).$$ By the above lemma, this is a closed substack
of $\BunNb'$, and it is easy to see that it does not depend on the choice
of $G'$.

\sssec{Variants}

We fix a point $x\in X$. For a coweight $\cnu$, let $_{\leq \cnu}\BunNb'$
denote a version of $\BunNb'$, where we allow each $\kappa^\lambda$
have a pole at $x$ of order $\leq\langle \lambda,\cnu\rangle$.

For $G'$ as above, due to the fact that the kernel of $G'\to G$ is connected,
we can find a preimage $\cnu'$ of $\cnu$ in the coweight
lattice $\cLambda'$ of $G'$, and we define $_{\leq \cnu}\BunNb\subset {}_{\leq \cnu}\BunNb'$
as the image of $_{\leq \cnu'}\BunNb'(G')$ under
\begin{equation} \label{isogeny stratum}
_{\leq \cnu'}\BunNb'(G')\to {}_{\leq \cnu}\BunNb'.
\end{equation}

As in \lemref{isogeny} one shows that the map in \eqref{isogeny stratum}
is a closed embedding. Moreover, its image is easily seen to be
independent of the choice of $\cnu'$ for a fixed $G'$, and of $G'$ itself.

\medskip

If $\cnu_1-\cnu_2\in \cLambda^{pos}$ we have a natural closed embedding
$_{\leq \cnu_2}\BunNb\hookrightarrow {}_{\leq \cnu_1}\BunNb$. We define
${}_\infty\BunNb$ as
$$\underset{\cnu\in \cLambda}{\underset{\longrightarrow}{lim}}\,
\Bigl({}_{\leq \cnu}\BunNb\Bigr)$$
with respect to the natural ordering on $\cLambda$ and the above closed embeddings.

By definition, $_\infty\BunNb$ splits into connected components, numbered by the
quotient of $\cLambda$ by the coroot lattice.

\medskip

Let $_{\cnu}\BunNb'$ is an open substack of $_{\leq \cnu}\BunNb'$ corresponding
to the condition that each $\kappa^\lambda$ has a pole of order exactly
$\leq\langle \lambda,\cnu\rangle$ at $x$. Set
$$_{\cnu}\BunNb:={}_{\cnu}\BunNb'\cap {}_{\leq \cnu}\BunNb.$$
One easily shows that $_{\cnu}\BunNb$ equals the image of
$_{\cnu'}\BunNb'(G')$ under the map of \eqref{isogeny stratum}.

Let us note that over each $_\cnu\BunNb$ there exists a canonical
$N^-[[t]]$-torsor, which we will denote by $_\cnu\CN$. We will denote by
$_\cnu^k\CN$ the induced $N^-([t]/t^k)$-torsor.

We will denote by $\oi_{\leq \cnu}$ (resp., $\oi_\cnu$) the closed (resp., locally closed)
embedding of $_{\leq \cnu}\BunNb$ (resp., $_\cnu\BunNb$) into $_\infty\BunNb$.
We have:
$$_\cnu\BunNb={}_{\leq \cnu}\BunNb-\underset{\cnu_1<\cnu}\cup\, {}_{\leq \cnu_1}\BunNb.$$

\medskip

We let $_\cnu\BunN$ denote the open sub-stack of $_\cnu\BunNb$, where we
demand that the maps $\kappa^\lambda$ have no zeroes away from $x$. This
substack is isomorphic to $\Bun_{B^-}\underset{\Bun_T}\times {pt}$, where the map
$pt\to \Bun_T$ corresponds to the point $\CP^0_T(\cnu\cdot x)$. We will denote
by $i_\nu$ the locally closed embedding of $_\cnu\BunN$ into $_\infty\BunNb$;
by \cite{FGV}, Sect. 3.3, the morphism $i_\nu$ is affine.

\medskip

Let $\ol{x}':=x'_1,...,x'_m$ be a collection of points on $X$, distinct from $x$. Let
$_\infty\BunNb^{\text{n.z.}\ol{x}'}$ be the open sub-stack of $_\infty\BunNb$ defined by the
condition that the maps $\kappa^\lambda$ have no zeroes at $x'_1,...,x'_m$. As
in \cite{FGV}, Sect. 3.2, one shows that over $_\infty\BunNb^{\text{n.z.}\ol{x}'}$ there exists a natural
torsor with respect to the group-scheme $\underset{j=1,...,m}\Pi\, N^-[[t'_j]]$,
denoted $\CN^{\ol{x}'}$, where $t'_j$ is a local coordinate at $x'_j$. Moreover,
$\CN^{\ol{x}'}$ carries an action of the group-indscheme $\underset{j=1,...,m}\Pi\, N^-((t'_j))$.

\medskip

For  an integer $k$ let $^k\Bun_G$ denote the principal $G^k$-bundle
over $\Bun_G$ corresponding to choosing a structure of level $k$ at $x$ in a $G$-bundle.
We will denote by $^k_\infty\BunNb$ the Cartesian product
$_\infty\BunNb\underset{\Bun_G}\times {}^k\Bun_G$.

We will denote by $_{\leq \cnu}^k\BunNb$, $_{\cnu}^k\BunNb$,
$_{\cnu}^k\BunN$ the corresponding stacks obtained by base change.
By a slight abuse of notation, we will use the symbols
$\oi_{\leq \cnu}$, $\oi_\cnu$ and $i_\cnu$ for the embeddings
of these stacks into $^k_\infty\BunNb$.
Similary, we introduce the stacks $_\infty^k\BunNb^{\text{n.z.}\ol{x}'}$,
$^k\CN^{\ol{x}'}$ as Cartesian products.

\medskip

Note that there is a natural isomorphism
\begin{equation} \label{descr open stratum}
_\cnu^k\BunNb\simeq G([t]/t^k)\overset{N^-([t]/t^k)}\times {}_\cnu^k\CN.
\end{equation}

In particular, we obtain a natural map
$$\ol{\on{ev}}_\cnu:{}_\cnu^k\BunNb\to \bigl(G/N^-\bigr)([t]/t^k)\to
\bigl(G/B^-\bigr)([t]/t^k).$$
The restriction of this map to $_\cnu^k\Bun_N$, denoted $\on{ev}_\cnu$, is smooth.

\sssec{}

For $\cmu\in \cLambda$ let $\Bun_B^\cmu$ be the corresponding connected component of
$\Bun_B$. We recall that $\Bun_B^\cmu$ can be interpreted as the stack classifying the
data of a principal  $G$-bundle $\CP_G$ on $X$, a $T$-bundle $\CP_T$, such that each
associated line bundle
$\CP_T^\lambda$ has degree $-\langle \lambda,\cmu\rangle$, and a collection of
bundle maps
$$\kappa^{\lambda,-}:\CP_T^\lambda\to V^\lambda_{\CP_G},$$
defined for $\lambda\in \Lambda^+$, which satisfy the Pl\"ucker relations.
(Here $\CP_T^\lambda$ denoted the line bundle associated with $\CP_T$
and the character $\lambda:T\to \BG_m$.)

Note that if $\cmu$ is such that $\langle \alpha,\cmu\rangle>(2g-2)$ for all positive roots
$\alpha$, then the map $\fp^{\cmu,-}:\BunBm^\cmu\to \Bun_G$ is smooth.

\medskip

Consider the Cartesian product $_\infty\BunNb\underset{\Bun_G}\times \Bun_B^\cmu$.
We will denote by $_\infty\CZ^\cmu$ the corresponding Zastava space, i.e. the open
substack of the above Cartesian product, defined by the condition that the reductions
to $N^-$ and $B$ are transversal at the generic point of the curve. This means that the
composed maps
$$\CP_T^\lambda\overset{\kappa^{\lambda,-}}\longrightarrow
V^\lambda_{\CP_G}\overset{\kappa^\lambda}\longrightarrow
\CO_X$$ are non-zero for all $\lambda\in \Lambda^+$.

We will denote by $_\infty^k\CZ^\cmu$ the stack obtained by adding a structure of
level $k$ to the $G$-bundle $\CP_G$ at $x$. All of the above stacks are acted on
by the group $T$.

\medskip

Let us denote by $_{\leq \cnu}^k\CZ^\cmu$ (resp., $_\cnu^k\CZ^\cmu$, $^k\CZ^\cmu$)
the preimage in $_{\infty}^k\CZ^\cmu$ of the substack $_{\leq \cnu}^k\BunNb$
(resp., $^k_\cnu\BunNb$, $^k\BunNb={}^k_{\leq 0}\BunNb$) of $_{\infty}^k\BunNb$.
Note that $_{\leq \cnu}^k\CZ^\cmu$ is empty unless
$\cnu+\cmu\in \cLambda^{pos}$.
By $_\cnu^k\oCZ^\cmu$ we will denote the open substack of $_\cnu^k\CZ^\cmu$
equal to the preimage of $_\cnu^k\BunN$.

\medskip

For $\cmu\in \cLambda^{pos}$, let $X^\cmu$ be the corresponding partially symmetrized power
of the curve. By definition, $X^\cmu$ classifies the data of a principal $T$-bundle $\CP_T$
and its generic trivialization, such that for $\lambda\in \Lambda^+$ the resulting maps
$\CP_T^\lambda\to \CO_X$ are all regular and the divisor of zeroes has degree
$\langle \lambda,\cmu\rangle$.

For $\cnu\in \cLambda$, let $_{\leq \cnu}X^\cmu$
be a version of $X^\cmu$, where the maps $\CP_T^\lambda\to \CO_X$ are allowed to
have poles at $x$ of order $\leq\langle \lambda,\cnu\rangle$ for $\lambda\in \Lambda^+$.
This space is empty unless $\cmu+\cnu\in \cLambda^{pos}$. If
$\cnu_1-\cnu_2\in \Lambda^{pos}$ we have a natural closed embedding
$$_{\leq \cnu_2}X^\cmu\hookrightarrow{}_{\leq \cnu_2}X^\cmu.$$
We define $_{\infty}X^\cmu$ as the ind-scheme
$$_{\infty}X^\cmu=\underset{\cnu\in \cLambda}{\underset{\longrightarrow}{lim}}\,
\Bigl({}_{\leq \cnu}X^\cmu\Bigr)$$
with respect to the usual ordering on $\cLambda$ and the above closed
embeddings. This space also splits into connected components numbered by
the quotient of $\cLambda$ by the coroot lattice.

\medskip

By construction, we have a natural map
$$_\infty\fs^\cmu:{}_\infty\CZ^\cmu\to {}_{\infty}X^\cmu.$$
We will denote the restriction of $_\infty\fs^\cmu$ to
$_{\leq \cnu}\CZ^\cmu$ (resp., $\CZ^\cmu={}_{\leq 0}\CZ^\cmu$) by
$_{\leq \cnu}\fs^\cmu$ (resp., $\fs^\cmu$). Note that
$_{\leq \cnu}\fs^\cmu$ maps to $_{\leq \cnu}X^\cmu$.

We will denote by $^k\fs^\cmu$ the
composition of $\fs^\cmu$ and the forgetful map
$_\infty^k\CZ^\cmu\to {}_\infty\CZ^\cmu$, and similarly for
$^k_{\leq \cnu}\fs^\cmu$, $^k\fs^\cmu$.

\medskip

Let $\oX$ denote the open curve $X-x$, and $\oX^\cmu$ be the corresponding open
subset of $X^\cmu$. For $\cmu_1,\cmu_2$ we will denote by
$\Bigl(\oX^{\cmu_1}\times {}_{\infty}X^{\cmu_2}\Bigr)_{disj}$ the open subset
in the product $\oX^{\cmu_1}\times {}_{\infty}X^{\cmu_2}$, corresponding to the condition
that the two divisors have disjoint support.

As in \cite{BFGM}, we have:

\begin{lem} \label{factorization}
For $\cmu_1+\cmu_2=\cmu$ there exist natural isomorphisms
$$_\infty^k\CZ^\cmu \underset{_{\infty}X^\cmu}\times
\Bigl(\oX^{\cmu_1}\times {}_{\infty}X^{\cmu_2}\Bigr)_{disj}\simeq
\Bigl(\CZ^{\cmu_1} \times {}_\infty^k\CZ^{\cmu_2}\Bigr)
\underset{X^{\cmu_1}\times {}_{\infty}X^{\cmu_2}}\times
\Bigl(\oX^{\cmu_1}\times {}_{\infty}X^{\cmu_2}\Bigr)_{disj}.$$
\end{lem}

\medskip

Let $x'$ be any point of the curve, and for $\cmu\in \cLambda^{pos}$, let
$\cmu\cdot x'$ be the corresponding element of $\oX^\cmu$. Then, by \cite{BFGM}, we have:
\begin{equation} \label{non-central fibers}
(\fs^\cmu)^{-1}(\cmu\cdot x')\simeq \left(N((t'))\cdot \cmu\right)\cap
\overline{\left(N^-((t'))\cdot 1_{\Gr}\right)},
\end{equation}
where $t'$ is a local coordinate at $x'$.

\medskip

In the same way we obtain that for an arbitrary element
$\cmu\in \cLambda$ and the point $\cmu\cdot x\in {}_\infty{}X^\cmu$
\begin{equation} \label{central fiber}
^k_{\infty}\fF^\cmu:=({}^k_{\infty}\fs^\cmu)^{-1}(\cmu\cdot x)\simeq
\Bigl(N((t))\cdot \cmu\Bigr)
\underset{\Gr}\times G((t))/G^k.
\end{equation}

\sssec{}

For an integer $m$, let $\on{Jets}^+(T)^m$ be the group-scheme over $X^{(m)}$,
whose fiber over a divisor $\Sigma\, m_j\cdot x_j$, where the points $x_j$ are
distinct, is $\underset{j}\Pi\, T[[t'_j]]$. More precisely, for a test-scheme $S$ and
an $S$-point $\varphi$ of $X^{(m)}$ its lift to an $S$-point of $\on{Jets}^+(T)^m$ is an
$X^{(m)}$-map
$$\widehat{\Gamma}_\varphi\to T,$$
where $\widehat{\Gamma}_\varphi\subset S\times X$ is the formal neighbourhood of the
preimage $\Gamma_\varphi$ of the incidence divisor in $X^{(m)}\times X$ under
$\varphi\times \on{id}$.

If $\cmu,\cnu\in \cLambda$ are two elements with
$\cmu+\cnu\in \cLambda^{pos}$, we have a natural map
$_{\leq\cnu}X^\cmu\to X^{(m)}$, where $m=l(\cmu+\cnu)$, and let
$_{\leq \cnu}\on{Jets}^+(T)^\cmu$ be the resulting group-scheme on $_{\leq\cnu}X^\cmu$.

\begin{propconstr}
The group-scheme $_{\leq \cnu}\on{Jets}^+(T)^\cmu$ acts naturally
on $_{\leq \cnu}\CZ^\cmu$.
\end{propconstr}

\begin{proof}

To simplify the notation, we will assume that $_\cnu=0$, and we will work with
the "usual" Zastava space $\CZ^\cmu$.

According to \cite{BFGM}, Sect. 2, given an $S$-point of  $\CZ^\cmu$, the resulting
$G$-bundle $\CP_G$ on $S\times X$ acquires a trivialization on $S\times X-\Gamma_\varphi$,
where $\varphi$ is the composition of the initial map to $\CZ^\cmu$ and
$$\CZ^\cmu\to X^\cmu\to X^{(m)}.$$

As usual in this situation, given a map $g_S:\widehat{\Gamma}_\varphi\to G$, we can
produce a new  $G$-bundle $\CP'_G$, by declaring it to be the same as $\CP_G$
on $S\times X-\Gamma_\varphi$ and $\widehat{\Gamma}_\varphi$ and changing
the gluing data on the formal punctured neighbourhood of $\Gamma_\varphi$ by
means of $g_S$.

If $g_S$ was a map $\widehat{\Gamma}_\varphi\to T$, then the data of $\kappa^\lambda$
and $\kappa^{\lambda,-}$ for $\CP_G$ give rise to well-defined data of
$(\kappa^\lambda)'$ and $(\kappa^{\lambda,-})'$ for $\CP'_G$. Thus, we obtain
a new point of $\CZ^\cmu$.

\end{proof}

Note that $_{\leq \cnu}\on{Jets}^+(T)^\cmu$ contains as a direct factor
the constant group-subscheme with fiber $T$.
Its action on $_{\leq \cnu}\CZ^\cmu$ coincides with the "global"
one, mentioned above.

\medskip

Let us consider now $_{\leq \cnu}^k\CZ^\cmu$. One can show that the above action of
$_{\leq \cnu}\on{Jets}^+(T)^\cmu$ on $_{\leq \cnu}\CZ^\cmu$ {\it does not} lift to
an action of $_{\leq \cnu}^k\CZ^\cmu$. However, we do have an action fiber-wise
over each point of $_{\leq \cnu}X^\cmu$. For example, the action of $T[[t]]$ on
$^k_{\cmu}\fF^\cmu$ is given in terms of isomorphism
$$^k_{\cmu}\fF^\cmu\simeq \Bigl(N((t))\cdot \cmu\cap
\overline{\left(N^-((t'))\cdot (-\cnu)\right)}\Bigr)
\underset{\Gr}\times G((t))/G^k,$$
by the natural action of $T((t))$ on $G((t))$ by left multiplication.

\medskip

We will use the following construction. Let us choose an identification
$T\simeq \BG_m^r$, and a point $y\in X-x$. For a string of positive
integers $\ol{m}=m_1,...,m_r$, consider the affine space consisting of
$r$-tuples of functions $(X-y)\to {\mathbb A}^1$, whose values at $x$ is $1$, and the
pole of the $i$-th function at $y$ is of order $\leq m_i$. We will denote this
space by $\on{Maps}(X,T)^{\ol{m}}$.

The Abel-Jacobi map gives rise to a morphism $\on{Maps}(X,T)^{\ol{m}}\to \oX^{\ol{m}}:=
\underset{i}\Pi\, \oX^{(m_i)}$, and we have a natural morphism
$$\left(\on{Maps}(X,T)^{\ol{m}}\times X\right)\underset{\oX^{\ol{m}}\times X}\times
\left(\oX^{\ol{m}}\times X\right)_{disj}\to T,$$
where $\left(\oX^{\ol{m}}\times X\right)_{disj}\subset \oX^{\ol{m}}\times X$ has the same meaning
as before--the complement to the incidence divisor. (This morphism explains the
notation $\on{Maps}(X,T)^{\ol{m}}$ for the above scheme.)

\begin{propconstr}
We have a natural map
$$\on{act}_T:\left(\on{Maps}(X,T)^{\ol{m}}\times {}
_{\leq \cnu}^k\CZ^\cmu\right)\underset{\oX^{\ol{m}}\times {}_{\leq \cnu}X^\cmu}\times
\left(\oX^{\ol{m}}\times {}_{\leq \cnu}X^\cmu\right)_{disj}\to {}_{\leq \cnu}^k\CZ^\cmu.$$
\end{propconstr}

\begin{proof}

We retain the notation from the proof of the previous proposition-construction.
The difference now is that the map $g_S$ is defined on a Zarisky-open of
$S\times X$ that contains $\Gamma_\varphi$ and $S\times x$. In particular,
the restrictions of $\CP_G$ and $\CP'_G$ to the formal neighborhood of
$x$ are identified. Hence, $\CP'_G$ is also equipped with a structure of
level $k$ at $x$.

\end{proof}

\ssec{A category of perverse sheaves}

\sssec{}   \label{category of sheaves}

For an integer $k$ we define the category $\Pervsmk$ to be the full subcategory
of the category of $T$-equivariant perverse sheaves on $^k_{\infty}\BunNb$, consisting of
objects satisfying the following three properties:

\medskip

\noindent{(1)}
For a finite collection $\ol{x}=x'_1,...,x'_m$ of points on $X$ distinct from $x$,
the pull-back of $\CF$ to $^k\CN^{\ol{x}'}$ is equivariant with respect to the
group-indscheme $\underset{j=1,...,m}\Pi\, N^-((t'_j))$.

\noindent{(2)} The factorization property:

We say that a perverse sheaf $\CF$ on $^k_{\infty}\BunNb$ is factorizable if
for any $\cmu_1,\cmu_2$, satisfying $\langle \alpha,\cmu_i\rangle>(2g-2)$
and $\cmu_2-\cmu_1\in \cLambda^{pos}$, the retsriction of the
pull-back  $\fp^{-,\cmu_2}{}^*(\CF)$ onto the left-hand side of
$$_\infty^k\CZ^{\cmu_2} \underset{_{\infty}X^{\cmu_2}}\times
\Bigl(\oX^{\cmu_2-\cmu_1}\times {}_{\infty}X^{\cmu_1}\Bigr)_{disj}\simeq
\Bigl(\CZ^{\cmu_2-\cmu_1} \times {}_\infty^k\CZ^{\cmu_1}\Bigr)
\underset{X^{\cmu_2-\cmu_1}\times {}_{\infty}X^{\cmu_1}}\times
\Bigl(\oX^{\cmu_2-\cmu_1}\times {}_{\infty}X^{\cmu_1}\Bigr)_{disj}$$
is isomorphic (up to a cohomological shift by the corresponding relative
dimensions) to the restriction onto the right-hand side of the
external product
$$\IC_{\CZ^{\cmu_2-\cmu_1}}\boxtimes\, \fp^{-,\cmu_1}{}^*(\CF).$$
(Note that both complexes in question are perverse sheaves, since the
maps $\fp^{-,\cmu_i}$, $i=1,2$ are smooth by assumption.)

\medskip

\noindent{(3)} If $\CF$ is supported on $^k_{\leq \cnu}\BunNb$, then
for $\cmu\in \Lambda$, satisfying $\langle \alpha,\cmu\rangle>(2g-2)$,
the pull-back of $\CF$ on $_{\leq \cnu}^k\CZ^\cmu$ is
$\on{Maps}(X,T)^{\ol{m}}$-equivariant for any $\ol{m}$. The latter
means that there exists an isomorphism between two pull-backs of
$(\fp^{\cmu,-})^*(\CF)|_{_{\leq \cnu}^k\CZ^\cmu}$
to
$$\left(\on{Maps}(X,T)^{\ol{m}}\times {}
_{\leq \cnu}^k\CZ^\cmu\right)\underset{\oX^{\ol{m}}\times {}_{\leq \cnu}X^\cmu}\times
\left(\oX^{\ol{m}}\times {}_{\leq \cnu}X^\cmu\right)_{disj},$$
which induces the identity map on the further restriction of both sides to the unit point
of $\on{Maps}(X,T)^{\ol{m}}$.

\medskip

\noindent{\it Remark.}
As we shall see, imposing property (1) is in fact superfluous, i.e., it follows
formally from the factorization property (2). In addition, some portion of
property (2) follows from (1).

In addition, if $k=1$ (which the main case of interest for this paper),
property (3) follows automatically.
\footnote{We remark also that property (3) has to do with the fact that our
category $\Pervsmk$ models perverse sheaves on $G((t))/N^-((t))\cdot T[[t]]$
rather than on $G((t))/N^-((t))$.}

In general, we shall see that property (3) is equivalent
to imposing the condition that either *- or !-restriction of
$\fp^{\cmu,-}(\CF)|_{_{\leq \cnu}^k\CZ^\cmu}$ to $^k_{\cmu}\fF^\cmu$ is
$T[[t]]$-equivariant.

\bigskip

In the sequel we will formulate a conjecture, from which it follows that
the category $\Pervsmk$ is independent of the curve $X$, and possesses the symmetries
expected from "the category of $G^k$-equivariant sheaves on
${\mathcal Fl}^{\frac{\infty}{2}}:=G((t))/N^-((t))\cdot T[[t]]$",
in particular, it will carry an action of the lattice $\cLambda\simeq T((t))/T[[t]]$ by translation
functors.

\sssec{}

Our present goal is to describe the irreducibles in $\Pervsmk$. Recall the
isomorphism \eqref{descr open stratum}, which realizes $_\cnu^k\BunNb$
as a fibration over the base $G/N^-([t]/t^k)$ with typical fiber $_\cnu^k\CN$.
(In fact, $_\cnu^k\BunNb$ is a principal $N^-([t]/t^k)$-bundle over the product
$G/N^-([t]/t^k)\times {}_\cnu\BunNb$.)

In particular, for a perverse sheaf $\CF'$ on $\bigl(G/N^-\bigr)([t]/t^k)$, we can form
the twisted external product
$$\CF'\tboxtimes \IC_{_\cnu^k\CN}\in \Perv({}_\cnu^k\BunNb).$$
Up to a cohomological shift, it is isomorphic to the pull-back of
$$\CF'\boxtimes \IC_{_\cnu\BunNb}\in \Perv(G/N^-([t]/t^k)\times {}_\cnu\BunNb).$$

\begin{prop}  \label{sheaves on strata} \hfill

\smallskip

\noindent{\em (1)}
For $\CF\in \Pervsmk$, all perverse cohomologies of the restriction
$\oi_\cnu^*(\CF)$ are of the form $\CF'\tboxtimes \IC_{_\cnu^k\CN}$,
where $\CF'$ is a perverse sheaf on $\bigl(G/N^-\bigr)([t]/t^k)$, that comes
as a pull-back from a perverse sheaf on $\bigl(G/B^-\bigr)([t]/t^k)$.

\smallskip

\noindent{\em (2)}
The perverse sheaf (resp., each perverse cohomology of)
$(\oi_\cnu)_{!*}(\CF'\tboxtimes \IC_{_\cnu^k\CN})$
(resp., $(\oi_\cnu)_{!}(\CF'\tboxtimes \IC_{_\cnu^k\CN})$)
for $\CF'$ as above is an object of $\Pervsmk$.

\smallskip

\noindent{\em (3)}
Perverse sheaves of the form $(\oi_\cnu)_{!*}(\CF'\tboxtimes \IC_{_\cnu^k\CN})$
for $\CF'$ as above are all the irreducible objects of $\Pervsmk$.
\end{prop}

The rest of this subsection is devoted to the proof of the proposition. Note, however,
that point (3) is a formal corollary of points (1) and (2).

\bigskip

The factorization isomorphisms of \lemref{factorization} respect the substacks
$_\cnu^k\oCZ^\cmu$, $_\cnu^k\CZ^\cmu$, $_{\leq \cnu}^k\CZ^\cmu$ of
$_{\infty}^k\CZ^\cmu$. Hence, it makes sense to introduce the category $_\cnu'\Pervsmk$,
which is a full subcategory of $\Perv\Bigl({}_\cnu^k\BunNb\Bigr)$, consisting
of objects, satisfying the same conditions (1), (2) and (3) as in the definition of
$\Pervsmk$.

\medskip

It is clear that for $\CF_1\in \Pervsmk$, the perverse cohomologies of the restriction
$\oi_\cnu^*(\CF_1)$ are objects of $_\cnu'\Pervsmk$, and vice versa: for
$\CF_2\in {}_\nu'\Pervsmk$, the perverse sheaf (resp., each perverse cohomology of)
$(\oi_\cnu)_{!*}(\CF_2)$ (resp., $(\oi_\cnu)_{!}(\CF'_2)$) belongs to
$\Pervsmk$.

\medskip

Therefore, the assertion of the proposition reduces to showing that the functor
$\CF'\mapsto \CF'\tboxtimes \IC_{_\cnu\BunNb}$ defines an equivalence
$$\Perv\Bigl(\bigl(G/B^-\bigr)([t]/t^k)\Bigr)\to {}_\nu'\Pervsmk.$$

\medskip

First, we claim that every object of $_\nu'\Pervsmk$ is the Goresky-MacPherson
extension of its restriction to the open sub-stack $_\cnu\Bun_N$. Indeed, if
it were not, we would be able to find $\cmu_1$ and $\cmu_2$ large enough,
so that either $!$ or $*$-restriction of $\CF$ to the closed sub-stack
$$\Bigl(\bigl(\CZ^{\cmu_1}-\oCZ^{\cmu_1} \bigr)\times
{}_\infty^k\CZ^{\cmu_2}\Bigr)
\underset{X^{\cmu_1}\times {}_{\infty}X^{\cmu_1}}\times
\Bigl(\oX^{\cmu_2}\times {}_{\infty}X^{\cmu_2}\Bigr)_{disj}$$
would have non-zero perverse cohomologies in positive (resp., negative)
degrees. However, this contradicts the factorizability property (2).

\medskip

Let us denote by $_\cnu\Pervsmk$ the corresponding full subcategory
of $\Perv\Bigl({}_\cnu^k\BunN\Bigr)$ consisting of perverse sheaves, satisfying
(1) and (3). We are reduced to showing that
$$\CF'\mapsto \CF'\tboxtimes \IC_{_\cnu\CN^k}:
\Perv\Bigl(\bigl(G/B^-\bigr)([t]/t^k)\Bigr)\to {}_\nu\Pervsmk$$
is an equivalence. Note that the latter functor is isomorphic,
up to a cohomological shift, to the pull-back functor under the smooth map
\begin{equation} \label{eval map}
_\cnu^k\BunN\overset{\on{ev}_\cnu}\longrightarrow \bigl(G/N^-\bigr)([t]/t^k)\to
\bigl(G/B^-\bigr)([t]/t^k).
\end{equation}

\bigskip

The fact that the functor in question is fully faithful is clear, since
the map in \eqref{eval map} has connected fibers. Hence, it remains to
show the essential surjectivity.

First, let us show that any $\CF\in \Perv\Bigl({}_\cnu^k\BunN\Bigr)$
is the pull-back under $\on{ev}_\cnu$ of some perverse sheaf
$\CF'$ on $\bigl(G/N\bigr)([t]/t^k)$.

For any non-empty collection of points $\ol{x}'$, distinct from $x$,
consider the pull-back of $\CF$ to
$_\cnu^k\BunN\underset{_\cnu^k\BunNb}\times {}^k\CN^{\ol{x}'}$.
By property (1), it is equivariant with respect to the group-indscheme
$\underset{j=1,...,m}\Pi\, N^-((t'_j))$.

This implies our assertion, since the above group-indscheme acts transitively
along the fibers of the composed map
\begin{equation} \label{aux5}
_\cnu^k\BunN\underset{_\cnu^k\BunNb}\times {}^k\CN^{\ol{x}'}\to
\bigl(G/N^-\bigr)([t]/t^k).
\end{equation}

\medskip

Thus, it remains to show that condition (3) on $\CF$ implies that
the perverse sheaf $\CF'$ on $\bigl(G/N^-\bigr)([t]/t^k)$
comes as a pull-back from a perverse sheaf on $\bigl(G/B^-\bigr)([t]/t^k)$.
I.e., we have to show that $\CF'$ is equivariant with respect to
$T([t]/t^k)$. Note that the equivariance with respect to the subgroup
$T\subset T([t]/t^k)$ follows from the assumption that $\CF$ on
$_\cnu^k\BunN$ was $T$-equivariant. Thus, it remains to check the
equivarince property with respect to the unipotent subgroup
$\on{ker}\bigl(T([t]/t^k)\to T\bigr)$.

\medskip

For $\cmu$ such that $\cmu+\cnu\in \cLambda^{pos}$, consider the
composed map
\begin{equation} \label{aux 4}
_\cnu^k\ofF^\cmu\to \bigl(G/N^-\bigr)([t]/t^k),
\end{equation}
where $_\cnu^k\ofF^\cmu$ is the fiber of $^k_{\cnu}\oCZ^\cmu$ over
$\cmu\cdot x\in {}_{\infty}X^\cmu$. The above map is equivariant
with respect to $T[[t]]$ acting on the two sides. Moreover, it
is surjective if $\cmu$ was chosen large enough.

\medskip

Let $k'\geq k$ be such that the action of $T[[t]]$ on $_\cnu^k\ofF^\cmu$
factors through $T([t]/t^{k'})$. Let $\ol{m}$ be large enough, so that the map
$\on{Maps}(X,T)^{\ol{m}}\to T([t]/t^{k'})$, given by Taylor expansion at $x$,
is surjective.

Property (3) for this $\ol{m}$ implies then that the restriction of $\CF$ to
$_\cnu^k\ofF^\cmu$ is $T([t]/t^{k'})$-equivariant. This implies that
$\CF'$ is also equivariant with respect to this group.

\sssec{}

We will now investigate the mutual dependence of conditions (1) and (2). For
a natural number $m$ consider the product $\oX^m\times {}_{\infty}^k\BunNb$,
and let $\Bigl(\oX^m\times {}_{\infty}^k\BunNb\Bigr)^{n.z.}$ denote the open
subset, corresponding to the condition that the zeros of the maps $\kappa^\lambda$
are away from the $m$ marked points of $X^m$. In other words, the fiber of this
space over a given $\ol{x}'\in \oX^m$ is the stack that we denoted by
$_{\infty}^k\BunNb^{\text{n.z.}\ol{x}'}$.

Over $X^m$ we have a group-scheme, denoted $\on{Jets}^+(N^-)^m$, whose fiver over
$\ol{x}'=\{x'_1,...,x'_m\}$ is $\Pi\, N^-[[t'_j]]$, where the product is taken over distinct
points among the $x'_i$'s. In addition, we have a group-indscheme, denoted
$\on{Jets}(N^-)^m$, whose fiber over the same collection of points is $\Pi\, N^-((t'_j))$.
Since $N^-$ is unipotent, this group-indscheme can be represented as a union
of its closed group-subschemes.

Finally, over $\Bigl(X^m\times {}_{\infty}^k\BunNb\Bigr)^{n.z.}$ there exists
a canonical $\on{Jets}^+(N^-)^m$-torsor, which we will denote by $^k\CN^m$.
The action of $\on{Jets}^+(N^-)^m$ on $^k\CN^m$ extends to an action of
$\on{Jets}(N^-)^m$.

\begin{lem}  \label{stronger equivariance}
Let $\CF$ be a perverse sheaf on $_\infty^k\BunNb$,
which satisfies property (1) of \secref{category of sheaves}.
Then the pull-back of $\CF$ to $^k\CN^m$ is equivariant with respect to
$\on{Jets}(N^-)^m$.
\end{lem}

This follows from the fact "fiber-wise equivariance" implies "equivariance"
for a unipotent group-scheme.

\medskip

\noindent{\it Remark.} Arguing as in \cite{FGV}, Sect. 6.2, one can show that condition (1)
is equivalent to the following, seemingly weaker, condition. Namely, it is
sufficient to impose the $N^-((t'))$-equivariance condition for just one
fixed point $x'$ distinct from $x$.

\medskip

Let us say that a perverse sheaf $\CF$ on $_\infty^k\BunNb$ has a weak factorization
property if, in the notation of \secref{category of sheaves},  the isomorphism between
$\fp^{-,\cmu_2}{}^*(\CF)$ and $\IC_{\CZ^{\cmu_2-\cmu_1}}\boxtimes\, \fp^{-,\cmu_1}{}^*(\CF)$
holds over the open subset
$$\Bigl(\oCZ^{\cmu_2-\cmu_1} \times {}_\infty^k\CZ^{\cmu_1}\Bigr)
\underset{X^{\cmu_2-\cmu_1}\times {}_{\infty}X^{\cmu_1}}\times
\Bigl(\oX^{\cmu_2-\cmu_1}\times {}_{\infty}X^{\cmu_1}\Bigr)_{disj}.$$

Since $\oCZ^{\cmu_2-\cmu_1}$ is smooth, this condition is equivalent to
the restriction of $\fp^{-,\cmu_2}{}^*(\CF)$ to the above open subset being
constant along the first factor.

\begin{prop}  \label{equiv of (1) and (2)}
For a perverse sheaf $\CF$ on $_\infty^k\BunNb$, property (1) is equivalent
to the weak factorization property.
\end{prop}

Before giving a proof let us make the following observation: we have two maps
$$\hl_{N^-},\hr_{N^-}:\on{Jets}(N^-)^m\overset{\on{Jets}^+(N^-)^m}\times {}^k\CN^m\to
{}_\infty^k\BunNb,$$
the first being the tautological projection, and the second is given by the action of
$\on{Jets}(N^-)^m$ on $^k\CN^m$.  If $\cmu_1,\cmu_2\in \cLambda$ are two elements,
with $\cmu_2-\cmu_1\in \cLambda^{pos}$ such that $m=l(\cmu_2-\cmu_1)$
there is a natural projection $\oX^m\to \oX^{\cmu_2-\cmu_1}$ and a map
\begin{align}  \label{Zastava vs. groupoid}
&\Bigl(\oCZ^{\cmu_2-\cmu_1} \times {}_\infty^k\CZ^{\cmu_1}\Bigr)
\underset{X^{\cmu_2-\cmu_1}\times {}_{\infty}X^{\cmu_1}}\times
\Bigl(\oX^{\cmu_2-\cmu_1}\times {}_{\infty}X^{\cmu_1}\Bigr)_{disj}
\underset{\oX^{\cmu_2-\cmu_1}}\times \oX^m \to \\
&\on{Jets}(N)^m\overset{\on{Jets}^+(N)^m}\times {}^k\CN^m,
\end{align}
such that its composition with $\hl_N$ is the projection
$$\Bigl(\oCZ^{\cmu_2-\cmu_1} \times {}_\infty^k\CZ^{\cmu_1}\Bigr)
\underset{X^{\cmu_2-\cmu_1}\times {}_{\infty}X^{\cmu_1}}\times
\Bigl(\oX^{\cmu_2-\cmu_1}\times {}_{\infty}X^{\cmu_1}\Bigr)_{disj}
\underset{\oX^{\cmu_2-\cmu_1}}\times \oX^m \to
{}_\infty^k\CZ^{\cmu_1}\to {}_\infty^k\BunNb,$$
and its composition with $\hr_N$ identifies via \lemref{factorization} with
$$_\infty^k\CZ^{\cmu_2} \underset{_{\infty}X^{\cmu_2}}\times
\Bigl(\oX^{\cmu_2-\cmu_1}\times {}_{\infty}X^{\cmu_1}\Bigr)_{disj}
\underset{\oX^{\cmu_2-\cmu_1}}\times \oX^m\to {}_\infty^k\CZ^{\cmu_2} \to
{}_\infty^k\BunNb.$$

Now let us proof the proposition:

\begin{proof}

Assume first that $\CF$ satisfies property (1), and hence, by
\propref{stronger equivariance}, its pull-back to $^k\CN^m$ is
$\on{Jets}(N^-)^m$-equivariant. We obtain that the restrictions of
$\hl^*_{N^-}(\CF)$ and $\hr^*_{N^-}(\CF)$ to any finite-dimensional subscheme of
$\on{Jets}(N^-)^m\overset{\on{Jets}^+(N^-)^m}\times {}^k\CN^m$ are isomorphic.
Then the weak factorizability of $\CF$ follows from the properties of the map
from \eqref{Zastava vs. groupoid} above.

\medskip

To prove the implication in the opposite direction, we reverse the steps.
We have to show that for a given finite collection of distinct points $\ol{x}'=\{x'_1,...,x'_m\}$,
the restrictions of $\hl_{N^-}$ and $\hr_{N^-}$ to the fiber of
$\on{Jets}(N^-)^m\overset{\on{Jets}^+(N^-)^m}\times {}^k\CN^m$ over
$\ol{x}'\in X^m$ are isomorphic over every finite-dimensional subscheme
of this ind-scheme.  Since each $N^-((t'_j))$ is a union of pro-unipotent
subgroups, it is sufficient to show that the isomorphism holds after
the base change with respect to
$$\Bigl({}_\infty^k\CZ^\cmu\underset{_\infty{X}^\cmu}\times {}_{\infty}(X-\ol{x}')^\cmu\Bigr)
\to {}_\infty^k\BunNb$$
for $\cmu$ large enough.

Note that the above fiber, base-changed to $_\infty^k\CZ^\cmu$, is isomorphic
to
$$\Bigl({}_\infty^k\CZ^\cmu\underset{_\infty{X}^\cmu}\times {}_{\infty}(X-\ol{x}')^\cmu\Bigr)
\times \underset{j}\Pi\, \Bigl(N^-((t'_j))\cdot 1_{\Gr}\Bigr).$$

Our assertion follows now from \eqref{non-central fibers}, since
$N^-((t'_j))\cdot 1_{\Gr}$ can be exhausted
by affine subspaces, each of which contains as a dense subset the intersection
$$\left(N((t'_j))\cdot \cmu'\right)\cap \left(N^-((t'_j))\cdot 1_{\Gr}\right)$$
for some $\cmu'$.

\end{proof}

As a corollary of the first assertion of the proposition, we obtain the following:

\begin{cor}  \label{1,2 stable under extensions}
Let $0\to \CF_1\to \CF\to \CF_2\to 0$ be a short exact sequence of objects
of $\Perv({}_\infty^k\BunNb)^T$, with $\CF_1,\CF_2$ satisfying properties
(1) and (2) of the definition of $\Pervsmk$. Then $\CF$ also satisfies
properties (1) and (2).
\end{cor}

\begin{proof}

Since the group $N^-((t'))$ is (ind)-pro-unipotent,
the only non-trivial condition to check is the factorizability property. For
$\CF$ as above, its pull-back to
$$_\infty^k\CZ^{\cmu_2} \underset{_{\infty}X^{\cmu_2}}\times
\Bigl(\oX^{\cmu_2-\cmu_1}\times {}_{\infty}X^{\cmu_1}\Bigr)_{disj}\simeq
\Bigl(\CZ^{\cmu_2-\cmu_1} \times {}_\infty^k\CZ^{\cmu_1}\Bigr)
\underset{X^{\cmu_2-\cmu_1}\times {}_{\infty}X^{\cmu_1}}\times
\Bigl(\oX^{\cmu_2-\cmu_1}\times {}_{\infty}X^{\cmu_1}\Bigr)_{disj}$$
is the Goresky-MacPherson extension from the open sub-space
$$\Bigl(\oCZ^{\cmu_2-\cmu_1} \times {}_\infty^k\CZ^{\cmu_1}\Bigr)
\underset{X^{\cmu_2-\cmu_1}\times {}_{\infty}X^{\cmu_1}}\times
\Bigl(\oX^{\cmu_2-\cmu_1}\times {}_{\infty}X^{\cmu_1}\Bigr)_{disj}.$$

However, the latter is constant along the $\oCZ^{\cmu_2-\cmu_1}$-factor
because of property (1), \lemref{stronger equivariance} and
\propref{equiv of (1) and (2)}. Along the $_\infty^k\CZ^{\cmu_1}$
factor it is isomorphic to $\fp^{-,\cmu_1}{}^*(\CF)$ by
\eqref{Zastava vs. groupoid}.

\end{proof}

\sssec{}

Our present goal is to establish the following:

\begin{prop}  \label{stable under extensions}
The category $\Pervsmk$, as a subcategory of the category of $T$-equivariant
perverse sheaves on $_\infty^k\BunNb$, is stable under extensions.
\end{prop}

The rest of the present subsection is devoted to this proposition. In view of
\corref{1,2 stable under extensions}, we have to show that if
$0\to \CF_1\to \CF\to \CF_2\to 0$ is a short exact sequence in $\Perv({}_\infty^k\BunNb)^T$
with $\CF_1,\CF_2\in \Pervsmk$, then $\CF$ satisfies property (3).

Consider the pull-back
$$\on{act}^*_T\left((\fp^{-,\cmu})^*(\CF)\right)\in \Perv\Bigl(
\left(\on{Maps}(X,T)^{\ol{m}}\times {}
_{\leq \cnu}^k\CZ^\cmu\right)\underset{\oX^{\ol{m}}\times {}_{\leq \cnu}X^\cmu}\times
\left(\oX^{\ol{m}}\times {}_{\leq \cnu}X^\cmu\right)_{disj}\Bigr).$$

Since $\on{Maps}(X,T)^{\ol{m}}$ is isomorphic to the affine space, it is sufficient
to show that the restriction of the above pull-back to the fiber over every
geometric point $\bz:=(\CP_G,\{\kappa^\lambda\},\{\kappa^{\lambda,-}\})\in {}_\infty^k\CZ^\cmu$
is a complex with constant cohomologies.

By the factorization property, it is sufficient to consider the case when the point
$^k_\infty\fs^\cmu(\bz)\in {}_\infty X^\cmu$ equals $\cmu\cdot x$. In this case,
the map $\on{act}^*_T$ factors through the action of
$\on{ker}\Bigl(T[[t]]\to T\Bigr)$ on $^k_\infty\fF^\cmu$.

Hence, it is sufficient to check that the restriction of $\CF$ to $^k_\infty\fF^\cmu$
is $\on{ker}\Bigl(T[[t]]\to T\Bigr)$-equivariant. But the above restriction is an extension of the restrictions
of $\CF_1$ and $\CF_2$. Since for $\ol{m}'$ large enough the map
$\on{Maps}(X,T)^{\ol{m}'}\to \on{ker}\Bigl(T[[t]]\to T\Bigr)$ is surjective with connected fibers,
the fact that $\CF_1$ and $\CF_2$ satisfy property (3) implies that their
restrictions to $^k_\infty\fF^\cmu$ are $\on{ker}\Bigl(T[[t]]\to T\Bigr)$-equivariant. This
proves our assertion, since $\on{ker}\Bigl(T[[t]]\to T\Bigr)$ is pro-unipotent, and hence
the equivariance property is stable under extensions.

\bigskip

Thus, \propref{stable under extensions} is proved. As a by-product we obtain the following
alternative way to spell out condition (3):

\begin{cor}
Let $\CF\in \Perv({}_\infty^k\BunNb)^T$ be a perverse sheaf, satisfying properties
(1) and (2) from the definition of $\Pervsmk$. The the following are equivalent:

\smallskip

\noindent{\em (1)} $\CF$ satisifies also property (3).

\smallskip

\noindent{\em (2)} The *- (or !-) restrictions of $\CF$ to every $_\cnu^k\Bun_{N^-}$
are such that their perverse cohomologies are pull-backs from $T([t]/t^k)$-equivariant
perverse sheaves on $G/N^-([t]/t^k)$.

\smallskip

\noindent{\em (3)} The *- (or !-) restrictions of $(\fp^{\cmu,-})^*(\CF)$ to every
$^k_\infty\fF^\cmu$ is $T[[t]]$-equivariant.

\end{cor}

\sssec{}

Recall that for $\cnu\in \cLambda$, we have introduced the category
$$_\cnu'\Pervsmk\subset \Perv({}_\cnu^k\BunNb),$$
which is equivalent to
$${}_\cnu\Pervsmk\subset \Perv({}_\cnu^k\Bun_N).$$

\begin{prop} \label{local exts}
Let $\CF'\in {}_\cnu'\Pervsmk$ be such that $(\oi_\cnu)_!(\CF')$ is a
perverse sheaf. Then for $\CF\in \Pervsmk$ the canonical map
$$Ext^1_{_\infty^k\BunNb}\Bigl((\oi_\cnu)_!(\CF'),\CF\Bigr)\to
Ext^1_{_\infty^k\BunNb}\Bigl((i_\cnu)_!(\CF'),\CF\Bigr)$$
is an isomorphism.
\end{prop}

Note that due to \propref{stable under extensions},
the above proposition can be reformulated as follows:
$$Ext^1_{\Pervsmk}\Bigl((\oi_\cnu)_!(\CF'),\CF\Bigr)\simeq
\on{R}^1\on{Hom}_{_\cnu^k\Bun_N}\Bigl(\CF',i_\cnu^!(\CF)\Bigr).$$

\begin{proof}
The fact that the map in question is injective is evident, since
$(i_\cnu)_!(\CF')$ surjects onto $(\oi_\cnu)_!(\CF')$, and $\CF$
has no sub-objects supported on $_{\cnu}^k\BunNb-{}_{\cnu}^k\BunN$.

\medskip

To prove the surjectivity we can replace $_{\infty}^k\BunNb$ by its
open sub-stack $_{\geq \cnu}^k\BunNb$, which is obtained by removing
from $_\infty^k\BunNb$ all $_{\leq \cnu'}^k\BunNb$
for $\cnu'<\cnu$. Evidently, $_\cnu^k\BunNb$ is closed in
$_{\geq \cnu}^k\BunNb$.

\medskip

Let
$$0\to \CF\to \CF_1\to (i_\cnu)_!(\CF')\to 0$$
be an extension. We have to show that it is induced from
an extension of $(\oi_\cnu)_!(\CF')$ by $\CF$.
Let $\tilde{\CF}_1$ be the perverse sheaf
on $_{\geq \cnu}^k\BunNb$ obtained as a Goresky-MacPherson
of the restriction of $\CF_1$ to the open substack
$$_{\geq \cnu}^k\BunNb-({}_\cnu^k\BunNb-{}_\cnu^k\BunN).$$
We claim that $\tilde{\CF}_1$ is the desired extension. Namely,
we have the maps
$$\CF\hookrightarrow \tilde{\CF}_1\twoheadrightarrow (\oi_\cnu)_!(\CF'),$$
and we claim that this is a short exact sequence.

To check this, by \propref{sheaves on strata}(3), it is enough
to show that $\tilde{\CF}_1$ is an object of the corresponding category
$_{\geq \cnu}\Pervsmk$. However, properties (1) and (3) are
automatic, and the factorization property (2) follows by
combining \propref{equiv of (1) and (2)} and the definition of
Goresky-MacPherson extension.

\end{proof}

The 5-lemma yields:
\begin{cor} \label{Ext 2 injects}
For $\CF$ as in the proposition, the natural map
$$Ext^2_{\Pervsmk}\Bigl((\oi_\cnu)_!(\CF'),\CF\Bigr)\to
\on{R}^2\on{Hom}_{_\cnu^k\Bun_N}\Bigl(\CF',i_\cnu^!(\CF)\Bigr)$$
is injective.
\end{cor}

\noindent{\it Remark.}
>From \propref{local exts} one can formally deduce that the maps
$$Ext^i_{\Pervsmk}\Bigl((\oi_\cnu)_!(\CF'),\CF\Bigr)\to
\on{R}^i\on{Hom}_{_\cnu^k\Bun_N}\Bigl(\CF',i_\cnu^!(\CF)\Bigr)$$
are isomorphisms for all $i$.

\ssec{The spherical case}

\sssec{}

Let $\Pervsmg$ denote $\Pervsmk$ for $k=0$; this is a full subcategory
in $_\infty\BunNb$. For $\cnu\in \cLambda$
we will denote by $\fIC_\cnu$ the corresponding irreducible, i.e.,
$$\fIC_\cnu\simeq (i_\cnu)_{!*}(\IC_{{}_\cnu\Bun_N})\simeq
(\oi_\cnu)_{!*}(\IC_{{}_\cnu\BunNb}).$$
These are the irreducible objects of $\Pervsmg$.

\begin{prop}  \label{semi-simplicity of sph}
The category $\Pervsmg$ is semi-simple.
\end{prop}

\begin{proof}

It would be enough to show that if $\fIC_{\cnu_1}$ and $\fIC_{\cnu_2}$
are two simple objects of $\Pervsmg$, whose support is contained
in some $_{\leq\cnu}\BunNb$, then over some open substack of
$_{\leq \cnu}\BunNb$, $Ext^1(\fIC_{\cnu_1},\fIC_{\cnu_2})$ is zero.

Let $\cnu_1,\cnu_2$ be two elements of $\cLambda$. In order for
$Ext^1(\fIC_{\cnu_1},\fIC_{\cnu_2})$ to be non-trivial, the support
of one sheaf must be contained in the closure of the support of the other.
This means that either $\cnu_1\leq \cnu_2$ or $\cnu_2\leq \cnu_1$.
By Verdier duality we can assume that $\cnu_1\leq \cnu_2$.

Consider the open sub-stack of  $_{\leq \cnu_2}\BunNb$ obtained
by removing the closed sub-stack $_{\leq \cnu_1}\BunNb-{}_{\cnu_1}\Bun_N$.
As in \cite{FGV}, Sect. 6.1.4,
$$Ext^1_{{}_{\leq \cnu_2}\BunNb}(\fIC_{\cnu_1},\fIC_{\cnu_2})\hookrightarrow
Ext^1_{{}_{\leq \cnu_2}\BunNb-({}_{\leq \cnu_1}\BunNb-{}_{\cnu_1}\Bun_N)}
(\fIC_{\cnu_1},\fIC_{\cnu_2}),$$
so it is enough to show that the latter is $0$. Since
$$_{\cnu_1}\Bun_N\subset {}_{\leq \cnu_2}\BunNb-({}_{\leq \cnu_1}\BunNb-{}_{\cnu_1}\Bun_N)$$
is closed, the latter $\on{Ext}^1$ is isomorphic to
$$\on{R}^1\on{Hom}{}_{_{\cnu_1}\Bun_N}\left(\IC{}_{_{\cnu_1}\Bun_N},
i_{\cnu_1}^!(\IC{}_{_{\cnu_2}\Bun_N})\right).$$

There are two cases: if $\cnu_1<\cnu_2$, then we are done by \cite{BFGM}, since
$i_{\cnu_1}^!(\IC_{_{\cnu_2}\Bun_N})$ lives in the cohomological degrees
$\geq 2$.

If $\cnu=\cnu_1=\cnu_2$, then the assertion follows from the fact that $_\cnu\Bun_N$
is simply-connected, cf. \cite{FGV}, Sect. 6.

\end{proof}

\sssec{}

Consider the object of ${\mathsf D}(_{\infty}\BunNb)$ equal to
$(\oi_{\cnu})_!(\IC_{_{\cnu}\BunNb})$. This is a complex that lives in
non-positive cohomological degrees, and each of its perverse
cohomologies is an object of $\Pervsmg$, by \propref{sheaves on strata}.

\begin{thm}   \label{naive extension}
The $-k$-th perverse cohomology of $(\oi_{\cnu})_!(\IC_{_{\cnu}\BunNb})$
is isomorphic to the direct sum over collections of $k$ distinct positive roots
$\{\beta_1,...,\beta_k\}$ of
$$\fIC_{\cnu-\underset{j}\Sigma\, \beta_j}.$$
\end{thm}

\begin{cor}  \label{top naive extension}
The complex $(\oi_{\cnu})_!(\IC_{_{\cnu}\BunNb})$ (resp., $(\oi_{\cnu})_*(\IC_{_{\cnu}\BunNb})$)
lives in the cohomological degrees $[-\dim(\fn),0]$ (resp., $[0,\dim(\fn)]$) and its
$-\dim(\fn)$- ($\dim(\fn)$-) degree cohomology is isomorphic to
$\fIC_{\cnu-2\crho}$.
\end{cor}

The rest of this subsection is devoted to the proof of the above theorem. For
$\cmu\in \cLambda^{pos}$ consider the stack $_{\leq \cnu}\BunNb^{\leq \cmu}$,
fibered over $X^\cmu$, classifying pairs $(D\in X^\cmu,\{\kappa^\lambda\})$
such that each $\kappa^\lambda$ factors as
$$V^\lambda_{\CP_G}\to \CO_X\Bigl(\lambda(\cnu\cdot x-D)\Bigr)\to \CO_X.$$
Let $_\cnu\BunN^\cmu$ be the open sub-stack of $_{\leq \cnu}\BunNb^{\leq \cmu}$,
corresponding to the condition that the maps
$V^\lambda_{\CP_G}\to \CO_X\Bigl(\lambda(\cnu\cdot x-D)\Bigr)$ above, are bundle maps.

It is easy to see that $_\cnu\BunN^\cmu$ is smooth over $X^\cmu$.
The projection $_{\leq \cnu}\BunNb^{\leq \cmu}\to X^\cmu$ is
ULA (universally locally acyclic) with respect to the IC sheaf on this stack,
by \cite{BG}, Sect. 5.2.

We let $\oi^{\leq \cmu}$ (resp., $i^\cmu$) denote the natural maps
from the above stacks to $_{\leq \cnu}\BunNb$. By \cite{BG},
$\oi^{\leq \cmu}$ is finite (and, in particular, proper), and $i^\cmu$ is a locally closed
embedding. Moreover, by \cite{FGV}, Sect. 3.3, $i^\cmu$ is affine. In particular, every
$i^\cmu_!(\IC_{_\cnu\BunN^\cmu})$ is a perverse sheaf.

\medskip

The following is a reformulation of the main result of \cite{FFKM} and \cite{BFGM}:

\begin{thm} \label{stalks of IC}
The $k$-th cohomology of $(i^\cmu)^*(\fIC_\cnu)$ is isomorphic to the direct sum
over the set of partitions $\fP$
$$\cmu=\Sigma\, m_j\cdot \beta^j,\, \beta^j\neq \beta^{j'},\, \Sigma\, m_j=k,$$
where $\beta^i$'s are positive roots,
of the direct images of the shifted by $[k]$ constant perverse sheaves on each
$$X^\fP\underset{X^\cmu}\times {}_\cnu\BunN^\cmu,$$
where $X^\fP\simeq \underset{j}\Pi\, X^{(m_j)}$, that maps naturally to $X^\cmu$.
\end{thm}

For each partition $\fP$ as above let $\CE_\fP$ be the perverse sheaf on
$X^\cmu$, equal to the direct image under $X^\fP\to X^\cmu$ of the irreducible
perverse sheaf obtained by taking the external product over $j$ of the $1$-dimensional
local systems on each $X^{(m_j)}-\text{Diag}$ with monodromy $-1$ around the diagonal.
By the ULA property of $_{\leq \cnu}\BunNb^{\leq \cmu}$ over $X^\cmu$, the tensor
product $\IC_{_{\leq \cnu}\BunNb^{\leq \cmu}}\otimes \CE_\fP[-k]$ is a perverse
sheaf.

The usual Koszul complex argument yields the following:

\begin{cor}  \label{spherical const, global}
Irreducible constituents of $(i_\nu)_!(\IC_{_\cnu\Bun_N})$ are the perverse
sheaves $\oi^{\leq \cmu}_*\Bigl(\IC_{_{\leq \cnu}\BunNb^{\leq \cmu}}\otimes \CE_\fP[-k]\Bigr)$
for all $\cmu\in \cLambda$ and partitions $\fP$, each appearing once.
\end{cor}

\medskip

Recall that $_{\leq \cnu}\fF^\cmu$ denotes the fiber of $_{\leq \cnu}\CZ^\cmu$ over
$\cmu\cdot x\in {}{_\infty}X^\cmu$. By \cite{BFGM}, we have:
\begin{equation} \label{cohomology of central fiber}
H_c\Bigl({}_{\leq \cnu}\fF^\cmu,\IC_{_{\leq \cnu}\CZ^\cmu}|_{{}_{\leq \cnu}\fF^\cmu}\Bigr)\simeq
U(\check\fn)_{\cmu+\cnu},
\end{equation}
in particular, the above cohomology is concentrated in cohomological degree $0$.

Combining this result with \corref{spherical const, global}, and taking into account that
the restriction of $\CE_\fP$ to the diagonal divisor is $0$ unless all $m_j=1$,
we obtain the following:
\begin{cor} \label{cohomology of open fibers}
The cohomology group
$$H^{-k}_c\Bigl({}_{\leq \cnu}\fF^\cmu,(i_\nu)_!(\IC_{_\cnu\oCZ^\cmu})|_{{}_{\leq \cnu}\fF^\cmu}\Bigr)$$
is the direct sum over $\clambda\in \cLambda^{pos}$ of $U(\check\fn)_{\cmu+\cnu-\clambda}$,
each appearing the number of times equal to the number of partitions of $\clambda$
as a sum of $k$ distinct positive roots.
\end{cor}

Let us note that the intersection $_\cnu\ofF^\cmu:={}_{\infty}\fF^\cmu\cap {}_\cnu\oCZ^\cmu$
is isomorphic to
$$\left(N((t))\cdot (\cnu+\cmu)\right)\cap \left(N^-((t))\cdot 1_{\Gr}\right).$$
Thus, \corref{cohomology of open fibers} gives an expression for
\begin{equation} \label{cohom open fib}
H^{-k+\langle 2\rho,\cmu+\cnu\rangle}_c\Bigl(\left(N((t))\cdot \cnu\right)\cap \left(N^-((t))\cdot
(-\cmu)\right),\uBC\Bigr)\simeq
H^{-k}_c\Bigl({}_{\leq \cnu}\fF^\cmu,(\oi_\nu)_!(\IC_{_\cnu\CZ^\cmu})|_{{}_{\leq \cnu}\fF^\cmu}\Bigr).
\end{equation}

\medskip

Now we can finish the proof of \thmref{naive extension}, essentially be reversing
the logic. We have to show that the
multiplicity $m_k(\clambda)$ of $\IC_{\cnu-\clambda}$ in the $-k$-th
perverse cohomology of $(\oi_\cnu)_!(\IC_{{}_\cnu\BunNb})$ equals the number of partitions
of $\clambda$ as a sum of $k$ distinct positive roots, i.e., $\dim(\Lambda^k(\fn)_\clambda)$.

We will argue by induction on $\clambda$, so we can assume that the assertion is
known for all $\clambda'<\clambda$. Consider the cohomology in \eqref{cohom open fib}
for $\cmu=\clambda-\cnu$. By \eqref{cohomology of central fiber}, the contributions of
different constituents do not cancel out, and we obtain an equality:
$$\underset{\clambda'\in \Lambda^{pos}}\Sigma\, \dim(\Lambda^k(\fn)_{\clambda'})\cdot
\dim(U(\fn)_{\clambda-\clambda'})=\underset{\clambda'<\lambda}\Sigma\,
\dim(\Lambda^k(\fn)_{\clambda'})\cdot
\dim(U(\fn)_{\clambda-\clambda'})+m_k(\lambda).$$
This implies the desired equality.

\ssec{The Iwahori case}

\sssec{}    \label{naive convolution}

Note that the stack $^k_{\infty}\BunNb$ is acted on by the group $G([t]/t^k)$. In particular, we have the
convolution functors:
$$\sD(G([t]/t^k))\times \sD({}^k_{\infty}\BunNb)\to \sD({}^k_{\infty}\BunNb):
\CS,\CF\mapsto \CS\starstar\CF \text{ and } \CS,\CF\mapsto \CS\starshriek\CF.$$
Moreover, these functors are defined on each of the subcategories
$$\sD({}^k_{\leq \cnu}\BunNb), \sD({}^k_{\cnu}\BunNb) \text{ and }\sD({}^k_{\cnu}\BunN),$$
so that the *-convolution commutes in the natural sense with the functors
$$(\oi_{\leq \cnu})_*, (\oi_{\cnu})_*, (i_{\cnu})_*,(\oi_{\leq \cnu})^!, (\oi_{\cnu})^!, (i_{\cnu})^!,
(\on{ev}_\cnu)^!$$
and the !-convolution commutes with the functors
$$(\oi_{\leq \cnu})_!=(\oi_{\leq \cnu})_*, (\oi_{\cnu})_!, (i_{\cnu})_!,(\oi_{\leq \cnu})^*, (\oi_{\cnu})^*, (i_{\cnu})^*, (\on{ev}_\cnu)^*$$

\begin{lem}
For $\CF\in \Pervsmk$ and any $\CS\in \sD(G([t]/t^k))$, the perverse cohomologies
of both $\CS\starstar\CF$ and $\CS\starshriek\CF$ belong to $\Pervsmk$.
\end{lem}

\begin{proof}

This follows immediately, since the action of $G([t]/t^k)$ extends to
$_\infty^k\CZ^\cmu$, respects the factorization isomorphisms, and commutes with
the action of the group-schemes involved in the definition of $\Pervsmk$.

\end{proof}

In what follows we will be interested in the case $k=1$.

\sssec{}

Let us denote by
$^{I}_{\infty}\BunNb$ (resp., $^{I^0}_{\infty}\BunNb$) the quotient stack
of $^1_{\infty}\BunNb$ by $B\subset G$ (resp., $N\subset G$).

We will denote by
$$\Pervsmi\subset \Perv({}^{I}_{\infty}\BunNb) \text{ and }
\Pervsmni \subset \Perv({}^{I^0}_{\infty}\BunNb)$$
the full subcategories of, consisting of objects, whose pull-back
to $^1_{\infty}\BunNb$ belongs to $\Pervsmk$, $k=1$.

For $\cnu\in \cLambda$, let us denote by $_\cnu^I\BunNb$ (resp.,
$_\cnu^{I^0}\BunNb$) the corresponding locally closed substack of
$^{I}_{\infty}\BunNb$ (resp., $^{I^0}_{\infty}\BunNb$), and by
$\ol{\on{ev}}_\cnu$ the map from it to $B\backslash G/N^-$ (resp.,
$N\backslash G/N^-$).

For an element $\tw\in W_{aff}$, written as $w\cdot \cnu$ with $w\in W$,
we will denote by $_\tw^I\BunNb$ (resp., $_\tw^{I^0}\BunNb$) the
preimage under $\ol{\on{ev}}_\cnu$ of the Schubert cell
$$B\backslash (B\cdot w\cdot N^-)/N^-\subset B\backslash G/N^-.$$
Let $_\tw^I\BunN$ (resp., $_\tw^{I^0}\BunN$) be the preimage of
the same Schubert cell under the map $\on{ev}_\cnu:{}_\cnu^I\Bun_N\to B\backslash G/N^-$.
We will denote by $\oi_\tw$ and $i_\tw$ the corresponding locally closed embeddings.

\medskip

We will denote by $\fIC_\tw\in \Perv({}^{I}_{\infty}\BunNb)$
the intersection cohomology sheaf on $_\tw^I\BunNb$. In other words,
$$\fIC_\tw\simeq (\oi_{\cnu})_{!*}\Bigl(\IC_{w,G/B^-}\tboxtimes
\IC_{{}_\cnu^I\CN}\Bigr),$$
in the notation of \propref{sheaves on strata}. In particular, we see that
$\fIC_\tw$ is an object of $\Pervsmi$, and these sheaves are all the irreducibles
of the categories $\Pervsmi$ and $\Pervsmni$.

\medskip

For $\tw=w\cdot \cmu$ as above, let us denote by $\nabla_\tw$ and
$\Delta_\tw$ the complexes
$$(\oi_{\cnu})_{!}\Bigl(j_{!,w\cdot w_0}\tboxtimes
\IC_{{}_\cnu^I\CN}\Bigr) \text{ and }
(\oi_{\cnu})_{*}\Bigl(j_{*,w\cdot w_0}\tboxtimes
\IC_{{}_\cnu^I\CN}\Bigr),$$
respectively, where $j_{!,w\cdot w_0}$ (resp., $j_{*,w\cdot w_0}$)
is the perverse sheaf on $G/B^-$ corresponding to the same-named
perverse sheaf under the isomorphism $G/B^-\to G/B$, given by
the right multiplication by $w_0$.

According to the above, we can act by objects of $\sD(G/B)^B$ (resp., $\sD(G/B)^N$) on
objects of $\Pervsmi$ and obtain complexes, whose cohomologies belong
to $\Pervsmi$ (resp., $\Pervsmni$). Evidently, we have:

\begin{equation} \label{standards from one another}
j_{!,w_1}\star \nabla_{w_2}\simeq \nabla_{w_1\cdot w_2} \text{ and }
j_{*,w_1}\star \Delta_{w_2}\simeq \Delta_{w_1\cdot w_2},
\end{equation}
provided that $l(w_1\cdot w_2)=l(w_1)+l(w_2)$.

\begin{prop}   \label{affineness of cells}
Both $\nabla_\tw$ and $\Delta_\tw$ are perverse sheaves.
\end{prop}

>From \propref{sheaves on strata} we obtain:

\begin{cor}
Both $\nabla_\tw$ and $\Delta_\tw$
are objects of $\Pervsmi$,
\end{cor}

\begin{proof}

Evidently, we have:
$$\nabla_\tw\simeq (\oi_{\tw})_!(\IC_{{}_\tw^I\BunNb}).$$
We claim that the morphism $\oi_{\tw}$ is affine. Clearly, this would imply the
proposition. To simplify the notation we will assume that $\cnu=0$; the proof
in the general case is the same.

\medskip

For an element $w\in W$ we can find a weight $\lambda$ and $B$-stable subspaces
$$'V_w^\lambda\subset V_w^\lambda\subset V^\lambda,$$
with $\dim(V_w^\lambda/{}'V_w^\lambda)=1$,
such that a point of $G/B^-$, thought of as a quotient line
$\ell^\lambda\twoheadleftarrow V^\lambda$, belongs
to $B\cdot w \cdot B^- /B^-$ if and only if
the composition
$$'V_w^\lambda\to V^\lambda\to \ell^\lambda$$ is zero, and
$V_w^\lambda\to V^\lambda\to \ell^\lambda$
is non-zero.

\medskip

Then, $_{\tw}^I\BunNb$, as a substack of $^{I}_{\leq 0}\BunNb$, corresponds to those
$\kappa^\lambda$, for which the map
\begin{equation} \label{Demazure1}
('V_w^\lambda)_{\CP_{G,x}}\to (V^\lambda_{\CP_G})_x\to \CO_x\simeq \BC
\end{equation}
is zero, and
\begin{equation} \label{Demazure2}
(V_w^\lambda)_{\CP_{G,x}}\to (V^\lambda_{\CP_G})_x\to \CO_x\simeq \BC
\end{equation}
is non-zero. Note that $(V_w^\lambda)_{\CP_{G,x}}$ and $('V_w^\lambda)_{\CP_{G,x}}$
make sense as subspaces of $(V^\lambda_{\CP_G})_x$,
since a part of the data of a point of $^{I}_{\infty}\BunNb$ is the reduction of the
fiber $\CP_{G,x}$ of $\CP_G$ at $x$ to $B$.

Hence, the closure of $_{\tw}^I\BunNb$ is contained in the closed sub-stack
of $^{I}_{\leq 0}\BunNb$, consisting
of all those points, for which the composition in \eqref{Demazure1} vanishes.
The locus of non-vanishing of \eqref{Demazure2} is the complement
to a Cartier divisor in this closed substack.

\end{proof}

\sssec{}

We will work with the abelian category $\Pervsmni$ and its derived category,
denoted $\Dsmni$.

By \propref{stable under extensions}, for $\CF_1,\CF_2\in \Pervsmni$,
$$Ext^1_{\Pervsmni}(\CF_1,\CF_2)\simeq
Ext^1_{{}_\infty^{I^0}\BunNb}(\CF_1,\CF_2)_T,$$
where the subscript $T$ stands for the $T$-equivariant category. Hence, the map
$$Ext^2_{\Pervsmni}(\CF_1,\CF_2)\to
Ext^2_{{}_\infty^{I^0}\BunNb}(\CF_1,\CF_2)_T$$
in injective.

\medskip

>From \corref{Ext 2 injects}, and using the fact that each $_\tw^{I^0}\Bun_N$
is contractible, we obtain:

\begin{cor}  \label{no Ext 2}
$Ext^i_{\Pervsmni}(\nabla_\tw,\Delta_{\tw'})=0$ for $i=1,2$ and any $\tw,\tw'\in W_{aff}$.
\end{cor}

\medskip

\noindent{\it Remark.}
>From \corref{no Ext 2} one can formally deduce that
$Ext^i_{\Pervsmni}(\nabla_\tw,\Delta_{\tw'})$ vanishes for all $i>0$
and any $\tw,\tw'\in W_{aff}$. More generally, for $\CF\in \Pervsmni$,
$$Ext^i_{\Pervsmni}(\nabla_\tw,\CF)\simeq H^i({}_\tw^{I^0}\Bun_N,\CF|_{{}_\tw^{I^0}\Bun_N}).$$

Note that by \propref{sheaves on strata}, the $!$-restriction of any
$\CF\in\Pervsmni$ to $_\tw^{I^0}\Bun_N$ is a complex with constant cohomologies.
Since $_\tw^{I^0}\Bun_N$ is contractible,
$H^\bullet({}_\tw^{I^0}\Bun_N,\uBC)\simeq \BC$,  so,
the above expression for $Ext^i$ amounts to taking
stalks of $\CF$ on the stratum $_\tw^{I^0}\Bun_N$.

\sssec{The baby Whittaker case}

Let $\Perv({}_\infty^1\BunNb)^{N^-,\psi}$ be the category of $(N^-,\psi)$-equivariant
perverse sheaves on $_\infty^1\BunNb$. We introduce the category
$$\Pervsmw\subset \Perv({}_\infty^1\BunNb)_{N^-,\psi},$$ as the full
subcategory, consisting of objects,
which belong to $\Pervsmk$, $k=1$, when regarded merely as objects of
$\Perv({}_\infty^1\BunNb)$. This category is stable under extensions by
\propref{stable under extensions}.

By \propref{sheaves on strata}, we can produce objects in $\Pervsmw$,
starting from objects of $\Perv(G/B^-)^{N^-,\psi}$. We will denote by
$\psi_{G/B^-}$ the unique irreducible in the latter category, which
corresponds to $\psi_{G/B}$ under
$$\Perv(G/B^-)^{N^-,\psi}\simeq \Pervfdw.$$

For $\cnu\in \cLambda$, set
$$\fIC^\psi_\cnu:=(\oi_{\cnu})_{!*}(\psi_{G/B^-}\tboxtimes \IC_{{}_\cnu^1\CN}),$$
and
$$\nabla_\cnu^{\psi}:=(\oi_{\cnu})_!(\psi_{G/B^-}\tboxtimes \IC_{{}_\cnu^1\CN}),\,\,
\Delta_\cnu^{\psi}:=(\oi_{\cnu})_*(\psi_{G/B^-}\tboxtimes \IC_{{}_\cnu^1\CN}).$$

Since the embedding of the corresponding locally closed subset into
$_\cnu^1\BunNb$ is affine (cf. the proof of \propref{affineness of cells}),
both $\nabla_\cnu^{\psi}$ and $\Delta_\cnu^{\psi}$ are perverse sheaves, and hence,
by \propref{sheaves on strata}, are objects of $\Pervsmw$.
In \secref{proof of Whit theorem} we will prove the following:

\begin{thm} \label{Whittaker on semiinf}
The canonical maps $\nabla_\cnu^{\psi}\to \fIC^\psi_\cnu\to\Delta_\cnu^{\psi}$ is
are isomorphisms.
 \end{thm}

Thus, the extension of $\psi_{G/B}\tboxtimes \IC_{{}_\cnu\BunNb}$ under
$\oi_{\cnu}$ is clean, and $\nabla_\cnu^{\psi}\simeq\Delta_\cnu^{\psi}$ is irreducible.
Hence, the category $\Pervsmw$ is semi-simple and equivalent to $\cT\mod$.

\sssec{}

Let us denote by $$\on{Av}_{!,N^-,\psi},\on{Av}_{!,N^-,\psi}:
\sD({}_{\infty}^1\BunNb)\to \sD({}_{\infty}^1\BunNb)^{I^-,\psi}$$
the functors, which are left and right adjoint, respectively, to
$\sD({}_{\infty}^1\BunNb)^{I^-,\psi}\to \sD({}_{\infty}^k\BunNb)$.
As in \propref{psi averaging} we obtain:

\begin{lem}
There exists an isomorphism of functors
$$\on{Av}_{!,N^-,\psi}[-\dim(\fn)]|_{\sD({}_{\infty}^{I^0}\BunNb)}\to
\on{Av}_{*,N^-,\psi}[\dim(\fn)]|_{ \sD({}_{\infty}^{I^0}\BunNb)}.$$
Moreover, the resulting functor
$\on{Av}_{N^-,\psi}:\sD({}_{\infty}^{I^0}\BunNb)\to \sD({}_{\infty}^1\BunNb)^{I^-,\psi}$
is exact.
\end{lem}

\medskip

Let us call an object of $\Pervsmni$ partially integrable if all of its
irreducible sub-quotients are of the form $\fIC_{w\cdot \cnu}$, $w\neq w_0$.
Thus, the only irreducibles, that are not partially integrable are
$\fIC_{w_0\cdot \cnu}$. Let us denote by $\Pervsmnif$ the resulting
quotient abelian category.

The following is parallel to \propref{whit and part int}.

\begin{prop}  \label{non-degen quotient on semi-inf} \hfill

\smallskip

\noindent{\em (1)}
The functor $$\on{Av}_{I^-,\psi}:\Pervsmni\to \Pervsmw$$
factors through $\Pervsmnif$.

\smallskip

\noindent{\em (2)}
The resulting functor $$\Pervsmnif\to \Pervsmw$$
is faithful.

\end{prop}

\begin{proof}

To prove the first statement we have to show that $\on{Av}_{I^-,\psi}(\fIC_{w\cdot \cnu})=0$
for $w\neq w_0$. This is nearly evident: such an irreducible is a pull-back
from the quotient stack $P_\imath\backslash{}_{\infty}^1\BunNb$, where $P_\imath$ is
some sub-minimal parabolic in $G$. Our assertion follows from the fact that
the direct image of $\psi_{G/B}$ under $G/B\to G/P_\imath$ vanishes.

To prove the second statement, it suffices to show that
$$\on{Av}_{I^-,\psi}(\fIC_{w_0\cdot \cnu})\simeq \fIC^\psi_\cnu.$$
We know that the left-hand side is a perverse sheaf, and the isomorphism
over the open part of the support, namely $_{\cnu}^1\BunNb$, is evident.
The fact that the left-hand side is a Goresky-MacPherson extension from
this sub-stack follows from the exactness of the functor $\on{Av}_{I^-,\psi}$,
and the fact that it commutes with all $\oi^*_{\cnu'}$ and $\oi^!_{\cnu'}$.

\end{proof}

\begin{cor}  \label{cosocle on semiinf} \hfill

\smallskip

\noindent{\em (1)}
The kernel of $\nabla_{w_0\cdot \cnu}\to \fIC_{w_0\cdot \cnu}$ is
partially integrable.

\smallskip

\noindent{\em (2)}
$\fIC_{w_0\cdot \cnu}$ is the cosocle of $\Delta_{\cnu}$ and socle
of $\nabla_\cnu$.

\smallskip

\noindent{\em (3)}
For any $w\in W$, $\fIC_{w_0\cdot \cnu}$ is the only non-partially
integrable constituent of $\nabla_{w\cdot \cnu}$.

\end{cor}

\begin{proof}

Evidently, we have
$$\on{Av}_{I^-,\psi}(\nabla_{w_0\cdot \cnu})\simeq \nabla^\psi_\cnu.$$
Combining this with \propref{non-degen quotient on semi-inf} and
\thmref{Whittaker on semiinf}, we arrive to the assertion of point (1).
Point (3) follows from point (1) by \eqref{standards from one another}.
Finally, point (2) follows from point (1) in the same way as in the
proof of \propref{cosocle of 0 costandard}.

\end{proof}

We will now introduce one more object of $\Pervsmni$. For $\cnu\in \cLambda$
set
$$\bPi_{!,\cnu}:=(\oi_{\cnu})_{!}\Bigl(\Xi\tboxtimes \IC_{{}_\cnu^I\BunNb}\Bigr) \text{ and }
\bPi_{*,\cnu}:=(\oi_{\cnu})_{*}\Bigl(\Xi\tboxtimes \IC_{{}_\cnu^I\BunNb}\Bigr),$$
where $\Xi$ is the perverse sheaf on $N\backslash G/B^-$, corresponding to the
same-named perverse sheaf on $G/B$.

\begin{thm} \label{cleanness of CXi}
The canonical map $\bPi_{!,\cnu}\to \bPi_{*,\cnu}$ is an
isomorphism.
\end{thm}

\begin{proof}

Consider the convolution with $\Xi$ as a functor
$\Pervsmi\to \Pervsmni$. As usual, this functor annihilates
all partially integrable objects.

Evidently,
$$\bPi_{!,\cnu}\simeq \Xi\star \nabla_{w_0,\cnu} \text{ and }
\bPi_{!*\cnu}\simeq \Xi\star \Delta_{w_0,\cnu}.$$

Our assertion follows now from \corref{cosocle on semiinf}, which
implies that the
cone of the map $\nabla_{w_0,\cnu} \to \Delta_{w_0,\cnu}$
is partially integrable.

\end{proof}

\sssec{}  \label{cosocle of other costandard, semiinf}

We will now establish the following fact, parallel to \propref{geometric cosocles}(2):

\begin{prop}
For $\cnu\in \cLambda$ there exists a non-zero map
$$\Delta_{w_0\cdot \cnu}\to \fIC_{\cnu-2\crho},$$
where $\fIC_{\cnu-2\crho}\in \Pervsmg$ is thought of as
an object of $\Pervsmi$.
\end{prop}

\begin{proof}

As in the proof of \propref{cosocle of another costandard}, we have the functor
$$\on{Av}_{!,G/B}:\sD({}_\infty^I\BunNb)\to \sD({}_\infty\BunNb),$$ left adjoint
to the forgetful functor. By definition,
\begin{equation}  \label{desired Hom, semiinf}
Hom_{\Pervsmi}(\Delta_{w_0\cdot \cnu},\fIC_{\cnu-2\crho})\simeq
Hom_{\sD({}_\infty\BunNb)}(\on{Av}_{!,G/B}(\Delta_{w_0\cdot \cnu}),\fIC_{\cnu-2\crho}).
\end{equation}

However, since $G/B$ is proper,
$$\on{Av}_{!,G/B}(\Delta_{w_0\cdot \cnu})\simeq (\oi_{\cnu})_*(\IC_{_\cnu\BunNb})[\dim(\fn)].$$

Hence, the assertion of the proposition follows from \corref{top naive extension}.

\end{proof}

\section{Convolution}

\ssec{Definition of convolution}

\sssec{}

Consider the Hecke stack for $G$ at $x$:
$$\Bun_G \overset{\hl_G}\longleftarrow \CH_{G,x}\overset{\hr_G}\longrightarrow \Bun_G,$$
and for two integers $k_1,k_2$ let $^{k_1,k_2}\CH_{G,x}$ denote its base change with
respect to
$$^{k_1}\Bun_G\times {}^{k_2}\Bun_G\to \Bun_G\times \Bun_G.$$
By a slight abuse of notation we will continue to denote by $\hl_G,\hr_G$ the
projections of $^{k_1,k_2}\CH_{G,x}$ on $^{k_1}\Bun_G$ and $^{k_2}\Bun_G$,
respectively.

We can regard $^{k_1,k_2}\CH_{G,x}$ over $^{k_2}\Bun_G$ as the space
associated with the canonical $G^{k_2}$-torsor $\CG^{k_2}_x$ over
$^{k_2}\Bun_G$ and the $G^{k_2}$-space $G((t))/G^{k_1}$:
$$^{k_1,k_2}\CH_{G,x}\simeq G((t))/G^{k_1}\overset{G^{k_2}}\times \CG^{k_2}_x.$$
We also have a symmetric picture:
$$^{k_1,k_2}\CH_{G,x}\simeq G((t))/G^{k_2}\overset{G^{k_1}}\times \CG^{k_1}_x.$$

Recall now that there exists a canonical equivalence of derived categories
$$\CS\mapsto \CS^{op}:
\sD_{G^{k_1}}(G((t))/G^{k_2})\simeq \sD_{G^{k_2}}(G((t))/G^{k_1}).$$
It is defined as follows.

First of all, it is clear that $G^{k_1}$-invariant sub-schemes of $G((t))/G^{k_2}$
are in bijection with $G^{k_2}$-invariant sub-schemes in $G^{k_1}\backslash G((t))$.
For $\CS\in \sD_{G^{k_1}}(G((t))/G^{k_2})$, let $\CY$ be the corresponding
finite-dimensional sub-scheme of $G^{k_1}\backslash G((t))$. There exists
an integer $k'_1>>0$, such that if we denote by $\CY'$ the preimage of
$\CY$ in $G^{k'_1}\backslash G((t))$, the map $\CY'\to G((t))/G^{k_2}$ is well-defined.
The pull-back $\CS'$ of $\CS$ to $\CY'$ is an $G^{k_1}/G^{k'_1}$-equivariant,
and, hence, descends to a well-defined $G^{k_2}$-equivariant object of
$\sD(G^{k_1}\backslash G((t)))$.

Finally, the desired functor is obtained by applying the inversion on $G((t))$.

\sssec{}

As in \cite{BG} we have a commutative diagram, in which both squares are Cartesian
$$
\CD
_{\infty}^{k_1}\BunNb @<{\hl'_G}<< {}^{k_1,k_2}\CH_{G,N^-,x} @>{\hr'_G}>>
_{\infty}^{k_2}\BunNb \\
@V{\fp}VV   @VVV  @V{\fp}VV \\
^{k_1}\Bun_G @<{\hl_G}<< {}^{k_1,k_2}\CH_{G,x} @>{\hr_G}>>
^{k_2}\Bun_G.
\endCD
$$

For a complex $\CF$ on $_{\infty}^{k_2}\BunNb$ and a $G^{k_1}$-equivariant
complex $\CS$ on $G((t))/G^{k_2}$ let
$\CS^{op}\tboxtimes \CF$ be the corresponding complex on $^{k_1,k_2}\CH_{G,N^-,x}$.
We set
$$\CS\starstar\CF:=(\hl'_G)_*(\CS^{op}\tboxtimes \CF),\,
\CS\starshriek\CF:=(\hl'_G)_!(\CS^{op}\tboxtimes \CF)\in {\mathsf D}({}_{\infty}^{k_1}\BunNb).$$

Evidently, when $k_1=k_2=k$, and $\CS$ is supported on $G[[t]]/G^k\subset G((t))/G^k$,
we arrive to the functors discussed in \secref{naive convolution}.

\medskip

The following is straightforward from the definitions:

\begin{lem}
For $\CS\in \sD_{G^{k_1}}(G((t))/G^{k_2})$ the functor
$$\CF\mapsto \CS\starshriek \CF: \sD({}_\infty^{k_2}\BunNb)\to
\sD({}_\infty^{k_1}\BunNb)$$ is the left adjoint of
$$\CF'\mapsto \BD(\CS^{op})\starstar \CF': \sD({}_\infty^{k_1}\BunNb)\to
\sD({}_\infty^{k_2}\BunNb).$$
\end{lem}

The above picture admits the following variants. First, we can replace the equivariance condition
on $\CS\in \sD\Bigl(G((t))/G^{k_2}\Bigr)$ with respect to $G^{k_1}$ by that of $I^0$, $I$ or
$(I^-,\psi)$. In the case the target will be the corresponding category
$\sD({}_{\infty}^{I^0}\BunNb)$, $\sD({}_{\infty}^I\BunNb)$ or
$\sD({}_{\infty}^1\BunNb)^{N^-,\psi}$.

\medskip

Secondly, instead $\sD\Bigl(G((t))/G^{k_2}\Bigr)$ we can consider $\sD(\Fl)$
or $\sD(\Gr)$. We obtain the convolution functors
$$\sD(\Fl)^{G^k}\times \sD({}_{\infty}^I\BunNb)\to \sD({}{}_{\infty}^k\BunNb) \text{ and }
\sD(\Gr)^{G^k}\times \sD({}_{\infty}\BunNb)\to \sD({}{}_{\infty}^k\BunNb).$$

In both these cases, the *-convolution coincides with the !-convolution, since
$\Fl$ and $\Gr$ are ind-proper. We will denote the resulting functor simply by $\star$.
Here again, the equivariance condition with respect
to $G^k$ can be replaced by any of $I^0$-, $I$- or $(I^-,\psi)$-equivariance conditions.

\sssec{}

We will now show that the convolution functors essentially preserve our category
$\Pervsm$.

\begin{prop}  \label{convolution respects}
If $\CF\in \Pervsmktwo$, then the perverse cohomologies of both
$\CS\starstar\CF$ and $\CS\starshriek\CF$ belong to $\Pervsmkone$.
\end{prop}

The rest of this subsection is devoted to the proof of this proposition.
First, let us notice that if $\CF$ satisfies condition (1), then so do the complexes
$(\hl'_G)_*(\CS^{op}\tboxtimes \CF)$ and $(\hl'_G)_!(\CS^{op}\tboxtimes \CF)$. Hence, by
\lemref{stronger equivariance} and \propref{equiv of (1) and (2)}, these complexes
satisfy the weak factorization property. Hence, to show that their perverse cohomologies
satisfy the full factorization property, it is enough to show that their pull-backs to
$$\Bigl(\CZ^{\cmu_2-\cmu_1}\times {}_{\infty}^{k_1}\CZ^{\cmu_1}\Bigr)
\underset{\oX^{\cmu_2-\cmu_1}\times {}_{\infty}X^{\cmu_1}}\times \Bigl(\oX^{\cmu_2-\cmu_1}
\times {}_{\infty}X^{\cmu_1}\Bigr)_{disj}$$
can be written as extensions of complexes, each of which has the form
$\IC_{\CZ^{\cmu_2-\cmu_1}}\boxtimes \CF'$, where $\CF'$ is some complex on
$_{\infty}^{k_1}\CZ^{\cmu_1}$.

\medskip

Let us denote by $\CY^{\cmu}$ the Cartesian product
$$_{\infty}^{k_1}\CZ^{\cmu}\underset{_{\infty}^{k_1}\BunNb}\times {}^{k_1,k_2}\CH_{G,N,x}.$$

As in \lemref{factorization}, we have a canonical isomorphism
$$\CY^{\cmu_2}\underset{_{\infty}X^{\cmu_2}}\times
\Bigl(\oX^{\cmu_2-\cmu_1}\times {}_{\infty}X^{\cmu_1}\Bigr)_{disj}\simeq
\Bigl(\CZ^{\cmu_2-\cmu_1}\times \CY^{\cmu_1}\Bigr)
\underset{\oX^{\cmu_2-\cmu_1}\times {}_{\infty}X^{\cmu_1}}
\times \Bigl(\oX^{\cmu_2-\cmu_1}\times {}_{\infty}X^{\cmu_1}\Bigr)_{disj}.$$

We claim that the pull-back under
$$\CY^{\cmu_2}\underset{_{\infty}X^{\cmu_2}}\times
\Bigl(\oX^{\cmu_2-\cmu_1}\times {}_{\infty}X^{\cmu_1}\Bigr)_{disj}\to
\CY^{\cmu_2}\to {}^{k_1,k_2}\CH_{G,N,x}$$
of $\CS^{op}\tboxtimes \CF$ is an extension of complexes, each of which has the form
$\IC_{\CZ^{\cmu_2-\cmu_1}}\boxtimes \CF''$, where $\CF''$ is some complex
on $\CY^{\cmu_1}$. This would clearly imply our assertion.

\medskip

Note that $\CY^\cmu$ can be represented as a union of locally closed sub-stacks
$_{\cnu}\CY^{\cmu}$ for $\cnu\in \cLambda$, where a point
$(\CP_G,\{\kappa^\lambda\},\{\kappa^{\lambda,-}\})$ belongs to $_{\cnu}\CY^{\cmu}$
if and only if each $\kappa^{\lambda,-}$ has a pole of order $\langle \lambda,\cnu\rangle$
at $x$.

Note that we have a natural map $_{\cnu}\CY^{\cmu}\to {}_{\infty}^{k_2}\CZ^{\cmu+\cnu}$,
that covers the map $$\hr'_G:{}^{k_1,k_2}\CH_{G,N,x}\to {}_{\infty}^{k_2}\BunNb.$$ Moreover,
the diagram
$$
\CD
_{\cnu}\CY^{\cmu_2}\underset{_{\infty}X^{\cmu_2}}\times
\Bigl(\oX^{\cmu_2-\cmu_1}\times {}_{\infty}X^{\cmu_1}\Bigr)_{disj} @>>>
\CZ^{\cmu_2-\cmu_1}\times {}_{\cnu}\CY^{\cmu_1}  \\
@VVV   @VVV \\
_{\infty}^{k_2}\CZ^{\cmu_2+\cnu}\underset{_{\infty}X^{\cmu_2+\cnu}}\times
\Bigl(\oX^{\cmu_2-\cmu_1}\times {}_{\infty}X^{\cmu_1+\cnu}\Bigr)_{disj}
@>>>
\CZ^{\cmu_2-\cmu_1}\times {}_{\infty}^{k_2}\CZ^{\cmu_1+\cnu}
\endCD
$$
is commutative. Hence, our assertion follows from condition (2) imposed on $\CF$.

\medskip

It remains to show that the perverse cohomologies of
$\CS\starstar\CF$ and $\CS\starshriek\CF$ satisfy condition (3). Since we have
to check an equivariance condition with respect to a unipotent group-scheme,
it is enough to show that their pull-backs to $_\infty^{k_1}\CZ^\cmu$ can be written
as extensions of complxes satisfying this equivariance condition. This follows
in the same way as above, by sub-dividing the stack $\CY^\cmu$ into the locally
closed substacks $_\cnu\CY^\cmu$.

\ssec{Exactness and smallness}

\sssec{}

Consider the convolution functor
$$\sD(\Gr)^{G^k}\times \Pervsmg\to \sD({}_\infty^k\BunNb).$$

Since $\Pervsmg$ is semi-simple, it is enough to evaluate the
above functor on the objects of the form $\fIC_\cnu$, $\cnu\in \cLambda$.

\begin{thm}  \label{exactness of convolution}
The functor
$$\CS\mapsto \CS\star \fIC_\cnu:\sD(\Gr)^{G^k}\to \sD({}_\infty^k\BunNb).$$
is exact.
\end{thm}

\begin{proof}

Since the situation is self-dual with respect to the Verdier duality,
it is sufficient to show that for $\CS\in \Perv(\Gr)^{G^k}$,
the convolution $\CS\star \fIC_\cnu$ is supported in non-positive
cohomological degrees. For that
it is sufficient to show that $\oi_\cmu^*(\CS\star \fIC_\cnu)$ is supported
in non-positive cohomological degrees for every $\cmu\in \cLambda$.

Consider the preimage  $(\hl'_G)^{-1}\Bigl({}_{\cmu}^k\BunNb\Bigr)\subset
{}^{k,0}\CH_{G,N,x}$. It admits a decomposition into locally closed pieces
\begin{equation} \label{finer strata}
(\hl'_G)^{-1}\Bigl({}_{\cmu}^k\BunNb\Bigr)\cap
(\hr'_G)^{-1}\Bigl({}_{\cmu'}\BunNb\Bigr)\cap {}^{k,0}\CH^\clambda_{G,N,x}
\end{equation}
for $\cmu'\in \cLambda$ and $\clambda\in \cLambda^+$, where
$^{k,0}\CH^\clambda_{G,N,x}$ is the preimage of the corresponding
locally closed sub-stack in $\CH_{G,x}$.

The statement of the theorem would follow once we prove the following:

\smallskip

\noindent(1) The dimension of fibers of the map
$$\hl'_G:(\hl'_G)^{-1}\Bigl({}_{\cmu}^k\BunNb\Bigr)\cap
(\hr'_G)^{-1}\Bigl({}_{\cmu'}\BunNb\Bigr)\cap {}^{k,0}\CH^\clambda_{G,N,x}\to
{}_\cmu^k\BunNb$$ is
$\leq \langle \cmu'-\cmu+\clambda,\rho\rangle$.

\smallskip

\noindent(2) The *-restriction of $\CS^{op}\tboxtimes \fIC_\cnu$ to
$(\hl'_G)^{-1}\Bigl({}_{\cmu}^k\BunNb\Bigr)\cap
(\hr'_G)^{-1}\Bigl({}_{\cmu'}\BunNb\Bigr)\cap {}^{k,0}\CH^\clambda_{G,N,x}$
lives in the cohomological degrees $\leq -\langle \cmu'-\cmu+\clambda,\rho\rangle$.

\medskip

The first assertion follows from the identification of the locally closed substack
from \eqref{finer strata}, projecting to $_\cmu^k\BunNb$ by means of $\hl'_G$,
with
\begin{equation} \label{fiber bundle left}
\Bigl(\Gr^{\clambda}\cap N^-((t))\cdot (\cmu-\cmu')\Bigr)\overset{N^-[[t]]}\times
{}_\cmu^k\CN,
\end{equation}
where $^k_\cmu\CN$ is the $N^-[[t]]$-torsor over $_\cmu^k\BunNb$
introduced before.

\medskip

To prove the second assertion let us view the locally closed sub-stack of
\eqref{finer strata} projecting to $_{\cmu'}\BunNb$ by means of $\hr'_G$;
it identifies with
$$p_k^{-1}\Bigl(\Gr^{-w_0(\clambda)} \cap N^-((t))\cdot (\cmu'-\cmu)\Bigr)\overset{N^-[[t]]}\times
_{\cmu'}\CN,$$
where $p_k$ is the projection $G((t))/G^k\to \Gr$.

The *-restriction of $\CS^{op}\tboxtimes \fIC_\cnu$ to it identifies with
$$\CS^{op}|_{p_k^{-1}\Bigl(\Gr^{-w_0(\clambda)} \cap
N^-((t))\cdot (\cmu'-\cmu)\Bigr)}\tboxtimes \fIC_\cnu|_{_{\cmu'}\BunNb}.$$

Hence, it is enough to show that the *-restriction of $\CS^{op}$ to
$$p_k^{-1}\Bigl(\Gr^{-w_0(\clambda)} \cap N^-((t))\cdot (\cmu'-\cmu)\Bigr)$$ lives
in the cohomological degrees $\leq -\langle \cmu'-\cmu+\clambda,\rho\rangle$.

\medskip

First, the restriction to $p_k^{-1}(\Gr^{-w_0(\clambda)})$ lives in non-positive degrees,
since $\CS$ was assumed perverse. By assumption, this complex is
$G[[t]]$-equivariant, and hence, universally locally acyclic over $\Gr^{-w_0(\clambda)}$,
since the latter is a $G[[t]]$-homogeneous space. Since
$$\on{codim}\Bigl(\Gr^{-w_0(\clambda)} \cap N^-((t))\cdot (\cmu'-\cmu),
\Gr^{-w_0(\clambda)}\Bigr)\geq \langle \cmu'-\cmu+\clambda,\rho\rangle,$$
our assertion follows.

\end{proof}

\sssec{Convolution in the spherical case}   \label{spherical convolution}

We will now study a particular case of the above situation, when the functor
we consider is:
$$\Sph\times \Pervsmg\to \Pervsmg.$$

\begin{prop}  \label{sph conv}
For $V\in \Rep(\cG)$ and $\cnu\in \cLambda$, there exists a canonical isomorphism
$$\CV\star \fIC_\cnu\simeq \underset{\cmu}\oplus\, \fIC_{\cnu+\cmu}\otimes \uV(\cmu).$$
Moreover, for $V,U\in \Rep(\cG)$, the diagram
$$
\CD
(\CU\star \CV)\star \fIC_\cnu @>{\sim}>> \underset{\cmu'}\oplus\,
(\CU\star \fIC_{\cnu+\cmu'})\otimes \uV(\cmu') \\
@V{\sim}VV    @V{\sim}VV  \\
\underset{\cmu}\oplus\, \fIC_{\cnu+\cmu}\otimes \underline{(U\otimes V)}(\cmu) @>{\sim}>>
\underset{\cmu',\cmu''}\oplus\, \fIC_{\cnu+\cmu'+\cmu''} \otimes \uU(\cmu'')\otimes \uV(\cmu')
\endCD
$$
commutes.
\end{prop}

Before giving the proof let us recall that the for $V\in \Rep(\cG)$ and the corresponding
object $\CV\in \Sph$, we have a canonical isomorphism
\begin{equation} \label{geometric weight spaces}
V(\cmu)\simeq H_c^{-\langle 2\rho,\cmu\rangle}\Bigl(N^-((t))\cdot \cmu,
\CV|_{N^-((t))\cdot \cmu}\Bigr).
\end{equation}

\begin{proof}

Note first that the result of the convolution $\CV\star \fIC_\cnu$ is an object of $\Pervsmg$,
and hence, is semi-simple. (Alternatively, semi-simplicity follows from the
decomposition theorem, since every $\CV\in \Sph$ is a direct sum of
intersection cohomology sheaves.)

By the proof of \thmref{exactness of convolution},
$$Hom(\fIC_{\cnu+\cmu},\CV\star \fIC_\cnu)\simeq
H_c^{-\langle 2\rho,\cmu\rangle}\Bigl(N^-((t))\cdot \cmu, \CV|_{N^-((t))\cdot \cmu}\Bigr),$$
which is exactly the expression that appears in \eqref{geometric weight spaces}.

The second assertion of the proposition follows from the definition of the structure
of the tensor functor on $V\mapsto \CV:\Rep(\cG)\to \Sph$, cf. \cite{MV} or \cite{BG1}.

\end{proof}

The commutativity of the following two diagrams also follows from \eqref{geometric weight spaces}:
\begin{equation} \label{comm sph1}
\CD
\IC_{\clambda,\Gr}\star \fIC_\cnu @>>> \IC_{\clambda,\Gr}\star \IC_{\cmu,\Gr}\star
\IC_{-w_0(\cmu),\Gr}\star \fIC_\cnu \\
@VVV  @VVV \\
\fIC_{\cnu+\clambda} @<<< \IC_{\clambda,\Gr}\star \IC_{\cmu,\Gr} \star \fIC_{\cnu-\cmu},
\endCD
\end{equation}
where the left vertical arrow comes from taking the direct summand corresponding
to $\uV^\clambda(\clambda)$, and the right vertical arrow comes from taking the
summand corresponding to $\uV^{-w_0(\cmu)}(-\cmu)$.

For the following diagram $V$ is an object of $\Rep(\cG)$ and $\clambda$ is a coweight
large compared to $V$:
\begin{equation} \label{comm sph2}
\CD
\CV\star \fIC_\cnu @>>>
(\IC_{-w_0(\clambda),\Gr}\star \IC_{\clambda,\Gr}\star \CV)\star \fIC_\cnu  \\
@VVV  @V{\sim}VV \\
\underset{\cmu}\oplus\, \fIC_{\cnu+\cmu}\otimes \uV(\cmu) & &
\underset{\cmu}\oplus\, \IC_{-w_0(\clambda),\Gr}\star\IC_{\clambda+\cmu,\Gr}\star
\fIC_\cnu \otimes \uV(\cmu) \\
@V{\on{id}}VV  @VVV  \\
\underset{\cmu}\oplus\, \fIC_{\cnu+\cmu}\otimes \uV(\cmu)
@<<< \underset{\cmu}\oplus\, \IC_{-w_0(\clambda),\Gr}\star \fIC_{\clambda+\cnu+\cmu}\otimes \uV(\cmu).
\endCD
\end{equation}

\ssec{Convolution with $\Pervgri$}

\sssec{}

We will now consider the convolution functor
$$\Perv_{G[[t]]}(\Fl)\times \Pervsmg\to \Pervsmi.$$

Recall the objects $\CL^w\in \Pervgri$ defined for $w\in W$
We will prove:

\begin{thm}  \label{conv is simple}
If $\CL^w=\IC_{w\cdot \clambda,\Gr}$, then
$$\CL^w\star \fIC_\cnu\simeq \fIC_{w\cdot (\clambda+\cnu)}$$
\end{thm}

The rest of this sub-section is devoted to the proof of this theorem.
We will retrace the argument proving \thmref{exactness of convolution}
and show that the map defining $\CL^w\star \fIC_\cnu$ is small
(vs. semi-small).

\medskip

First, to calculate the top (=$0$-th) cohomology of $\CL^w\star \fIC_\cnu$
we only need to consider the locally closed sub-stack of $^{I,0}\CH_{G,N^-,x}$
isomorphic to
$$\BO_w \overset{N^-[[t]]}\times {_\cnu}^I\CN,$$
and the constant perverse sheaf on it, where $\BO_\tw$ is the open $G[[t]]$-orbit
in the support of $\left(\CL^w\right)^{op}$ on $\Fl$. Its intersection with the
preimage of $_{w'\cdot \cmu}^I\BunNb$ under $\hl'_G$ can be described as
follows.

Note that the pull-back $_{w'\cdot \cmu}^I\BunNb\underset{_\cmu^I\BunNb}\times
{}_{\cmu}\CN$ of the $N^-[[t]]$-torsor $_\cmu\CN$ to $_{w'\cdot \cmu}^I\BunNb$
admits a reduction to the subgroup $N^-[[t]]\cap \on{Ad}_{(w')^{-1}}(I)$.
Then the above intersection identifies with the total space of the bundle
associated with the $N^-[[t]]\cap \on{Ad}_{(w')^{-1}}(I)$-space
\begin{equation} \label{int int}
N^-((t))\cdot (\cmu-\cnu)\cap \Bigl(\on{Ad}_{(w')^{-1}}(I)\cdot
((w')^{-1}\cdot w(\clambda))\Bigr)\subset N^-((t))\cdot (\cmu-\cnu)\cap \Gr^\clambda.
\end{equation}

Evidently, when $\cmu=\clambda+\cnu$ and $w'=w$ the above intersection
is the point-scheme. This means that $\fIC_{w\cdot (\clambda+\cnu)}$
indeed appears as a direct summand in the convolution $\CL^w\star \fIC_\cnu$.
It remains to show that if $\cmu\neq \clambda+\cnu$ or $w'\neq w$,
then the scheme in \eqref{int int} is of dimension strictly less than
$\langle \cnu-\cmu+\clambda,\rho\rangle$.

\medskip

We will deduce this from \thmref{geometric steinberg}. Let us take $\cmu_1$ to be a
large dominant coweight and set $\cnu_1=\cmu_1+\cmu-\cnu$. We will show that
if the dimension of \eqref{int int} violated the above inequality, the perverse sheaf
$\IC_{w'\cdot \cnu_1,\Gr}$ would appear as a direct summand of
$\IC_{w\cdot \clambda,\Gr} \star \IC_{\cmu_1,\Gr}$. For that end, it is sufficient to
show that the fiber of
$$\Bigl(I\cdot (w\cdot \clambda)\Bigr)\star \Gr^{\cmu_1}$$ over the point
$w'\cdot \cnu_1$ is of dimension $\geq \langle \cnu-\cmu+\clambda,\rho\rangle$.
We claim that the above fiber contains a subscheme is isomorphic to the scheme
\eqref{int int}.

\medskip

Consider the orbit of the group $\on{Ad}_{w'}N^-((t))$ passing through
$w'\cdot \cnu_1\in \Gr$. Its preimage in $\Bigl(I\cdot (w\cdot \clambda)\Bigr)\star \Gr^{\cmu_1}$
is the union over parameters $\cnu'_1$ of the schemes

\smallskip

\begin{equation}  \label{preimage of orbit}
\Bigl(\left(\on{Ad}_{w'}N^-((t))\cdot (w'\cdot \cnu'_1)\right)\cap
\left(I\cdot (w\cdot \clambda)\right)\Bigr)\star
\Bigl(\left(\on{Ad}_{w'}N^-((t))\cdot (w'\cdot (\cnu_1-\cnu'_1))\right)
\cap \Gr^{\cmu_1} \Bigr),
\end{equation}
each of which is fibered over
\begin{equation} \label{another intersection of orbits}
\Bigl(\on{Ad}_{w'}N^-((t))\cdot (w'\cdot \cnu'_1)\Bigr)\cap
\Bigl(I\cdot (w\cdot \clambda)\Bigr)
\end{equation}
with a typical fiber
$$\Bigl(\on{Ad}_{w'}N^-((t))\cdot (w'\cdot (\cnu_1-\cnu'_1))\Bigr)
\cap \Gr^{\cmu_1}.$$

Let us take $\cnu'_1=\cnu_1-\cmu_1$. We claim that that the intersection of
\eqref{preimage of orbit} with the preimage of the point $w'\cdot \cnu_1$ in
$\Bigl(I\cdot (w\cdot \clambda)\Bigr)\star \Gr^{\cmu_1}$ surjects onto
the scheme in \eqref{another intersection of orbits}. This would imply
our assertion, since the schemes \eqref{int int} and
\eqref{another intersection of orbits} are isomorphic for the above choice
of $\cnu'_1$.

\medskip

This amounts to showing that the subscheme
$$(-\cnu_1\cdot (w')^{-1})\cdot \Bigl(\left(\on{Ad}_{w'}N^-((t))\cdot
(w'\cdot (\cmu-\cnu))\right)
\cap \left(I\cdot (w\cdot \clambda)\right)\Bigr)$$
is contained in $\Gr^{-w_0(\cmu_1)}$.

Let $N^{?}$ be the group-subscheme of $N^-((t))$, such that
$$\left(\on{Ad}_{w'}N^-((t))\cdot (w'\cdot (\cmu-\cnu))\right)
\cap \left(I\cdot (w\cdot \clambda)\right)$$ is contained in
$$\left(\on{Ad}_{w'}(N^{?})\cdot (w'\cdot (\cmu-\cnu))\right)
\cap \left(I\cdot (w\cdot \clambda)\right).$$

We have to show that
$$(-\cnu_1)\cdot N^?\cdot (\cmu-\cnu)\subset \Gr^{-w_0(\cmu_1)},$$
which is equivalent to
$$\on{Ad}_{-\cnu_1}N^? \cdot (-\cmu_1)\subset \Gr^{-w_0(\cmu_1)}.$$
However, the latter containment is valid, whenever $\cnu_1$ is dominant enough so
that $\on{Ad}_{-\cnu_1}(N^?)\subset N^-[[t]]$.

\medskip

\noindent{\it Remark.}
Let us note that the fiber of $\Bigl(I\cdot (w\cdot \clambda)\Bigr)\star \Gr^{\cmu_1}$
over $w'\cdot \cnu_1$ is in fact entirely contained in the subscheme
\eqref{preimage of orbit} with $\cnu'_1=\cnu_1-\cmu_1$,
and it maps to the scheme \eqref{another intersection of orbits}
isomorphically.

To prove the first assertion note that there are only finitely
many $\cnu'_1$'s, for which the base \eqref{another intersection of orbits}
is non-empty. For any $\cnu'_1$ other than $\cnu_1-\cmu_1$ the subscheme
$$(-\cnu_1\cdot (w')^{-1})\cdot \Bigl(\left(\on{Ad}_{w'}N^-((t))\cdot
(w'\cdot \cnu'_1)\right)\cap \left(I\cdot (w\cdot \clambda)\right)\Bigr)$$
will have an empty intersection with $\Gr^{-w_0(\cmu_1)}$, because
eventually
$$\Bigr(\on{Ad}_{-\cnu_1}(N^?)\cdot  (\cnu'_1-\cnu_1)\Bigl)\,\cap \Gr^{-w_0(\cmu'_1)}
=\emptyset.$$

The second assertion is evident, since every fiber of $\pi:\Gr\star \Gr\to \Gr$
embeds into the base $\Gr$.

\sssec{The baby Whittaker case}  \label{proof of Whit theorem}

Our present goal is to prove \thmref{Whittaker on semiinf}. By Verdier duality,
it is sufficient to show that the map
$$\nabla_\cnu^\psi\to \fIC_\cnu^\psi$$
is an isomorphism. Suppose it is not, and let us look at the quotient
perverse sheaf; let $\cnu'$ be the maximal element of $\cLambda$,
such that this quotient is non-zero when restricted to $_{\cnu'}^1\BunNb$.
Then this restriction (either *- or !-) is a perverse sheaf, and its further
restriction onto the locally closed sub-stack of $_{\cnu'}^1\BunN$
equal to $(\on{ev}_{\cnu'})^{-1}(N^-\cdot w_0)$, is a local system.

Hence, we deduce that the Euler characteristic of the *-restriction of
$\fIC_\cnu^\psi$ to some $(\on{ev}_{\cnu'})^{-1}(N^-\cdot w_0)$ with $\cnu'\neq \cnu$
is non-zero. We are going to show that this is impossible by comparing the present
situation with the one for $\Pervgrw$.

\medskip

Let us recall that for any $\cmu\in \cLambda^+$, the perverse sheaf
$\IC^\psi_{\Gr}\star \IC_{\cmu,\Gr}$ is irreducible and is isomorphic to
the clean extension of the character sheaf on the $I^-$-orbit of the point
$w_0\cdot (\cmu+\crho')\in \Gr$, by \thmref{Cass-Shal for I}.

\medskip

We have the convolution functor
$$\Pervgrw\times \Pervsmg\to \Pervsmw.$$

\begin{thm} \label{convolution of Whittaker}
$\IC^\psi_{\Gr}\star \fIC_\cnu=\fIC^\psi_{\cnu+\crho'}$.
\end{thm}

We omit the proof, since it essentially repeats the proof of \thmref{conv is simple},
where instead of the fact that $\CL^w \star \IC_{\cmu,\Gr}$ is irreducible for
$\cmu\in \cLambda^+$, we use the above mentioned fact about
$\IC^\psi_{\Gr}\star \IC_{\cmu,\Gr}\in \Pervgrw$.

\medskip

We claim that the fiber of  $\IC^\psi_{\Gr}\star \fIC_\cnu$ at a point of
$(\on{ev}_{\cnu'})^{-1}(N^-\cdot w_0)$ can be written as an extension of certain complexes
$\CK_{\cnu''}$, and the fiber of $\IC^\psi_{\Gr}\star \IC_{\cmu,\Gr}$ at
a point of $I^-\cdot (w_0\cdot (\cmu'+\crho'))$ for $\cnu-\cnu'=\cmu-\cmu'$ can be
written as an exetension of the same complexes.

This would imply our assertion about Euler characteristics, since the fibers
of the convolution $\IC^\psi_{\Gr}\star \IC_{\cmu,\Gr}$ over
$I^-\cdot (w_0\cdot (\cmu'+\crho'))$ are zero
unless $\cmu'=\cmu$ by cleanness.

\bigskip

For $\cnu''$ the complex $K_{\cnu''}$ is defined as the fiber of the direct image
under
$$\hl'_G:(\hl'_G)^{-1}\Bigl({}_{\cnu'}^1\BunNb\Bigr)\cap
(\hr'_G)^{-1}\Bigl({}_{\cnu''}\BunNb\Bigr)\cap {}^{1,0}\CH^{\crho'}_{G,N^-,x}\to
{}_{\cnu'}^k\BunNb$$
of the *-restriction of $\IC^\psi_{\Gr}\tboxtimes \fIC_\cnu$ to the above substack.

Hence, $K_{\cnu''}$ is the cohomology with compact supports along the scheme
$$N((t))\cdot w_0\cdot (\cnu'-\cnu'')\cap I^-\cdot (w_0\cdot \crho')\subset \Gr^{\crho'}$$
of the complex equal to the tensor product of the character sheaf along
$I^-\cdot (w_0\cdot \crho')$ and the constant complex equal to the
stalk of $\fIC_{\cnu}$ on $_{\cnu''}\BunN$.

\medskip

Let us now calculate the fiber of $\IC^\psi_{\Gr}\star \IC_{\cmu,\Gr}$
at a point of $I^-\cdot (w_0\cdot \cmu')$ for $\cmu$ large and
$\cnu-\cnu'=\cmu-\cmu'$. For that we will intersect the fiber of the
convolution diagram over $w_0\cdot \cmu'$ with the subschemes of
the form $$\Bigl(I^-\cdot (w_0\cdot \crho)\Bigr)\star \Gr^{\cmu''}.$$

As we saw above, each of these intersections is isomorphic to
$$N((t))\cdot w_0\cdot (\cmu'-\cmu'') \cap I^-\cdot (w_0\cdot \crho)\subset \Gr^{\crho'}.$$
For each such $\cmu''$ the complex that we have to integrate is
the tensor product of the character sheaf along $I^-\cdot (w_0\cdot \crho')$
and the stalk of $\IC_{\cmu,\Gr}$ at $\Gr^{\cmu''}$.

\medskip

We set up the bijection between $\cnu''$ and $\cmu''$ so that
$\cnu''-\cnu'=\cmu''-\cmu'$. Our assertion follows from the
fact that for $\clambda$ small comared with $\cmu$ and $\cnu$
the stalk of $\fIC_{\cnu}$ on $_{\cnu-\clambda}\BunN$ is
isomorphic to the stalk of $\IC_{\cmu,\Gr}$ on $\IC_{\cmu-\clambda,\Gr}$.
This follows by combining \cite{FFKM,BFGM} and \cite{Lu,Soe}.

\ssec{Action of convolution on standard objects}

\sssec{}

We will now prove the following assertion, parallel to \corref{lambda invariance}:

\begin{prop} \label{semiinf lambda inv}
If $\clambda$ is dominant there is a canonical isomorphism
$$j_{!,\clambda}\star \nabla_\cnu\simeq \nabla_{\cnu+\clambda}.$$
\end{prop}

\begin{proof}

Using \propref{convolution respects}, it is sufficient to show that
the stalk of $j_{!,\clambda}\star \nabla_\cnu$ is $0$ on any $_{\tw'}^I\Bun_N$
for $\tw'\neq \clambda+\cnu$, and that it is canonically $\BC$ the
latter case. This follows in a rather straightforward way from the definition of
convolution.

Consider the stack
$$(\hl'_G)^{-1}\Bigl({}_{w'\cdot \cnu'}^I\BunNb\Bigr)\cap
(\hr'_G)^{-1}\Bigl({}_{w\cdot \cnu}^I\BunNb\Bigr)\cap {}^{I,I}\CH^\clambda_{G,N^-,x},$$
projecting to $_{w'\cdot \cnu'}^I\BunNb$ by means of $\hl'_G$. In the
above formula $\CH^\clambda_{G,N,x}$ is the locally closed substack of
$^{I,I}\CH^\clambda_{G,N^-,x}$, corresponding to the $I$-orbit
$I\cdot \clambda\subset \Fl$.

The fiber of the above stack over a point of $_{w'\cdot \cnu'}^I\BunNb$ is
isomorphic to
\begin{equation} \label{Iwahori intersection}
\Bigl(N^-((t))\cdot (\cnu'-\cnu)\cdot w^{-1}\Bigr) \cap
\Bigl((w')^{-1}\cdot I\cdot \clambda\Bigr) \subset \Fl.
\end{equation}

Set $w=1$, and we claim that the above intersection is empty unless
$\cnu'=\cnu+\clambda$ and $w'=1$, and that in the latter case, this is a
point-scheme.

\medskip

The latter assertion is evident. To prove the first one, we will use the following:

\begin{lem}
For $\clambda$ dominant,
$$N^-((t))\cdot B[[t]] \supset \on{Ad}_{\clambda}(I)\subset B[[t]]\cdot N^-((t)).$$
\end{lem}

Using the lemma, it is enough to show that
$$\Bigl(w'\cdot (\cnu'-\cnu)\cdot N^-((t))\Bigr)\cap \Bigl(N^+((t))\cdot \clambda\Bigr)\subset G((t))$$
is non-empty only if $w'=1$ and $\cnu'-\cnu=\clambda$, which is evident from the Bruhat
decomposition.

\end{proof}

\sssec{}

Let us now exhibit a compatibility relation between the isomorphisms of
\propref{semiinf lambda inv} and \propref{sph conv}. Namely, we claim that
for $\clambda\in \cLambda^+$ the diagrams
\begin{equation}  \label{Miura on semiinf}
\CD
\IC_{\Gr,\clambda}\star \fIC_\cnu @>>> j_{*,\clambda}\star \fIC_\cnu @>>>
j_{*,\clambda}\star \Delta_\cnu \\
@V{\sim}VV  & & @V{\sim}VV \\
\underset{\cmu}\oplus\, \fIC_{\cnu+\cmu}\otimes \uV^\clambda(\cmu) @>>>
\fIC_{\cnu+\clambda} @>>> \Delta_{\cnu+\clambda}.
\endCD
\end{equation}
and
\begin{equation} \label{ind sys1}
\CD
j_{*,\clambda}\star \fIC_\cnu @>>>
j_{*,\clambda}\star \IC_{\cmu,\Gr}\star \IC_{-w_0(\cmu),\Gr} \star \fIC_\cnu \\
@VVV   @VVV  \\
j_{*,\clambda}\star \Delta_\cnu & & j_{*,\clambda}\star \IC_{\cmu,\Gr}\star \fIC_{\cnu-\cmu} \\
@V{\sim}VV   @VVV  \\
\Delta_{\cnu+\clambda} @<<<  j_{*,\clambda+\cmu}\star \fIC_{\cnu-\cmu} ,
\endCD
\end{equation}
are commutative. This follows from the definition of the isomorphisms in both
cases.

\medskip

Note that \propref{semiinf lambda inv} implies that for $\cnu$ dominant
$$j_{!,-\cnu}\star \nabla_{w_0\cdot \cnu'}\simeq \nabla_{w_0\cdot (\cnu'-w_0(\cnu))},$$
and hence
\begin{equation} \label{other lambda inv}
j_{*,\cnu}\star \nabla_{w_0\cdot \cnu'}\simeq \nabla_{w_0\cdot (\cnu'+w_0(\cnu))}.
\end{equation}

Consider now the morphism
\begin{equation}   \label{socle of standard, semiinf}
\fIC_{\cnu}\to \nabla_{w_0\cdot (\cnu+2\crho)},
\end{equation}
obtained by Verdier duality from \secref{cosocle of other costandard, semiinf}.
By construction, the space of such morphisms for every $\cnu$ is a
$1$-dimensional vector space, canonically independent of $\cnu$. From the
construction one infers the following:

\begin{lem} \label{commutation for other st}
For $\cnu\in \cLambda$, $\clambda,\cmu\in \cLambda^+$ the diagrams
$$
\CD
\IC_{\clambda,\Gr}\star \fIC_\cnu @>>> j_{*,\clambda}\star \fIC_\cnu @>>>
 j_{*,\clambda}\star \nabla_{w_0\cdot (\cnu+2\crho)} \\
@V{\sim}VV & & @V{\sim}VV \\
\underset{\cnu'}\oplus\, \fIC_{\cnu+\cnu'}\otimes \uV^\clambda(\cnu') @>>>
\fIC_{\cnu+w_0(\clambda)} @>>> \nabla_{w_0\cdot (\cnu+2\crho+w_0(\clambda))}
\endCD
$$
and
$$
\CD
j_{*,\clambda}\star \fIC_\cnu @>>>
j_{*,\clambda}\star \IC_{\cmu,\Gr}\star \IC_{-w_0(\cmu),\Gr}\star \fIC_\cnu  \\
@VVV  @VVV \\
j_{*,\clambda}\star \nabla_{w_0\cdot (\cnu+2\crho)} & &
j_{*,\clambda}\star  \IC_{\cmu,\Gr}\star \fIC_{\cnu-w_0(\cmu)} \\
@V{\sim}VV  @VVV  \\
\nabla_{w_0\cdot (\cnu+2\crho+w_0(\clambda))}  @<<<
j_{*,\clambda+\cmu}\star  \fIC_{\cnu-w_0(\cmu)}
\endCD
$$
are commutative.
\end{lem}

\section{The equivalence}

\ssec{The functor}

\sssec{}

Let $\tCS$ be an object of $\Catgeomdk$. We attach to it a covariant functor
on $\Pervsmkb$ as follows. To an object $\CF\in \Pervsmkb$ we assign
the set of collections of morphisms $\tCS_\clambda\star \fIC_{\clambda}\to \CF$, such that for
any $V\in \Rep(\cG)$ and $\cmu\in \cLambda$, the diagram
$$
\CD
\tCS_\clambda \star \CV \star \fIC_{\clambda-\cmu}\otimes (\uV(\cmu))^*  @>>>
\left(\tCS_{\clambda-\cmu}\otimes \uV(\cmu)\right)
\star \fIC_{\clambda-\cmu}\otimes (\uV(\cmu))^*  \\
@VVV   @VVV  \\
\tCS_\clambda \star \fIC_{\clambda} \otimes \uV(\cmu)\otimes (\uV(\cmu))^*  & &
\tCS_{\clambda-\cmu}\star \fIC_{\clambda-\cmu} \\
@VVV  @VVV \\
\tCS_\clambda \star \fIC_{\clambda} @>>> \CF
\endCD
$$
commutes, where the upper horizontal arrow is given by the Hecke morphism for
$\tCS$, and the left vertical arrow by \propref{sph conv}.

\medskip

It is easy to see that the above functor is representable by
$$\on{co-eq}\Bigl(\underset{\clambda,\cmu,V}\oplus\,
\tCS_\clambda \star \CV \star \fIC_{\clambda-\cmu}\otimes (\uV(\cmu))^* \rightrightarrows
\underset{\cnu}\oplus\, \tCS_\cnu \star \fIC_{\cnu}\Bigr),$$
where the two arrows correspond to the two circuits of the above commutative
diagram.

\medskip

We denote the resulting functor $\Catgeomdk\to \Pervsmkb$ by
$\ConvHecke$. By construction, $\ConvHecke$ is right-exact.

\begin{prop}  \label{convolution of induced}
For  $\tCS=\CS\star \tRg\{\cmu\}\in \Catgeomdk$ the object
$\ConvHecke(\tCS)$ is canonically isomorphic to
$\CS\star \fIC_{-\cmu}$.
\end{prop}

\begin{proof}

For a morphism $\ConvHecke(\tCS)\to \CF$,
by taking its component
$\tCS_{-\cmu}\star \fIC_{-\cmu}\to \CF$ we obtain a map
$\CS\star \fIC_{-\cmu}\to \CF$, since
$$\tCS_{-\cmu}\simeq \CS\star \tRg\{\cmu\}_{-\cmu}\simeq \CS\star (\tRg)_0,$$
and it contains $\CS$ as a direct summand.

\medskip

Vice versa, having a map $\CS\star \fIC_{-\cmu}\to \CF$,
for every $V\in \Rep(\cG)$ and $\clambda$ we
define a map
$$\left(\CS\star \CV \otimes \uV^*(\clambda+\cmu)\right)
\star \fIC_{\clambda}\to \CF$$ by
$$\left(\CS\star \CV \otimes \uV^*(\clambda+\cmu)\right)\star \fIC_{\clambda}\to
\CS\star \fIC_{-\cmu}\otimes \uV(-\clambda-\cmu)\otimes \uV^*(\clambda+\cmu)
\to \CS\star \fIC_{-\cmu}\to \CF.$$

The fact that the resulting system of maps satisfies the defining
condition follows from the second assertion in \propref{sph conv}.

\end{proof}

We also have the following assertion that follows from \propref{sph conv}(2):

\begin{lem}
For $\tCS\in \Catgeomdk$ and $\cmu\in \cLambda$,
$$\ConvHecke(\tCS)\simeq
\on{co-eq}\Bigl(\underset{\uV}\oplus\, \tCS_{\cmu}\star \CV\star \fIC_{\cmu}\otimes V(0)^*
\rightrightarrows \tCS_{\cmu}\star \fIC_{\cmu}\Bigr).$$
\end{lem}

\sssec{}

We propose the following:

\begin{conj}  \label{general}
The functor $$\ConvHecke:\Catgeomdk\to \Pervsmkb$$ is exact and fully-faithful.
\end{conj}

In fact, we think that $\ConvHecke$ is very close to be an equivalence of categories.
Unfortunately, we cannot formulate a precise conjecture, due to our lack of
understanding of Noetherian properties of both categories. In  any case, we think
that one can express $\Pervsmk$ completely in terms of $\Catgeomdk$, which would
then supply a local (in particular, independent of the global curve $X$) description of
$\Pervsmk$.

\medskip

In what follows we are going to discuss a version of the above conjecture, where instead
of the level $G^k$ we take $I^0$. In this case it would be possible to formulate and
prove a more precise result.

\begin{thm}  \label{main}
The functor $$\ConvHecke:\Catgeomdni\to \Pervsmnib$$ is exact,
and it defines an equivalence between the sub-categories of
Artinian objects on both sides.
\end{thm}

Since the subcategory $\Catgeomdni_{Art}$ of Artinian objects in $\Catgeomdni$
is equivalent to $\tu_\ell\mod_0$, as a corollary we obtain:

\begin{thm}   \label{main quantum}
The category $\tu_\ell\mod_0$ is equivalent to the category of
Artinian objects in $\Pervsmni$.
\end{thm}

\sssec{}  \label{BK action}

Here we would like to add the following observation.

As we saw above, the category $\Catgeomdni$ is acted on by the group
$W_{aff}$ by self-equivalences: the elements of $\cLambda$ act by
shifting the grading, and $w\in W$ by the twisting functors
$\tCS\mapsto {}^w\tCS$ (which on the level of $\tu_\ell\mod$ correspond
to the functors $\sF_w$). Evidently, these functors
preserve the subcategory $\Catgeomdni_{Art}$, and, hence, the carry over
to the category of Artinian objects in $\Pervsmni$.

Let us describe how these functors act on the irreducibles of $\Pervsmni$.
For $w\in W$ let $\CL^w=\IC_{w\cdot \clambda,\Gr}$ be the corresponding
"restricted" irreducible in $\Pervgri$. By \thmref{conv is simple}
and \propref{convolution of induced},
\begin{equation} \label{image of simples}
\ConvHecke(\CL^w\star \tRg\{\cmu\})\simeq \fIC_{w\cdot (\clambda-\cmu)}.
\end{equation}

Hence,
$$\Bigl(\fIC_{w\cdot \cnu}\Bigr)\{\cmu\}\simeq \fIC_{w\cdot (\cnu-\cmu)}$$
and
\begin{equation} \label{Fourier}
\sF_{w'}\Bigl(\fIC_{w\cdot (\clambda-\cmu)}\Bigr)\simeq \fIC_{w\cdot (\clambda-w'(\cmu))}.
\end{equation}

\medskip

Recall that the $\BC$-linearized Grothendieck group of the category of Artinian objects in
$\Pervsmni$ identifies with Lusztig's periodic module over the affine Hecke algebra
(cf. \cite{FFKM})\footnote{For this to be formally true we have to pass to the category
of mixed D-modules of Hodge-Tate type in $\Pervsmni$},
and hence, also with the space of Iwahori-invariant functions in
the Schwarz space of \cite{BK}. Equation \eqref{Fourier} implies that the maps
on the Grothendieck group, induced by the functors $\sF_w$, are equal to the
Fourier transform operators, introduced in \cite{BK}.

\medskip

The rest of the paper is devoted to the proof of \thmref{main}.



\ssec{Proof of the equivalence}

\sssec{}

As a first step we prove the following:

\begin{prop}  \label{exactness of Hecke conv}
The functor $$\ConvHecke:\Catgeomdni\to \Pervsmnib$$ is exact.
\end{prop}

The present subsection is devoted to the proof of this proposition.

\medskip

Since $\Catgeomdni$ is the ind-completion of the subcategory of its Artinian objects, it is
sufficient to prove that $\ConvHecke$ restricted to $\Catgeomdni_{Art}$
is exact.

Let
$$0\to \tCS_1\to \tCS_2\to \tCS\to 0$$
be a short exact sequence of objects of $\Catgeomdni_{Art}$.
We have to show that $\ConvHecke(\tCS_1)\to \ConvHecke(\tCS_2)$
is injective. For that we may assume that $\tCS$ is simple. By
\secref{char of simples on gr}, $\tCS$ is then isomorphic to
$\CS\star \tRg\{\cmu\}$ for $\CS\in \Pervgrni$.

\medskip

We can find an object $\CS'\in \Pervgrni$ with a surjection $\CS'\twoheadrightarrow \CS$,
and a map $\CS'\to (\tCS_2)_{-\cmu}$ in $\Pervgrnib$, such that the diagram
$$
\CD
(\tCS_2)_{-\cmu} @>>> (\tCS)_{-\cmu} \\
@AAA  @AAA \\
\CS' @>>> \CS
\endCD
$$
is commutative. Hence, we obtain a map $\CS'\star \tRg\{\cmu\}\to \tCS_2$. Let $\tCS'_2$
be the Cartesian product of $\tCS_2$ and $\CS'\star \tRg\{\cmu\}$ over
$\CS\star \tRg\{\cmu\}$. We have a commutative diagram:
$$
\CD
& & & &  0 & & 0 \\
& & & & @AAA @AAA \\
0 @>>>  \tCS_1 @>>>  \tCS_2 @>>> \CS\star \tRg\{\cmu\} @>>> 0 \\
& & @A{\on{id}}AA @AAA  @AAA & \\
0 @>>>  \tCS_1 @>>>  \tCS'_2 @>>> \CS'\star \tRg\{\cmu\} @>>> 0 \\
& & & & @AAA @AAA \\
& &  & & \CS''\star \tRg\{\cmu\} @>{\on{id}}>> \CS''\star \tRg\{\cmu\} \\
& & & & @AAA @AAA \\
& & & & 0 & & 0 .\\
\endCD
$$

It is enough to show that the map
$$\ConvHecke\Bigl(\tCS_1 \oplus \CS''\star \tRg\{\cmu\} \Bigr)\to
\ConvHecke\Bigl( \tCS'_2\Bigr)$$ is injective. However, by construction,
$\tCS'_2$ splits as a direct sum $\tCS_1\oplus \CS'\star \tRg\{\cmu\}$.
Hence, it is enough to show that the map
$$\ConvHecke(\CS''\star \tRg\{\cmu\})\to \ConvHecke(\CS'\star \tRg\{\cmu\})$$
is injective. But the latter results from \propref{convolution of induced} combined
with \thmref{exactness of convolution}, since the map in question comes from
a map $\CS''\to \CS'$ in $\Pervgrni$.

\sssec{}

Recall that the Verdier duality functor $\BD$ is defined on $\Catgeomdni_{Art}$.
In this subsection we will prove the following:

\begin{prop}  \label{Hecke conv and duality}
The functor $\ConvHecke$ commutes with the Verdier duality.
\end{prop}

\medskip

Recall that if $\tCS\in \Catgeomdni_{Art}$ is an object represented as
$$\on{coker}\Bigl(\CS_1\star \tRg\{\cmu_1\}\to \CS_2\star \tRg\{\cmu_2\}\Bigr),$$
then $\BD(\tCS)$ is described as follows:

The map $\CS_1\star \tRg\{\cmu_1\}\to \CS_2\star \tRg\{\cmu_2\}$ comes from
a map $\alpha:\CS_1\to \CS_2\star \CV\otimes \uV^*(\cmu_2-\cmu_1)$ defined
for some $V\in \Rep(\cG)$. By adjunction, we have a map
$$\CS_1\star \BD(\CV^{op})\otimes \uV(\cmu_1-\cmu_2)\to \CS_2,$$
and applying the Verdier duality we obtain a map
$$\BD(\alpha): \BD(\CS_2)\to \BD(\CS_1)\star \CV^{op} \otimes \uV^*(\cmu_2-\cmu_1).$$

Recall that the functor $\CV\mapsto \BD(\CV^{op})$ corresponds on the level
of $\Rep(\cG)$ to the dualization functor $V\mapsto V^*$, whereas
$\CV\mapsto \BD(\CV)$ corresponds to the contragredient duality $V\mapsto V^\vee$.
In particular, $\uV(\cmu)\simeq \uV^{op}(-\cmu)$.

We then obtain a morphism
$$\BD(\CS_2)\star \tRg\{\cmu_2\}\to \BD(\CS_1)\star \tRg\{\cmu_1\},$$ whose kernel ia
$\BD(\tCS)$.

\medskip

For $\tCS$ as above,
by \propref{convolution of induced}, $\ConvHecke(\tCS)\simeq \on{coker}(\beta)$, where
$\beta$ is the map
\begin{align*}
&\CS_1\star \fIC_{-\cmu_1}\to (\CS_2\star \CV\star \fIC_{-\cmu_1})\otimes
\uV^*(\cmu_2-\cmu_1)\to  \\
&\CS_2 \star \fIC_{-\cmu_2}\otimes \uV(\cmu_1-\cmu_2)
\otimes \uV^*(\cmu_2-\cmu_1)\to \CS_2 \star \fIC_{-\cmu_2}.
\end{align*}

By \propref{exactness of Hecke conv},
$\ConvHecke(\BD(\tCS))\simeq \on{ker}(\gamma)$, where $\gamma$ is the map
\begin{align*}
&\BD(\CS_2)\star \fIC_{-\cmu_2}\to
(\BD(\CS_1)\star \CV^{op} \star \fIC_{-\cmu_2})\otimes \uV^*(\cmu_2-\cmu_1)\to \\
&\BD(\CS_1)\star \fIC_{-\cmu_1} \otimes \uV^{op}(\cmu_2-\cmu_1)
\otimes \uV^*(\cmu_2-\cmu_1)
\to \BD(\CS_1)\star \fIC_{-\cmu_1}.
\end{align*}

To prove the proposition it remains to see that the morphisms $\beta$ and
$\gamma$ are transformed into one-another by Verdier duaility. This is evident
when $V$ is the trivial representation. By transitivity, this reduces the assertion
to the case when $\CS_1\simeq \CS_2\star \CV\otimes \uV^*(\cmu_2-\cmu_1)$.

In the latter case, both arrows $\BD(\beta)$ and $\gamma$ are obtained from the
corresponding arrows for $\CS_2$ replaced by $\delta_{1,\Gr}$ by convolution with
$\CS_2$. The case $\CS_2=\delta_{1,\Gr}$ is a straightforward verification.

\sssec{}

We will now state a crucial result, from which we will deduce \thmref{main}.

\begin{thm}  \label{Baby to baby}
For $w\in W$ and $\cmu\in \cLambda$,
$$\ConvHecke\Bigl(\tCM^{w\cdot \cmu}\Bigr) \simeq \Delta_{w\cdot \cmu}.$$
\end{thm}

We will now deduce \thmref{main} from \thmref{Baby to baby}. Consider now the following
general set-up:

\medskip

\noindent Let $\CC$ be an abelian Artinian category; let $A$ be the set parametrizing
its irreducibles; for $a\in A$ we will denote by $\CL^a$ the corresponding object.
Assume also that for each $a\in A$ there exist objects $\nabla^a$ and $\Delta^a$,
such that $\CL^a$ is the cosocle of $\nabla^a$ and the socle of $\Delta^a$.
Assume, moreover, that $Ext^i(\nabla^{a'},\Delta^{a''})=0$ for $i=1,2$, and
$Hom(\nabla^{a'},\Delta^{a''})=0$ unless $a'=a''$, and in the latter case it is
$1$-dimensional (which implies that any element in $Hom(\nabla^1,\Delta^a)$
factors through $\CL^a$).

Let now $\CC_1$ and $\CC_2$ be two such categories with the same set of
irreducibles $A$. Let $\sG:\CC_1\to \CC_2$ be an exact functor, such that
$\sG(\CL^a_1)\simeq\CL^a_2$, $\sG(\nabla^a_1)\simeq \nabla^a_2$,
$\sG(\Delta^a_1)\simeq \Delta^a_2$.

\begin{lem}   \label{hered}
Under the above circumstances, $\sG$ is an equivalence of categories.
\end{lem}

\thmref{main} follows from this lemma, using \corref{no Ext 2}, \lemref{geom hered},
\propref{exactness of Hecke conv},
\eqref{image of simples}, \propref{Baby to baby} and \propref{Hecke conv and duality}.

\sssec{Proof of \lemref{hered}}  Note first of all that the assumption implies
that $\sG$ is faithful.

\smallskip

\noindent{\bf Step 1.} For $a,a'\in A$ consider the long exact sequences
\begin{align*}
&0\to Hom(\CL_i^a,\Delta_i^{a'})\to Hom(\nabla_i^a,\Delta_i^{a'})\to
Hom(\on{ker}(\nabla_i^a\to \CL_i^a),\Delta_i^{a'})\to \\
&Ext^1(\CL_i^a,\Delta_i^{a'})\to Ext^1(\nabla_i^a,\Delta_i^{a'})=0
\end{align*}
for $i=1,2$. Since $Hom(\nabla_1^a,\Delta_1^{a'})\to Hom(\nabla_2^a,\Delta_2^{a'})$
is an isomorphism, comparing the two, we infer that
$Ext^1(\CL_1^a,\Delta_1^{a'})\to Ext^1(\CL_2^a,\Delta_2^{a'})$ is injective.

\smallskip

\noindent{\bf Step 2.}
Consider now the long exact sequence
\begin{align*}
&0\to Hom(\CL_1^a,\CL^{a'}_1)\to Hom(\CL_1^a,\Delta^{a'}_1)\to
Hom(\CL_1^a,\Delta^{a'}_1/\CL^{a'}_1)\to \\
&Ext^1(\CL_1^a,\CL^{a'}_1)\to Ext^1(\CL_1^a,\Delta^{a'}_1)\to
Ext^1(\CL_1^a,\Delta^{a'}_1/\CL^{a'}_1)
\end{align*}
for $i=1,2$. $Hom(\CL_1^a,\Delta_1^{a'})\to Hom(\CL_2^a,\Delta_2^{a'})$
is an isomorphism and using Step 1, we find that
$Ext^1(\CL_1^a,\CL^{a'}_1)\to Ext^1(\CL_2^a,\CL^{a'}_2)$ is injective.

\smallskip

\noindent{\bf Step 3.}
Let $\CF'$ be any object of $\CC_1$. Using Step 3, by induction on the length
of $\CF'$, we find that the map
$Hom(\CL_1^a,\CF')\to Hom(\CL^a_2,\sG(\CF'))$ is an isomorphism.

\smallskip

\noindent{\bf Step 4.}
Returning to the long exact sequence of Step 1, we find that the map
$Ext^1(\CL_1^a,\Delta_1^{a'})\to Ext^1(\CL_2^a,\Delta_2^{a'})$ is an
isomorphism.

\smallskip

\noindent{\bf Step 5.}
Again, by induction on the length, using Step 3, we show that the map
$Ext^1(\CL_1^a,\CF')\to Ext^1(\CL^a_2,\sG(\CF'))$ is injective.

\smallskip

\noindent{\bf Step 6.}
By the exact sequence of Step 2, from Step 4 we find that
$Ext^1(\CL_1^a,\CL^{a'}_1)\to Ext^1(\CL_2^a,\CL^{a'}_2)$ is an isomorphism.

\smallskip

\noindent{\bf Step 7.}
Let $\CF$ be an object of $\CC_1$, and $\CF'$ some other object. By
induction on the length of $\CF$, from Step 5 we obtain that
$Hom(\CF,\CF')\to Hom(\sG(\CF),\sG(\CF'))$ is an isomorphism.

\medskip

Hence, $\sG$ is fully-faithful. To finish the proof of the lemma, we have
to show that $\sG$ induces isomorphsims on the level of $Ext^1(\cdot,\cdot)$.

\medskip

\noindent{\bf Step 8.}
By induction on the length of $\CF$, from Step 5 and Step 7 we obtain that
$Ext^1(\CF,\CF')\to Ext^1(\sG(\CF),\sG(\CF'))$ is injective.

\smallskip

\noindent{\bf Step 9.}
For $a,a'\in A$ consider the long exact sequences
\begin{align*}
&...0=Ext^1(\nabla_i^a,\Delta_i^{a'})\to
Ext^1(\on{ker}(\nabla_i^a\to \CL_i^a),\Delta_i^{a'})\to \\
&Ext^2(\CL_i^a,\Delta_i^{a'})\to Ext^2(\nabla_i^a,\Delta_i^{a'})=0
\end{align*}
for $i=1,2$. From Step 8 we infer that
$Ext^2(\CL_1^a,\Delta_1^{a'})\to Ext^2(\CL_2^a,\Delta_2^{a'})$ is
injective.

\smallskip

\noindent{\bf Step 10.}
Consider the long exact sequence
\begin{align*}
&Ext^1(\CL_1^a,\CL^{a'}_1)\to Ext^1(\CL_1^a,\Delta^{a'}_1)\to
Ext^1(\CL_1^a,\Delta^{a'}_1/\CL^{a'}_1)\to \\
&Ext^2(\CL_1^a,\CL^{a'}_1)\to Ext^2(\CL_1^a,\Delta^{a'}_1)\to
Ext^2(\CL_1^a,\Delta^{a'}_1/\CL^{a'}_1)
\end{align*}
By Step 4, Step 8 and Step 9, the map
$Ext^2(\CL_1^a,\CL^{a'}_1)\to Ext^2(\CL_2^a,\CL^{a'}_2)$ is injective.

\smallskip

\noindent{\bf Step 11.}
By induction on length, from Step 6, we obtain that
$Ext^1(\CL_1^a,\CF')\to Ext^1(\CL^a_2,\sG(\CF'))$ is an isomorphism.

\smallskip

\noindent{\bf Step 12.}
Again, by induction on the length, from Step 10 and Step 6, we obtain that the map
$Ext^2(\CL_1^a,\CF')\to Ext^2(\CL^a_2,\sG(\CF'))$
is injective.

\smallskip

\noindent{\bf Step 13.}
Finally, by induction on the length of $\CF$, from Steps 11 and 12 we infer that the
map $Ext^1(\CF,\CF')\to Ext^1(\sG(\CF),\sG(\CF'))$ is an isomorphism.

\ssec{Identification of the image of baby co-Verma modules}

\sssec{}

In this subsection we will prove \thmref{Baby to baby}. Note that it suffices
to show that
$$\ConvHecke(\tCM^1)\simeq \Delta_0,$$
since all other isomorphisms will then hold by \eqref{Baby Vermas from one another},
\eqref{lambda invariance} and \eqref{semiinf lambda inv}.

\medskip

We construct a map
\begin{equation} \label{desired map}
\ConvHecke(\tCM^1)\to \Delta_0
\end{equation}
as follows. We need to construct the maps
$$\underset{\clambda}{\underset{\longrightarrow}{lim}}\,
j_{*,\clambda+\cmu}\star \IC_{-w_0(\clambda),\Gr}\star \fIC_{-\cmu}\to \Delta_0$$
for every $\cmu$.

For $\clambda\in \cLambda^+$ as above we have a map
\begin{align*}
&j_{*,\clambda+\cmu}\star \IC_{-w_0(\clambda),\Gr}\star \fIC_{-\cmu} \simeq
\underset{\cnu}\oplus \, j_{*,\clambda+\cmu} \star \fIC_{-\cmu+\cnu}\otimes \uV^*(\cnu)
\to j_{*,\clambda+\cmu}\star \fIC_{-\cmu-\clambda}\to \\
& j_{*,\clambda+\cmu}\star \Delta_{-\cmu-\clambda}\simeq \Delta_0,
\end{align*}
The fact that these maps are compatible with the maps in the inductive system
that defines $\tCM^1$, follows from the commutativity of the diagrams
\eqref{Miura on semiinf} and \eqref{ind sys1}. The fact that the resulting system of maps
$$\tCM^1_\cmu\star \fIC_{\cmu}\to \Delta_0$$
factors through $\ConvHecke(\tCM^1)$ follows from \eqref{comm sph2}.

\sssec{}

Now, we claim that the map $\ConvHecke(\tCM^1)\to \Delta_0$ constructed above
is non-zero in the quotient category $\Pervsmnif$.
Using \propref{non-degen quotient on semi-inf}, it is enough to show that the map
$$\on{Av}_{N^-,\psi}\Bigl(\ConvHecke(\tCM^1)\Bigr)\to \on{Av}_{N^-,\psi}(\Delta_0)$$
is non-zero. The latter reduces to showing that for $\clambda$ dominant and regular,
the map
$$j_{*,\clambda}\star \fIC_{\cnu}\to \Delta_{\clambda+\cnu}$$
gives rise to a non-zero map
$$\on{Av}_{I^-,\psi}(\CW^{*,\clambda})\star \fIC_{\cnu}\to
\on{Av}_{N^-,\psi}(\Delta_{\clambda+\cnu}).$$
However, the latter is straightforward from the definition of
convolution.

\medskip

In particular, by \corref{cosocle on semiinf}(2), we obtain that the map
of \eqref{desired map} is surjective. Moreover, it is an
isomorphism in the quotient category $\Pervsmnif$ by
\propref{geometric cosocles}(1).

We claim that in order to finish the proof of the theorem, it suffices
to show that there exists a non-zero map
\begin{equation} \label{desired opp dir}
\Delta_0\to \ConvHecke(\tCM^1).
\end{equation}

Indeed, if such a map exists, its image in $\Pervsmnif$ cannot be $0$
by \corref{cosocle on semiinf}, and hence the composition
$$\Delta_0\to \ConvHecke(\tCM^1)\to \Delta_0$$
is non-zero. Then the above composition is the identity map on
$\Delta_0$, up to a scalar.

\medskip

Hence, it would remain to show that $\ConvHecke(\tCM^1)$ is indecomposable.
We claim that it in fact does not admit irreducible quotients besides the
canonical map
$$\ConvHecke(\tCM^1)\to \ConvHecke(\CL^{w_0}\star \tRg\{\crho'\}).$$
This is so because $\ConvHecke(\tCM^1)$ cannot map to any partially
integrable irreducible object of $\Pervsmni$ by the same argument as in the proof
of \propref{cosocle of 0 costandard}, and by \corref{geometric cosocles}(1),
$\ConvHecke(\CL^{w_0}\star \tRg\{\crho'\})$
is the only non-partially integrable constituent of $\ConvHecke(\tCM^1)$.

\sssec{}

Thus, our goal is to construct a map as in \eqref{desired opp dir}.
By \propref{dual descr}
and \propref{Hecke conv and duality}, it suffices to construct a map
$$\ConvHecke\Bigl(({}^{w_0}\tCM^{w_0})\{2\crho\}\Bigr)\to \nabla_0,$$
or, equivalently, a map
$$\ConvHecke\Bigl(({}^{w_0}\tCM^1)\{2\crho\}\Bigr)\to \nabla_{w_0}.$$

Consider the inductive system that defines $\left(({}^{w_0}\tCM^1)\{2\crho\}\right){}_\cmu$,
viewed as an object of $\Pervgrnib$:
$$\underset{\clambda}{\underset{\longrightarrow}{lim}}\,
j_{*,\clambda+\crho'-w_0(\cmu)+2\crho}\star \IC_{-w_0(\clambda),\Gr}.$$

For every such $\cmu$ and $\clambda$,  we define the map
$$j_{*,\clambda-w_0(\cmu)+2\crho}\star \IC_{-w_0(\clambda),\Gr}\star \fIC_{\cmu}\to
\nabla_{w_0}$$
as the composition:
\begin{align*}
&j_{*,\clambda-w_0(\cmu)+2\crho}\star \IC_{-w_0(\clambda),\Gr}\star \fIC_{\cmu}\simeq
\underset{\cnu}\oplus\, j_{*,\clambda+2\crho-w_0(\cmu)}\star \fIC_{\cmu+\cnu}\otimes
(\uV^\clambda)^*(\cnu)\to \\
& \to j_{*,\clambda+2\crho-w_0(\cmu)}\star \fIC_{\cmu-w_0(\clambda)}
\to j_{*,\clambda+2\crho-w_0(\cmu)}\star \nabla_{w_0\cdot (\cmu-w_0(\clambda)+2\crho)}
\simeq \nabla_{w_0},
\end{align*}
where the third arrow comes from \eqref{socle of standard, semiinf}, and the last arrow
comes from \eqref{other lambda inv}.

The fact that these maps for various $\clambda$ are compatible with the maps in the
inductive system follows from \lemref{commutation for other st}. The fact that the resulting
map $$\on{Conv}\Bigl(({}^{w_0}\tCM^1)\{2\crho\}\Bigr)\to \nabla_{w_0}$$
factors through $\ConvHecke\Bigl(({}^{w_0}\tCM^1)\{2\crho\}\Bigr)\to \nabla_{w_0}$
follows from \eqref{comm sph2}.


\begin{thebibliography}{199}

\bibitem[ABG]{ABG}
S.~Arkhipov, R.~Bezrukavnikov, V.~Ginzburg, {\em Quantum groups, the loop Grassmannian,
and the Springer resolution,} J. Amer. Math. Soc. {\bf 17} (2004), 595--678.

\bibitem[AB]{AB} S.~Arkhipov, R.~Bezrukavnikov, {\em Perverse sheaves on affine flags
and Langlands dual group}, math.RT/0201073.

\bibitem[AG]{AG} S.Arkhipov, D.~Gaitsgory, {\em Another realization of the category
of modules over the small quantum group}, Adv. in Math {\bf 173} (2003), 114--143.

\bibitem[BBM]{BBM} R.~Bezrukavnikov, A. ~Braverman, I.~Mirkovi\'c,
{\em Some results about geometric Whittaker model},
Adv. Math. {\bf 186} (2004) 143-152.

\bibitem[BG]{BG} A.~Braverman, D. Gaitsgory, {\em Geometric Eisenstein series},
Inv. Math. {\bf 150} (2002), 287--384.

\bibitem[BG1]{BG1} A.~Braverman, D.~Gaitsgory, {\em Crystals via the affine Grassmannian,}
Duke Math. J. {\bf 107}  (2001), 561--575.

\bibitem[BFGM]{BFGM} A.~Braverman, M.~Finkelberg, D.~Gaitsgory, I.~Mirkovic,
{\em Intersection cohomology of Drinfeld's compactifications},
Selecta Math. (N.S.) {\bf 8} (2002), 381--418.

\bibitem[BK]{BK} A.~Braverman, D.~Kazhdan,
{\em On the Schwartz space of the basic affine space,} Selecta Math. (N.S.) {\bf 5} (1999), 1--28.

\bibitem[FFKM]{FFKM} B~Feigin, M.~Finkelberg, A.~Kuznetsov, I.~Mirkovic,
{\em Semi-infinite flags. II. Local and global intersection cohomology
of quasimaps' spaces.  Differential topology, infinite-dimensional Lie algebras, and applications},
113--148, Amer. Math. Soc. Transl. Ser. 2, 194, Amer. Math. Soc., Providence, RI, 1999.

\bibitem[FF]{FF} B.~Feigin, E.~Frenkel, {\em Affine Kac-Moody algebras and semi-infinite flag
manifolds,} Comm. Math. Phys. {\bf 128} (1990) 161--189.

\bibitem[FM]{FM} M.~Finkelberg, I.~Mirkovic,
{\em Semi-infinite flags. I. Case of global curve $\mathbf P\sp 1$.
Differential topology, infinite-dimensional Lie algebras, and applications,}
81--112, Amer. Math. Soc. Transl. Ser. 2, 194, Amer. Math. Soc., Providence, RI, 1999.

\bibitem[FG]{FG} E.~Frenkel, D.~Gaitsgory, {\em Tamely ramified local Langlands correspondence for affine
Kac-Moody algebras}, forthcoming.

\bibitem[FGV]{FGV} E.~Frenkel, D.~Gaitsgory, K.~Vilonen,
{\em Whittaker patterns in the geometry of moduli spaces of bundles
on curves,} Ann. Math. {\bf 153} (2001), 699-748.

\bibitem[KL]{KL} D.~Kazhdan, G.~Lusztig,
{\em Tensor structures arising from affine Lie algebras. III, IV}, J. Amer. Math. Soc.
{\bf 7} (1994), 383--453, 335--381.

\bibitem[KT]{KT} M.~Kashiwara, T.~Tanisaki, {\em The Kazhdan-Lusztig conjecture for affine algebras
with negative level}, Duke Math. J. {\bf 77} (1995), 383-453.

\bibitem[Lu]{Lu} G. ~Lusztig,
{\em Singularities, character formulas, and a $q$-analog of weight multiplicities,}
Analysis and topology on singular spaces, II, III (Luminy, 1981), 208--229, Astérisque {\bf 101-102}.

\bibitem[Lu1]{Lu1} G.~Lusztig, {\em Introduction to quantum groups,}
Progress in Mathematics, 110. Birkhäuser Boston, Inc., Boston, MA, 1993.

\bibitem[MV]{MV} I.~Mirkovic, K.~Vilonen,
{\em Geometric Langlands duality and representations of algebraic groups over commutative rings},
math.RT/0401222.

\bibitem[Soe]{Soe} W.~Soergel, {\em Character formulas for tilting modules over Kac-Moody algebras,}
Rep. Theory  {\bf 2} (1998), 432--448

\end{thebibliography}
\end{document}